\documentclass[a4paper,12pt]{amsart}

\numberwithin{equation}{section}
\theoremstyle{plain}
\newtheorem{theorem}{\sc \bf Theorem}[section]
\newtheorem{lemma}[theorem]{\sc \bf Lemma}
\newtheorem{corollary}[theorem]{\sc \bf Corollary}
\newtheorem{proposition}[theorem]{\sc \bf Proposition}
\newtheorem{claim}{\sc \bf Claim}
\theoremstyle{definition}

\newtheorem{remark}[theorem]{\sc \bf Remark}

\usepackage{amssymb}
\usepackage{amsfonts}
\usepackage{amsmath}
\usepackage{enumerate}

\usepackage{tikz}

\newcommand{\K}{{\mathbb{K}}}
\newcommand{\R}{{\mathbb{R}}}
\newcommand{\C}{{\mathbb{C}}}
\newcommand{\Z}{{\mathbb{Z}}}

\newcommand{\Si}{\Sigma}
\newcommand{\si}{\sigma}
\newcommand{\Ph}{\varPhi}
\newcommand{\ph}{\varphi}
\newcommand{\phe}{\varphi(\epsilon,\cdot)}
\newcommand{\lam}{\lambda}
\newcommand{\e}{\epsilon}
\newcommand{\p}{{\prime}}

\newcommand{\IR}{\mathcal{I}}

\newcommand{\LR}{\mathcal{L}}
\let\MR=\relax
\newcommand{\MR}{\mathcal{M}}
\newcommand{\NR}{\mathcal{N}}

\newcommand{\PR}{\mathcal{P}}
\newcommand{\RR}{\mathcal{R}}
\newcommand{\SR}{\mathcal{S}}
\newcommand{\XR}{\mathcal{X}}

\newcommand{\ZR}{\mathcal{Z}}

\newcommand{\tLR}{\tilde{\mathcal{L}}}

\newcommand{\ali}[1]{\begin{align*}#1\end{align*}}
\newcommand{\alil}[1]{\begin{align}#1\end{align}}
\newcommand{\ite}[1]{\begin{enumerate}[(1)]#1\end{enumerate}}
\let\prop=\relax
\let\cor=\relax
\newcommand{\prop}[1]{\begin{proposition}#1\end{proposition}}
\newcommand{\cor}[1]{\begin{corollary}#1\end{corollary}}
\newcommand{\cors}{\begin{corollary}}
\newcommand{\core}{\end{corollary}}
\newcommand{\thm}[1]{\begin{theorem}#1\end{theorem}}
\newcommand{\thms}{\begin{theorem}}
\newcommand{\thme}{\end{theorem}}
\newcommand{\lem}[1]{\begin{lemma}#1\end{lemma}}
\newcommand{\lems}{\begin{lemma}}
\newcommand{\leme}{\end{lemma}}
\newcommand{\rem}[1]{\begin{remark}\normalfont #1\end{remark}}
\newcommand{\cla}[1]{\begin{claim}#1\end{claim}}
\newcommand{\ncla}[1]{\setcounter{claim}{0}\begin{claim}#1\end{claim}}

\newcommand{\pros}{\begin{proof}}
\newcommand{\proe}{\end{proof}}

\newcommand{\case}[1]{\begin{cases}#1\end{cases}}
\newcommand{\cd}{\cdot}
\newcommand{\up}{\upsilon}
\newcommand{\om}{\omega}
\newcommand{\ti}[1]{\tilde{#1}}

\newcommand{\di}{\displaystyle}

\newcommand{\diam}{{\mathrm{diam}}}

\newcommand{\sign}{\mathrm{sign}}
\let\sgn=\relax
\newcommand{\sgn}{\mathrm{sgn}}

\newcommand{\qqqqqquad}{\qquad\qquad\qquad\qquad}
\newcommand{\qqqqqqquad}{\qquad\qquad\qquad\qquad\qquad}
\newcommand{\qqqqqqqquad}{\qquad\qquad\qquad\qquad\qquad\qquad}

\newcommand{\qqqqqqqqqquad}{\qquad\qquad\qquad\qquad\qquad\qquad\qquad\qquad}

\newcounter{constants}
\makeatletter
\def\addconst{
\addtocounter{constants}{1}
\def\@currentlabel{\arabic{constants}}
\@currentlabel
}
\makeatother

\newcommand{\adl}[1]{\addconst\label{c:#1}}
\newcommand{\adr}[1]{\ref{c:#1}}

\begin{document}
\author[H. Tanaka]{Haruyoshi Tanaka}
\keywords{asymptotic perturbation theory \and topological pressure  \and Hausdorff dimension \and graph iterated function systems}
\subjclass[2010]{37F35 \and 37D35 \and 47A55 \and 37C45}
\title[Asymptotic solution of Bowen equation]{Asymptotic solution of Bowen equation for perturbed potentials on shift spaces with countable states}
\address{
{\rm Haruyoshi Tanaka}\\
Department of Mathematics and Statistics\\
Wakayama Medical University\\
580, Mikazura, Wakayama-city, Wakayama, 641-0011, Japan
}
\email{htanaka@wakayama-med.ac.jp}
\begin{abstract}
We study the asymptotic solution of the equation of the pressure function $s\mapsto P(s\ph(\e,\cd)+\psi(\e,\cd))$ for perturbed potentials $\ph(\e,\cd)$ and $\psi(\e,\cd)$ defined on the shift space with countable state space. In our main result, we give a sufficient condition for the solution $s=s(\epsilon)$ of $P(s\ph(\e,\cd)+\psi(\e,\cd))=0$ to have the $n$-order asymptotic expansion for the small parameter $\e$.
In addition, we also obtain the case where the order of the expansion of the solution $s=s(\epsilon)$ is less than the order of the expansion of the perturbed potentials. Our results can be applied to problems concerning asymptotic behaviors of Hausdorff dimensions obtained from Bowen formula: conformal graph directed Markov systems, an infinite graph directed systems with contractive infinitesimal similitudes mappings, and other concrete examples.
\end{abstract}
\maketitle
\tableofcontents
\section{Main results}\label{sec:intro}
Let $G=(V,E,i(\cd),t(\cd))$ be a directed multigraph endowed with countable vertex set $V$, countable edge set $E$, and two maps $i(\cd)$ and $t(\cd)$ from $E$ to $V$. For each $e\in E$, $i(e)$ is called the initial vertex of $e$ and $t(e)$ called the terminal vertex of $e$. Denoted by $E^{\infty}$ the one-sided shift space $\{\om=\om_{0}\om_{1}\cdots\in \prod_{k=0}^{\infty}E\,:\,t(\om_{n})=i(\om_{n+1}) \text{ for any } n\geq 0\}$ endowed with the shift transformation $\si\,:\,E^{\infty}\to E^{\infty}$ defined as $(\si\om)_{n}=\om_{n+1}$ for any $n\geq 0$. 
For $\theta\in (0,1)$, a metric $d_{\theta}$ on $E^{\infty}$ is given by $d_{\theta}(\om,\up)=\theta^{\inf\{n\geq 0\,:\,\om_{n}\neq \up_{n}\}}$.
The {\it incidence} matrix $A$ of $E^{\infty}$ is defined by $A=(A(ee^\p))_{E\times E}$ with $A(ee^\p)=1$ if $t(e)=i(e^\p)$ and $A(ee^\p)=0$ if $t(e)\neq i(e^\p)$. The matrix $A$ is {\it finitely irreducible} if there exists a finite subset $F$ of $\bigcup_{n=1}^{\infty}E^{n}$ such that for any $e,e^\p\in E$, $ewe^\p$ is a path on the graph $G$ for some $w\in F$.
A function $f\,:\,E^{\infty}\to \K$ is called {\it weakly $d_{\theta}$-Lipschitz continuous} if the number $\sup_{e\in E}\sup_{\om,\up\in [e]\,:\,\om\neq \up}|f(\om)-f(\up)|/d_{\theta}(\om,\up)$ is finite. A function $f\,:\,E^{\infty}\to \K$ is a weakly H\"older continuous function if this is a weakly $d_{\theta}$-Lipschitz continuous for some $\theta\in (0,1)$. Denoted by $\|\cd\|_{\infty}$ the supremum norm defined as $\|f\|_{\infty}=\sup_{\om\in E^{\infty}}|f(\om)|$.
\smallskip
\par
To state our main result, we introduce some conditions for potentials.
Let $n$ be a nonnegative integer. We consider the following conditions (g.1)-(g.5) for function $g(\e,\cd)\,:\,E^{\infty}\to \R$ with small parameter $\e\in (0,1)$:
\ite
{
\item[(g.1)] A function $g(\e,\cd)\,:\,E^{\infty}\to \R$ has the form $g(\e,\cd)=g+g_{1}\e+\cdots+g_{n}\e^{n}+\ti{g}_{n}(\e,\cd)\e^{n}$ for some real-valued weakly H\"older continuous functions $g, g_{1},\dots, g_{n}, \ti{g}_{n}(\e,\cd)$ with $\lim_{\e\to 0}\|\ti{g}_{n}(\e,\cd)\|_{\infty}=0$.
\item[(g.2)] $g(\om)\neq 0$ for each $\om\in E^{\infty}$ and $\|g\|_{\infty}<1$.
\item[(g.3)] $|g(\om)-g(\up)|\leq c_{\adr{g3}}|g(\om)|d_{\theta}(\om,\up)$ for $\om,\up\in E^{\infty}$ with $\om_{0}=\up_{0}$ for some $c_{\adl{g3}}>0,\theta\in (0,1)$.
\item[(g.4)] $|g_{k}(\om)|\leq c_{\adr{g4}}|g(\om)|^{t_{k}}$ and $|g_{k}(\om)-g_{k}(\up)|\leq c_{\adr{g4_2}}|g(\om)|^{t_{k}}d_{\theta}(\om,\up)$ for any $\om,\up\in E^{\infty}$ with $\om_{0}=\up_{0}$ for some constants $c_{\adl{g4}},c_{\adl{g4_2}}>0$ and $t_{k}\in (0,1]$ for $k=1,2,\dots, n$.
\item[(g.5)] $|\ti{g}_{n}(\e,\om)|\leq c_{\adr{g5}}(\e)|g(\om)|^{\ti{t}}$ for any $\om\in E^{\infty}$ for some constants $\ti{t}\in (0,1]$ and $c_{\adl{g5}}(\e)>0$ with $c_{\adr{g5}}(\e)\to 0$.
}
Moreover, we assume that function $\psi(\e,\cd)\,:\,E^{\infty}\to \R$ satisfy the following conditions ($\psi.1$)-($\psi.4$):
\ite
{
\item[($\psi.1$)] A function $\psi(\e,\cd)\,:\,E^{\infty}\to \R$ has the form $\psi(\e,\cd)=\psi+\psi_{1}\e+\cdots+\psi_{n}\e^{n}+\ti{\psi}_{n}(\e,\cd)\e^{n}$ for some real-valued weakly H\"older continuous functions $\psi, \psi_{1},\dots, \psi_{n}, \ti{\psi}_{n}(\e,\cd)$ with $\lim_{\e\to 0}\|\ti{\psi}_{n}(\e,\cd)\|_{\infty}=0$.
\item[($\psi.2$)] $\psi(\om)>0$ for any $\om\in E^{\infty}$.
\item[($\psi.3$)] $|\psi_{k}(\om)|\leq c_{\adr{psi3}}|\psi(\om)|$ and $|\psi_{k}(\om)-\psi_{k}(\up)|\leq c_{\adr{psi3_2}}|\psi(\om)|d_{\theta}(\om,\up)$ for any $\om,\up\in E^{\infty}$ with $\om_{0}=\up_{0}$ and for some $c_{\adl{psi3}},c_{\adl{psi3_2}}>0$ for $k=1,2,\dots, n$.
\item[($\psi.4$)] $|\ti{\psi}_{n}(\e,\om)|\leq c_{\adr{psi4}}(\e)|\psi(\om)|$ for any $\om\in E^{\infty}$ for some $c_{\adl{psi4}}(\e)>0$ with $c_{\adr{psi4}}(\e)\to 0$.
}
Let 
\alil
{
\underline{p}=\inf\{p\geq 0\,:\,P(p\log |g|+\log\psi)<+\infty\},
}
where $P(f)$ means the topological pressure of $f$ which is defined by (\ref{eq:toppres}). Put
\alil
{
p(n)=
\case
{
\underline{p}/\ti{t},&n=0\\
\max\big(\underline{p}+n(1-t_{1}), \underline{p}+n(1-t_{2})/2, \cdots, \underline{p}+n(1-t_{n})/n,\\
\qqqqqqquad\underline{p}/t_{1},\ \underline{p}/t_{2},\ \cdots, \underline{p}/t_{n},\ \underline{p}+1-\ti{t},\ \underline{p}/\ti{t}\big),& n\geq 1.
}\label{eq:t(n)=...} 
}
\par
Now we are in a position to state our main result:
\thms
\label{th:asympsol_Beq_Mgene}
Fix nonnegative integer $n$. Assume that the incidence matrix of $E^{\infty}$ is finitely irreducible and the conditions $(g.1)$-$(g.5)$ and $(\psi.1)$-$(\psi.4)$ are satisfied. Choose any $s(0)\in (p(n),+\infty)$ and any compact neighborhood $I\subset (p(n),+\infty)$ of $s(0)$. Let $p_{0}=P(s(0)\log|g|+\log\psi)$. Then there exist numbers $\e_{0}>0, s_{1},\dots, s_{n} \in \R$ such that the equation
\alil
{
P(s\log|g(\e,\cd)|+\log\psi(\e,\cd))=p_{0}\label{eq:exp(P(telog|ge|+logpsie))=}
}
has a unique solution $s=s(\e)\in I$ for each $0<\e<\e_{0}$, and $s(\e)$ has the asymptotic expansion
\alil
{
s(\e)=s(0)+s_{1}\e+\cdots+s_{n}\e^{n}+\ti{s}_{n}(\e)\e^{n}\label{eq:t(e)=t0+t1e+...}
}
and $|\ti{s}_{n}(\e)|\to 0$ as $\e\to 0$. In particular
\alil
{
\ti{s}_{n}(\e)=
\case
{
\di-\frac{\nu(\e,\tLR_{n,s(\e)}(\e,h))}{\nu(h\log|g|)}+O(\e),& n\geq 1\\
\di-\frac{\nu(\e,\tLR_{0,s(\e)}(\e,h))}{\nu(h\log|g|)}+o(\|\tLR_{0,s(\e)}(\e,h)\|_{\infty}),& n=0,
}\label{eq:tsne=}
}
where $h$ is the Perron eigenfunction of the eigenvalue $e^{p_{0}}$ of the Ruelle operator of $s(0)\log|g|+\log\psi$, $\nu$ is the Perron eigenvector of this dual operator with $\nu(h)=\nu(E^{\infty})=1$, and 
$\nu(\e,\cd)$ is the Perron eigenvector of the dual of the Ruelle operator of $s(\e)\log|g(\e,\cd)|+\log\psi(\e,\cd)$. Moreover, $\tLR_{n,s(\e)}(\e,\cd)$ is an operator as the remainder of the expansion of the Ruelle operator (see Lemma \ref{lem:asymp_Rop_ge^t_Mgene} for detail).
\thme
Note that the coefficients $s_{k}$ and the remainder $\ti{s}_{n}(\e)$ are precisely given in (\ref{eq:tn=...}) and (\ref{eq:ttne=...n>=1}), respectively.
The following follows immediately from this theorem:
\cors
\label{cor:asympsol_Beq_Mgene}
Under the same conditions of the above theorem, assume also that there exists $s(0)>p(n)$ such that $P(s(0)\log|g|+\log \psi)=0$. Then the equation $P(s\log|g(\e,\cd)|+\log\psi(\e,\cd))=0$ for $s\in \R$ has a unique solution $s=s(\e)$ for any small $\e>0$ and $s(\e)$ has an $n$-order asymptotic expansion as well as (\ref{eq:t(e)=t0+t1e+...}).
\core
Next, we give sufficient conditions for a situation that $s(\e)$ does not have $(n+1)$-order asymptotic expansion though it has $n$-order asymptotic behaviour. We introduce the following conditions:
\ite
{
\item[(g.6)] $E$ is infinitely countable.
\item[(g.7)] $\ti{g}_{n}(\e,\cd)\equiv 0$ (i.e. $\ti{t}=1$) and $\psi(\e,\cd)\equiv 1$.
\item[(g.8)] There exist $t_{0}\in [t_{1}, (nt_{1}+1)/(n+1))$ and $c_{\adl{g9}}>0$ such that $p(n)<\underline{p}+(1-t_{0})(n+1)$ and 
$g_{1}(\om)\sign(g(\om))\geq c_{\adr{g9}}|g(\om)|^{t_{0}}$ for any $\om\in E^{\infty}$.
\item[(g.9)] The numbers $t_{2},\dots, t_{n}$ satisfy $t_{0}\leq t_{k}$ for any $k=2,\dots, n$.
}
\prop
{\label{prop:asympsol_Beq_est_rem2}
Assume that the conditions (g.1)-(g.9) with $n\geq 1$ are satisfied. Choose any $s(0)\in (p(n), \underline{p}+(n+1)(1-t_{0}))\setminus\{1,2,\dots, n\}$ and put $p_{0}=P(s(0)\log|g|)$. Then a unique solution $s=s(\e)$ of the equation $P(s\log|g(\e,\cd)|)=p_{0}$ has the form $s(\e)=s(0)+s_{1}\e+\cdots+s_{n}\e^{n}+\ti{s}_{n}(\e)\e^{n}$ with $\lim_{\e\to 0}|\ti{s}_{n}(\e)|/\e= +\infty$.
}
\rem
{
\ite
{
\item[(1)] Theorem \ref{th:asympsol_Beq_Mgene} is proved by developing the proofs which considered the case of finite case in \cite{T2011,T2016}. A difficult point between the finite state case and the infinite state case is that even if the remainder $\ti{g}_{n}(\e,\cd)$ equals zero, the solution $s(\e)$ do not may have $(n+1)$-order expansion in the sense that $|\ti{s}_{n}(\e)|/\e\to +\infty$ as $\e\to 0$ (Proposition \ref{prop:asympsol_Beq_est_rem2}). In particular, this remainder can become any small fractional order (see Theorem \ref{th:ex_5}(2)). This fact suggests that the general analytic perturbation theory cannot be applied to the asymptotic behaviour of the solution $s(\e)$ in the infinite state case. In fact, in the finite state case $\sharp E<+\infty$, if $\ti{g}_{n}(\e,\cd)\equiv 0$ and the conditions (g.1), (g.2), ($\psi$.1) and ($\psi$.2) are satisfied, then the solution $s(\e)$ has asymptotic expansion with any integer order $m$ with $m\geq n$ \cite[Theorem 2.6]{T2011}. By the extension of the theorem of asymptotic solutions for topological pressures in \cite{T2011} (the proof of Theorem \ref{th:asymp_dim}), and by the generalization of asymptotic perturbation theory of bounded linear operators of \cite{T2011} (Theorem \ref{th:asymp_e.vec}), the main theorem is proved.
\item[(2)] Our results can be directly applied to problems concerning asymptotic behaviors of Hausdorff dimensions obtained from Bowen formula (Section \ref{sec:ex}). In particular, the coefficient and the remainder of the solution $s(\e)$ can be numerically calculated (Section \ref{sec:ex_fail(P2)_2}). Note that though the functions $\psi(\e,\cd)$ may be equal to $\psi(\e,\cd)\equiv 1$ in our examples of this paper, we shall treat the case of $\psi(\e,\cd)\not\equiv 1$ to study of a multifractal analysis of a perturbed system in future work. Another future studying is to estimate the dimensions of limit sets of nonconformal graph iterated function systems with infinite state. 
}
}
In the next section \ref{sec:ex}, we will illustrate asymptotic perturbations of dimensions of various examples: conformal graph directed Markov system, infinite graph-directed systems on Banach spaces and Continued fractions. 
Furthermore, we will demonstrate an example of linear countable IFS such that the coefficients and the remainder of the solution are explicit calculated (Section \ref{sec:ex_fail(P2)_2}).
In the section \ref{sec:proof}, we mention the proofs of our all results. In the appendices, we shall introduce some facts necessary for the proof of the main theorem.
In Appendix \ref{sec:thermo}, we recall the notion of thermodynamic formalism and the Ruelle operators acting on suitable function space in the infinite graph. In particular, a version of Ruelle-Perron-Frobenius Theorem of this operators is described (Theorem \ref{th:exGibbs}). We state in Appendix \ref{sec:abstract_asymp} the general theory of asymptotic behaviors of the eigenvalue and the corresponding eigenvector of bounded linear operators. This result is obtained by generalizing the results of \cite[Theorem 2.1]{T2011}.
Finally, we shall give upper bound of intermediate point of the binomial expansion in Appendix \ref{sec:interme_binom} which plays an important role in giving the proof of our results.
\medskip
\\
\noindent
{\it Acknowledgment.}\ 
This study was supported by JSPS KAKENHI Grant Number 20K03636.
\section{Examples}\label{sec:ex}
\subsection{Conformal graph directed Markov systems}\label{sec:CGDMS}
Let $G=(V,E,i(\cd),t(\cd))$ be a directed multigraph which $V$ is finite and $E$ is countable.
In this section, we consider asymptotic behaviours of the Hausdorff dimensions of the limit sets of perturbed graph directed Markov systems introduced by \cite{MU}.
We begin with the definition of this system. Let $D$ be a positive integer, $\beta\in (0,1]$ and $r\in (0,1)$. We introduce a set $(G,(J_{v}),(O_{v}),(T_{e}))$ satisfying the following conditions (i)-(iv):
\ite
{
\item[(i)] For each $v\in V$, $J_{v}$ is a compact and connected subset of $\R^{D}$ satisfying that the interior $\mathrm{int} J_{v}$ of $J_{v}$ is not empty, and $\mathrm{int} J_{v}$ and $\mathrm{int}J_{v^\p}$ are disjoint for $v^\p\in V$ with $v\neq v^\p$.
\item[(ii)] For each $v\in V$, $O_{v}$ is a bounded, open and connected subset of $\R^{D}$ containing $J_{v}$.
\item[(iii)] For each $e\in E$, a function $T_{e}\,:\,O_{t(e)}\to T_{e}(O_{t(e)})\subset O_{i(e)}$ is a $C^{1+\beta}$-conformal diffeomorphism with $T_{e}(\mathrm{int}J_{t(e)})\subset \mathrm{int}J_{i(e)}$ and $\sup_{x\in O_{t(e)}}\|T_{e}^\p(x)\|\leq r$, where $\|T_{e}^\p(x)\|$ means the operator norm of $T_{e}^\p(x)$. Moreover, for any $e,e^\p\in E$ with $e\neq e^\p$ and $i(e^\p)=i(e)$, $T_{e}(\mathrm{int}J_{t(e)})\cap T_{e^\p}(\mathrm{int}J_{t(e^\p)})=\emptyset$, namely the open set condition (OSC) is satisfied.
\item[(vi)] (Bounded distortion) There exists a constant $c_{\adl{Gbd}}>0$ such that for any $e\in E$ and $x,y\in O_{t(e)}$, $|\|T_{e}^\p(x)\|-\|T_{e}^\p(y)\||\leq c_{\adr{Gbd}}\|T_{e}^\p(x)\| |x-y|^{\beta}$, where $|\cd|$ means a norm of any Euclidean space. 
\item[(v)] (Cone condition) If $\sharp E=\infty$ then there exist $\gamma, l>0$ with $\gamma<\pi/2$ such that for any $v\in V$, $x\in J_{v}$, there is $u\in \R^{D}$ with $|u|=1$ so that the set $\{y\in \R^{D}\,:\,0<|y-x|<l \text{ and }(y-x,u)>|y-x|\cos \gamma\}$ is in $\mathrm{int} J_{v}$, where $(y-x,u)$ denotes the inner product of $y-x$ and $u$.
}
Under these conditions (i)-(v), we call the set $(G,(J_{v}), (O_{v}), (T_{e}))$ a graph directed Markov system (GDMS for short).
The Hausdorff dimension of the limit set of this system has been mainly studied by many authors \cite{MU1999,MU,RU,RU2}.
\smallskip
\par
The coding map $\pi\,:\,E^{\infty}\to \R^{D}$ is defined by $\pi\om=\bigcap_{n=0}^{\infty}T_{\om_{0}}\cdots T_{\om_{n}}(J_{t(\om_{n})})$ for $\om \in E^{\infty}$. Put $K=\pi(E^{\infty})$. This set is called the limit set of the GDMS.
We define a function $\ph\,:\,E^{\infty}\to \R$ by
\ali
{
\ph(\om)=\log \|T_{\om_{0}}^\p(\pi\si\om)\|.
}
Put
\ali
{
\underline{s}=\inf\{s\geq 0\,:\,P(s\ph)<+\infty\}.
}
We call the GDMS {\it regular} if $P(s\ph)=0$ for some $s\geq \underline{s}$. 
The GDMS is said to be {\it strongly regular} if $0<P(s\ph)<+\infty$ for some $s\geq \underline{s}$ (see \cite{MU,RU} for the terminology).
It is known that the general Bowen's formula is satisfied:
\thm
{[\cite{RU}]\label{th:MU}
Let $(G,(J_{v}),(O_{v}),(T_{e}))$ be a graph directed Markov system. Assume that $E^{\infty}$ is finitely irreducible. Then $\dim_{H}K=\inf\{t\in \R\,:\,P_{G}(t\ph)\leq 0\}$. In addition to the above condition, we also assume that the potential $\ph$ is regular. Then $s=\dim_{H}K$ if and only if $P_{G}(s\ph)=0$.
}
Now we formulate an asymptotic perturbation of graph directed Markov systems. Fix integers $n\geq 0$, $D\geq 1$ and a number $\beta\in (0,1]$. Consider the following conditions $(G.1)_{n}$ and $(G.2)_{n}$:
\ite
{
\item[$(G.1)_{n}$] The code space $E^{\infty}$ is finitely irreducible. The set $(G,(J_{v}),(O_{v}),(T_{e}))$ is a GDMS on $\R^{D}$ with strongly regular and the limit set $K$ has positive dimension. Moreover, the function $T_{e}$ is of class $C^{1+n+\beta}(O_{t(e)})$ for each $e\in E$.
\item[$(G.2)_{n}$] The set $\{(G,(J_{v}),(O_{v}),(T_{e}(\e,\cd)))\,:\,\e>0\}$ is a GDMS with a small parameter $\e>0$ satisfying the following (i)-(vi):
\ite
{
\item[(i)] For each $e\in E$, the function $T_{e}(\e,\cd)$ has the $n$-asymptotic expansion:
\ali
{
T_{e}(\e,\cd)=T_{e}+T_{e,1}\e+\cdots+T_{e,n}\e^{n}+\ti{T}_{e,n}(\e,\cd)\e^{n}\ \text{ on }J_{t(e)}
}
for some functions $T_{e,k}\in C^{1+n-k+\beta}(O_{t(e)},\R^{D})$ $(k=1,2,\dots, n)$ and $\ti{T}_{e,n}(\e,\cd)\in C^{1+\beta(\e)}(O_{t(e)},\R^{D})$ $(\beta(\e)>0)$ satisfying $\sup_{e\in E}\sup_{x\in J_{t(e)}}|\ti{T}_{e,n}(\e,x)|\to 0$.
\item[(ii)] There exist constants $t(l,k)\in (0,1]$ ($l=0,1,\dots, n$,\ $k=1,\dots, n-l+1$) such that the function $x\mapsto T_{e,l}^{(k)}(x)/\|T_{e}^{\p}(x)\|^{t(l,k)}$ is bounded, $\beta$-H\"older continuous and its H\"older constant is bounded uniformly in $e\in E$.
\item[(iii)] $c_{\adl{Gbd3}}(\e):=\sup_{e\in E}\sup_{x\in J_{t(e)}}(\|\frac{\partial}{\partial x}\ti{T}_{e,n}(\e,x)\|/\|T_{e}^\p(x)\|^{\ti{t}_{0}})\to 0$ as $\e\to 0$ for some $\ti{t}_{0}\in (0,1]$.
\item[(iv)] $\dim_{H}K/D>p(n)$, where $p(n)$ is taken by (\ref{eq:t(n)=...}) with
\alil
{
t_{k}:=&\min\{\frac{1}{D}\sum_{p=1}^{D}t(i_{p},j_{p}+1)\,:\,i:=i_{1}+\cdots+i_{D}\text{ and }j:=j_{1}+\cdots+j_{D}\text{ satisfy}\label{eq:tk=}\\
&\qqqqqquad i=k \text{ and }j=0 \text{ or }0\leq i<k \text{ and }1\leq i+j\leq k\}\nonumber\\
\ti{t}:=&\min\left\{t_{n},\ \ti{t}_{0},\ \frac{\ti{t}_{0}}{D}+\frac{D-1}{D}t(1,1),\dots,\ \frac{\ti{t}_{0}}{D}+\frac{D-1}{D}t(n,1)\right\}\label{eq:tt=}\\
\underline{p}:=&\underline{s}/D.\nonumber
}
}
}
Note that if the edge set $E$ is finite, then the conditions (ii) and (iv) are always satisfied because $\|T_{e}^\p(x)\|$ is uniformly bounded away from zero, and $p(n)$ becomes zero by taking $t(l,k)\equiv 1$. Moreover, $c_{\adr{Gbd3}}(\e)$ in (iii) can be taken as $\sup_{e\in E}\sup_{x\in J_{t(e)}}\|\frac{\partial }{\partial x}\ti{T}_{e,n}(\e,x)\|$ when $E$ is finite. Let $K(\e)$ be the limit set of the perturbed GDMS $(G,(J_{v}),(O_{v}),(T_{e}(\e,\cd)))$.
Then we obtain the following result:
\thm
{\label{th:asymp_dim}
Assume that the conditions $(G.1)_{n}$ and $(G.2)_{n}$ are satisfied with fixed integer $n\geq 0$. Then the perturbed GDMS $(G,(J_{v}),(O_{v}),(T_{e}(\e,\cd)))$ is strongly regular for any small $\e>0$, and 
there exist $s_{1},\dots, s_{n}\in \R$ such that the Hausdorff dimension $\dim_{H}K(\e)$ of the limit set $K(\e)$ of the perturbed system has the form $\dim_{H}K(\e)=\dim_{H}K+s_{1}\e+\cdots+s_{n}\e^{n}+o(\e^{n})$ as $\e \to 0$.
}
\rem
{
Roy and Urba\'nski \cite{RU} considered continuous perturbation of infinitely conformal iterated function systems given as a special GDMS. They also studied analytic perturbation of GDMS with $D\geq 3$ in \cite{RU2}. We investigated an asymptotic perturbation of GDMS with finite graph in \cite{T2016}. Theorem \ref{th:asymp_dim} is an infinite graph version of this previous result in \cite{T2016}
}
\subsection{Infinite graph-directed systems on Banach spaces}
In this section, we consider an abstract graph-directed systems introduced in \cite{Priyadarshi,NPL}. We assume the following conditions (I.1)-(I.4).
\ite
{
\item[(I.1)] Let $G=(V,E,i(\cd),t(\cd))$ be a directed multigraph satisfying $V$ is finite, $E$ is infinitely countable and $E^{\infty}$ has finitely irreducible incidence matrix.
\item[(I.2)] For $v\in V$, $(J_{v},d_{v})$ is a compact perfect metric space with $J_{v}\subset J$, where $J$ is a metric space with metric $d$ and where  $J_{v}$ is {\it perfect} if for any $x\in J$ there exists a sequence $x_{k}\in J_{v}$ with $x_{k}\neq x$ for $k\geq 1$ such that $d_{v}(x_{k},x)\to 0$ as $k\to \infty$.
\item[(I.3)] For any $e\in E$, $T_{e}\,:\,J_{t(e)}\to J_{i(e)}$ is a map satisfying the following:
\ite
{
\item[(i)] Each $T_{e}$ is a Lipschitz map satisfying  $\sup_{e\in E}\mathrm{Lip}(T_{e})=:r<1$, where $\mathrm{Lip}(T_{e}):=\sup_{x,y\in J_{t(e)}\,:\,x\neq y}d_{i(e)}(T_{e}(x),T_{e}(y))/d_{t(e)}(x,y)$.
\item[(ii)] Each $T_{e}$ is an {\it infinitesimal similitude}, i.e. for each $x\in J_{t(e)}$, for any sequences $(x_{k})$ and $(y_{k})$ on $J_{t(e)}$ with $x_{k}\neq y_{k}$ for each $k\geq 1$ and $x_{k}\to x$ and $y_{k}\to y$, the limit
\ali
{
\lim_{k\to \infty}\frac{d_{i(e)}(T_{e}(x_{k}),T_{e}(y_{k}))}{d_{t(e)}(x_{k},y_{k})}=:DT_{e}(x)
}
exists in $\R$ and is independent of the particular sequences $(x_{k})$ and $(y_{k})$.
\item[(iii)] $DT_{e}(x)>0$ For any $e\in E$ and $x\in J$.
\item[(iv)] There exist constants $c>0$ and $\beta>0$ such that for any $e\in E$ and $x,y\in J_{t(e)}$, $|DT_{e}(x)-DT_{e}(y)|\leq c|DT_{e}(x)| d_{t(e)}(x,y)^{\beta}$.
\item[(v)] There exist $t>0$ such that for any $v$, there is $x\in J_{v}$ such that $\sum_{e\in E\,:\,t(e)=v}DT_{e}(x)^{t}<+\infty$.
\item[(vi)] For any $\eta>0$ there exists $c(\eta)\geq 1$ such that for each $e\in E$ and $x,y\in J_{t(e)}$ with $0<d_{t(e)}(x,y)<\eta$,
\ali
{
c(\eta)^{-1}DT_{e}(x)\leq \frac{d(T_{e}(x),T_{e}(y))}{d_{t(e)}(x,y)}\leq c(\eta) DT_{e}(x)
}
and $\lim_{\eta\to +0}c(\eta)=1$.
}
}
We write $E$ as $E=\{e_{1},e_{2},\dots\}$. Under the condition (I.1), there exists $N\geq 1$ such that $E_{k}:=\{e_{1},e_{2},\dots, e_{k}\}$ satisfies $i(E)\cup t(E)=i(E_{k})\cup t(E_{k})$ and $(E_{k})^{\infty}$ has the finite irreducible incidence matrix. Then there exist nonempty unique subsets $K_{k,v}\subset J_{v}$ $(v\in V)$ such that $K_{k,v}=\bigcup_{e\in E_{k}\,:\,i(e)=v}T_{e}(K_{k,t(e)})$ for any $v\in V$. We call the set $K_{k}:=\bigcup_{v\in V}K_{k,v}$ the limit set of $(G,(T_{e})_{e\in E_{k}})$. We also assume the following.
\ite
{
\item[(I.4)] For each $k\geq N$, the limit set $K_{k}$ of the finite graph-directed system $(G,(T_{e})_{e\in E_{k}})$ satisfies that $T_{e}(K_{k})\cap T_{e^\p}(K_{k})=\emptyset$ for each $e\neq e^\p$, and $T_{e}|_{K_{k}}\,:\,K_{k}\to J$ is one to one for $e\in E_{k}$.
}
Such a system is firstly introduced by \cite{NPL} and developed by \cite{Priyadarshi}.

The coding map $\pi\,:\,E^{\infty}\to J$ is well defined as $\pi\om=\bigcap_{k=0}^{\infty}T_{\om_{0}\cdots \om_{k}}(J)$ and the limit set of the system $(G,(T_{e}))$ is given by $K=\pi(E^{\infty})$. Put 
\ali
{
\ph(\om)=\log DT_{\om_{0}}(\pi \si\om)
}
and $\underline{p}=\inf\{p\geq 0\,:\,P(p\ph)<\infty\}$. 
We see that this function $\ph$ is $d_{\theta}$-Lipschitz continuous function on $E^{\infty}$ with $\theta:=r^{\beta}$ by the condition (I.3)(iv). A version of Bowen's formula follows from this theorem:
\thm
{[\cite{NPL}]\label{thm:NPL_dim2}
Assume that $(G,(T_{e}))$ is a system satisfying the conditions (I.1)-(I.4) and $K$ is its limit set. 
Then we have $\dim_{H}K=\inf\{s>0\,:\,P(s\ph)<0\}$. In particular, if the system $(G,(T_{e}))$ is strongly regular, i.e. there is $p>\underline{p}$ so that $P(s\ph)=0$, then $P(s\ph)=0$ iff $s=\dim_{H}K$.
}
\pros
Let
\ali
{
\underline{p}^{*}=\inf\{t>0\,:\,\sum_{v\in V}\sum_{e\in E\,:\,i(e)=v}DT_{e}(x_{v})^{t}<+\infty \text{ for all }(x_{v})\in \prod_{v\in V}J_{v}\}.
}
For $p>\underline{p}^{*}$, we define an operator $\MR_{p}\,:\,\prod_{v\in V}C(J_{v},\R)\to \prod_{v\in V}C(J_{v},\R)$ by
\alil
{
(\MR_{p} f((x_{v})))_{v}=\sum_{e\in E\,:\,t(e)=v}DT_{e}(x)^{p}f_{i(e)}(T_{e}x)\label{eq:GBoweneq_P}
}
for $f=(f_{v})\in \prod_{v\in V}C(J_{v},\R)$, $v\in V$ and $x\in J_{v}$, where $C(J_{v},\R)$ denotes the set of all $\R$-valued continuous functions on $J_{v}$.
Then by virtue of \cite[Theorem 5.5]{Priyadarshi}, the equation
\ali
{
\dim_{H}K=\inf\{p\geq 0\,:\,\text{the spectral radius of }\MR_{p}\text{ is less than }1\}
}
holds. First we check the equality $\underline{p}^{*}=\underline{p}$ for $p>\underline{p}$. For $f\in \prod_{v\in V}C(J_{v},\R)$ and $\om\in E^{\infty}$, we define a map $\Pi(f,\om)$ by $\Pi(f,\om)=f_{i(\om_{0})}(\pi\om)$. We notice that for each $p>\underline{p}^{*}$, $f\in \prod_{v\in V}C(J_{v})$ and $k\geq 1$,
\alil
{
\LR_{p\ph}^{k}(\Pi(f,\cd))(\om)=&\sum_{w=w_{1}\cdots w_{k}\in E^{k}\,:\,t(w_{k})=i(\om_{0})}(\prod_{i=1}^{k}DT_{w_{i}}(\pi \si^{i-1}w\cd\om))^{p}f_{i(w_{1})}(\pi w\cd\om)\nonumber\\
=&\sum_{w=w_{1}\cdots w_{k}\in E^{k}\,:\,t(w_{k})=i(\om_{0})}(\prod_{i=1}^{k}DT_{w_{i}}(T_{w_{i}\cdots w_{k}}\pi\om))^{p}f_{i(w_{1})}(T_{w}\pi \om)
=(\MR_{p}^{k}f)_{i(\om_{0})}(\pi\om)\label{eq:LRt^k=MRt^k}
}
with the fact $\pi(e\cd\om)=T_{e}\pi(\om)$ for $e\in E$. Since $\|\LR_{p\ph}1\|_{\infty}\leq \|\MR_{p}1\|:=\max_{v\in V}\|(\MR_{p}1)_{v}\|<+\infty$, we see $\underline{p}\leq \underline{p}^{*}$. Conversely, when $p>\underline{p}$, $(\MR_{p}1)_{v}(x)$ is finite for some $x\in J_{v}$ for any $v\in V$, ans then $(\MR_{s}1)_{v}(x)$ is finite for any $s\geq p$, $x\in J_{v}$ and $v\in V$ (\cite[Lemma 3.3]{Priyadarshi}). Therefore $\underline{p}^{*}\leq \underline{p}$ and Thus $\underline{p}^{*}=\underline{p}$.

Finally we will show that the spectral radius $r(\LR_{p\ph})$ coincides with the spectral radius $r(\MR_{p})$ for each $p>\underline{p}$. By the equation (\ref{eq:LRt^k=MRt^k}), $\|\LR_{p\ph}^{k}\|_{\infty}=\|\LR_{p\ph}^{k}1\|_{\infty}\leq \|\MR_{p}^{k}\|$ holds and thus $r(\LR_{p\ph})\leq r(\MR_{p})$. To see the converse, \cite[Lemma 3.7]{Priyadarshi} yields the eigenfunction $u=(u_{v})_{v\in V}\in \prod_{v\in V}C(J_{v},\R)$ of the eigenvalue $r(\MR_{p})$ such that $\MR_{p}u=r(\MR_{p})u$ and $u_{i}>0$. This also implies $\LR_{p\ph}\Pi(u,\om)=(\MR_{p}u)_{i(\om_{0})}(\pi\om)=r(\MR_{p})u_{i(\om_{0})}(\pi\om)=r(\MR_{p})\Pi(u,\om)$. Therefore, $r(\MR_{p})$ is an eigenvalue of $\LR_{p\ph}$ and thus $r(\LR_{t\ph})\geq r(\MR_{p})$. We obtain $r(\LR_{p\ph})= r(\MR_{p})$. Consequently, we can apply the operator $\LR_{p\ph}$ to Ruelle-Perron-Frobenius type theorem and in particular, we obtain $P(p\ph)=\log r(\LR_{p\ph})=r(\MR_{p})$.
Hence the proof is complete by the equation (\ref{eq:GBoweneq_P}).
\proe
In order to formulate an asymptotic perturbation of NPL systems, we consider the following conditions:
\ite
{
\item[$(J.1)_{n}$] A pair $(J,(T_{e}))$ is a strongly regular NPL system satisfying that $J$ is a compact subset of a Banach space $(X,\|\cd\|)$. Moreover, there exists an open connected subset $O$ of $X$ containing $J$ such that each $T_{e}$ is extended to a map of $C^{n+1}(O,X)$ and $DT_{e}$ is extended to a map of $C^{n+\beta}(O,\R)$.
\item[$(J.2)_{n}$] A pair $(J,(T_{e}(\e,\cd)))$ is a NPL system with a small parameter $\e>0$ satisfying that the following (i)-(iv) hold:
\ite
{
\item[(i)] For each $e\in E$, $T_{e}(\e,\cd)$ and $DT_{e}(\e,\cd)$ have $n$-order asymptotic expansions:
\ali
{
T_{e}(\e,\cd)=&\textstyle T_{e}+\sum_{k=1}^{n}T_{e,k}\e^{k}+\ti{T}_{e,n}(\e,\cd)\e^{n}\ \text{ on }J\\
DT_{e}(\e,\cd)=&\textstyle DT_{e}+\sum_{k=1}^{n}S_{e,k}\e^{k}+\ti{S}_{e,n}(\e,\cd)\e^{n}\ \text{ on }J
}
for some $T_{e,k}\in C^{n-k+1}(O)$, $S_{e,k}\in C^{n-k+\beta}(O)$ $(1\leq k\leq n)$, $\ti{T}_{e,n}(\e,\cd)\in C^{1}(O)$ and $\ti{S}_{e,n}(\e,\cd)\in C^{\beta(\e)}(O)$ with $\beta(\e)>0$.
\item[(ii)] There exist $t(k,m)\in (0,1]$ $(0\leq k\leq n, 0\leq m\leq n-k)$ such that the function $x\mapsto S_{e,k}^{(m)}/DT_{e}(x)^{t(k,m)}$ is bounded, $\beta$-H\"older continuous and its H\"older constant is bounded uniformly in $e\in E$.
\item[(iii)] $\sup_{e\in E}\sup_{x\in J}\|\ti{T}_{e,n}(\e,x)\|\to 0$ and $\sup_{e\in E}\sup_{x\in J}|\ti{S}_{e,n}(\e,x)|/DT_{e}(x)^{\ti{t}_{0}}\to 0$ as $\e\to 0$ for some $\ti{t}_{0}\in (0,1]$.
\item[(iv)] $\dim_{H}K>p(n)$, where $p(n)$ is taken as (\ref{eq:t(n)=...}) with
\ali
{
t_{k}:=&\min\{t(i,j)\,:\,i=k \text{ and }j=0 \text{ or }0\leq i<k \text{ and }0\leq i+j\leq k\}\quad (1\leq k\leq n)\\
\ti{t}:=&\min\{t_{n},t(n,0),t(n-1,1),\dots, t(0,n),\ti{t}_{0}\}.
}
}
}
Then we have the following:
\thms\label{th:AP_NPL_MOREgene}
Fix $n\geq 0$. Assume that the conditions $(J.1)_{n}$ and $(J.2)_{n}$ are satisfied. Then the Hausdorff dimension of the limit set $K(\e)$ of the system $(J,(T_{e}(\e,\cd)))$ has an /$n$-asymptotic expansion for $\e$.
\thme
\subsection{Continued fractions}
Consider a graph $G$ with singleton vertex set $V=\{v\}$ and with infinite edge set $E\subset \{k\in \Z\,:\,k\geq 2\}$. Put $J_{v}=[0,1]$ and $O_{v}=(-\eta,1+\eta)$ for a fixed small number $\eta>0$. Fix $a\in \R$ with $a\neq 0$. We define a perturbed map of continued fractions
\ali
{
T_{e}(\e,x)=\frac{1}{e+x+a \e}
}
for $e\in E$, $\e>0$ and $x\in \R$. Consider a GDMS $(G,J_{v},O_{v},(T_{e}(\e,\cd)))$ such that this unperturbed GDMS is strongly regular. Note that such a system exists \cite{MU1999} by choosing edge set $E$. The function $T_{e}(\e,x)$ has the expansion
\ali
{
T_{e}(\e,x)=\frac{1}{e+x}-\frac{a}{(e+x)^{2}}\e+\cdots+(-1)^{n}\frac{a^{n}}{(e+x)^{n+1}}\e^{n}+\ti{T}_{n}(\e,x)\e^{n}
}
with $T_{e,k}(x)=(-a)^{k}(e+x)^{-k-1}$ and $\ti{T}_{n}(\e,x)=\e(-1)^{n+1}a^{n+1}/((e+x)^{n+1}(e+x+a\e))$. It is not hard to check that the conditions $(G.1)_{n}$ and $(G.2)_{n}$ are fulfilled putting $\underline{p}=0$, $t(l,k)=\ti{t}=1$, and therefore $p(n)=0$. Thus we have the $n$-order asymptotic expansion of the limit set of the GDMS by Theorem \ref{th:asymp_dim}.
\subsection{Linear countable IFS (1)}\label{sec:ex_fail(P2)_2}
In this section, we will give the coefficient and the estimate of the remainder for a concrete GDMS. Let $a>1$ and $E=\{1,2,\dots\}$. We take an infinite graph $G=(\{v\},E,i(\cd),t(\cd))$ with $i(e)=t(e)=v$ for $e\in E$, $J_{v}=[0,1]$ and $O_{v}=(-\eta,1+\eta)$ for a small $\eta>0$. For $e\in E$ and $\e\geq 0$, we define a function $T_{e}(\e,\cd)$ by
\ali
{
T_{e}(\e,x)=\Big(\frac{1}{5^{e}}+\frac{1}{a^{e}}\e\Big)x+b(e).
}
Here we choose $b(e)$ so that the set $(G,(J_{v}),(O_{v}),T_{e}(\e,\cd))$ satisfies the open set condition for any small $\e>0$. It is not hard to check that the condition $(G.1)_{n}$-(i)(iii) are valid with $T_{e}(x)=x/5^{e}+b(e)$, $T_{e,1}(x)=x/a^{e}$, $T_{e,k}\equiv 0$ $(k\geq 2)$ and $\ti{T}_{e,n}(\e,\cd)\equiv 0$. To see $(G.2)_{n}$-(ii)(iv), we remark that $|T_{1,e}^\p(x)|/|T_{e}^\p(x)|^{t}=(5^{t}/a)^{e}$ is bounded uniformly in $e$ if and only if $t\leq \log a/\log 5$. Therefore we put $t(1,1)=\min(\log a/\log 5,1)$ and otherwise $t(l,k)=1$ for any $(l,k)\neq (1,1)$ when $n\geq 1$. Let $\ph(\om)=\log(1/5^{\om_{0}})$.
Moreover, $\underline{s}=\inf\{s\geq 0\,:\,P(s\ph)<+\infty\}$ is equal to $0$. Thus $p(n)=n(1-\min(\log a/\log 5,1))$ for any $n\geq 0$.
We see $P(s(0)\ph)=0$ if and only if $\sum_{e\in E}(1/5^{e})^{s(0)}=1$ if and only if $s(0)=\dim_{H}K=\log 2/\log 5$ by the Bowen's formula. Then we obtain the following:
\thms
\label{th:ex_5}
Assume the above conditions for $T_{e}(\e,\cd)$.
\ite
{
\item[(1)] If $a\geq 5$ then the Hausdorff dimension $s(\e)=\dim_{H}K(\e)$ of the limit set of this GDMS has $n$-asymptotic expansion $s(\e)=\log 2/\log 5+s_{1}\e+\cdots+s_{n}\e^{n}+\ti{s}_{n}(\e)\e^{n}$ with $\ti{s}_{n}(\e)\to 0$ for any $n\geq 0$. Each coefficient $s_{k}$ $(k=1,2,\dots, n)$ is decided as
\alil
{
s_{k}=\frac{1}{2\log 5}\sum_{0\leq v\leq k,0\leq q\leq k-v\,:\,\atop{(v,q)\neq (0,1)}}\sum_{j=0}^{\min(v,q)}s_{q,k-v}\frac{a_{v,j,s(0)}}{(q-j)!}(-\log 5)^{q-j}\sum_{e=1}^{\infty}e^{q-j}\left(\frac{5^{v}}{2 a^{v}}\right)^{e},\label{eq:sk=}
}
where constants $s_{q,k-v}$ and $a_{v,j,s(0)}$ are defined by (\ref{eq:tki=}) and (\ref{eq:aljs=}), respectively.
\item[(2)] If $1<a<5$ then take the largest integer $k\geq 0$ satisfying $a\leq 5/2^{1/(k+1)}$. In this case, $s(\e)$ has the form
\ali
{
s(\e)=
\case
{
s(0)+s_{1}\e+\cdots +s_{k}\e^{k}+\hat{s}(\e)\e^{\frac{\log 2}{\log(5/a)}},&a< 5/2^{1/(k+1)}\\
s(0)+s_{1}\e+\cdots +s_{k}\e^{k}+\hat{s}(\e)\e^{k+1}\log\e,&a=5/2^{1/(k+1)}\\
}
}
with $|\hat{s}(\e)|\asymp 1$ as $\e\to 0$, where $s_{1},\dots,s_{k}$ are given by (\ref{eq:sk=}) and where $b(\e)\asymp c(\e)$ means $d^{-1}c(\e)\leq b(\e)\leq dc(\e)$ for any small $\e>0$ for some constant $d\geq 1$. Note that $k<\log 2/\log(5/a)\leq k+1$ is satisfied,
}
\thme
In particular, the numbers $s_{1}$ and $s_{2}$ are given by 
\ali
{
s_{1}=&\frac{\log 2}{(\log 5)^{2}}\frac{5}{4a-10}\\
s_{2}=&\frac{25\log 2}{(\log 5)^{3}}\left(\frac{1}{2(2a-5)^{2}}-\frac{a\log 2}{(2a-5)(4a^2-5)^2}+\frac{\log(2/5)}{8a^2-100}\right).
}
\subsection{Linear countable IFS (2)}\label{sec:ex_CIFS2}
Under the same notation $G,V,E,J_{v},O_{v}$ in Section \ref{th:ex_5}, we define a concrete function $T_{e}(\e,\cd)$ by
\ali
{
T_{e}(\e,x)=\Big(\frac{1}{5^e}+\frac{1}{4^e}\e+\frac{1}{3^e}\e^{2}\Big)x+b(e),
}
where $b(e)$ is suitable chosen so that the OSC is satisfied. By virtue of Theorem \ref{th:asympsol_Beq_Mgene}, we see the following:
\prop
{\label{prop:IIFS_2-asymp}
Under the above function, the Hausdorff dimension of the limit set $K(\e)$ of $(G,(J_{v}),(O_{v}),(T_{e}(\e,\cd))$ has at least $2$-order asymptotic expansion for $\e$.
}
\section{Proofs}\label{sec:proof}
In this section, the theorems and propositions in Section \ref{sec:intro} and Section \ref{sec:ex} are all proved. Recall the notation given in Section \ref{sec:ex}. For later convenience we introduce some function spaces. Let $\K$ be a numerical space or a complex space. Denoted by $C_{b}(E^{\infty},\K)$ the set of all $\K$-valued continuous functions $f$ on $E^{\infty}$ with $\|f\|_{\infty}<+\infty$, by $F_{\theta}(E^{\infty},\K)$ the set of all $\K$-valued weakly $d_{\theta}$-Lipschitz continuous on $E^{\infty}$, and by
 $F_{\theta,b}(E^{\infty},\K)$ the set of all $f\in F_{\theta}(E^{\infty},\K)$ with $\|f\|_{\theta}:=\|f\|_{\infty}+[f]_{\theta}<+\infty$, where we put $[f]_{\theta}:=\sup_{e\in E}\sup_{\om,\up\in [e]\,:\,\om\neq \up}|f(\om)-f(\up)|/d_{\theta}(\om,\up)$. If $\K$ is equal to $\C$, then we may drop the notation `$\K$' from these function spaces.
\subsection{Proof of Theorem \ref{th:asympsol_Beq_Mgene}}
To show our main result, we need to prove some auxiliary propositions. We begin with the following short proposition:
\prop
{\label{prop:nk-Lip_ex}
Let $(U,d)$ be a bounded metric space and $(E,\|\cd\|)$ a Banach algebra. Assume that functions $f_{k}, h_{k}$ $(0\leq k\leq n)$ from $U$ to $E$ satisfy $\|f_{k}(x)\|\leq c_{\adr{nk-L_fk1}}(k)\|h_{k}(x)\|$ for any $x\in U$ and $\|f_{k}(x)-f_{k}(y)\|\leq c_{\adr{nk-L_fk2}}(k)\|h_{k}(x)\|d(x,y)$ for any $x,y\in U$ for some constants $c_{\adl{nk-L_fk1}}(k), c_{\adl{nk-L_fk2}}(k)>0$. Then
\ali
{
{\textstyle\|\prod_{k=0}^{n}f_{k}(x)-\prod_{k=0}^{n}f_{k}(y)\|\leq c_{\adr{nk-L_fk3}}\prod_{i=0}^{n}\|h_{i}(x)\|d(x,y)}
}
with $c_{\adl{nk-L_fk3}}=\sum_{i=0}^{n}(\prod_{j=0}^{i-1}c_{\adr{nk-L_fk1}}(j))c_{\adr{nk-L_fk2}}(i)(\prod_{j=i+1}^{n}(c_{\adr{nk-L_fk1}}(j)+c_{\adr{nk-L_fk2}}(j)\diam U))$.
}
\pros
By the assumption, we note
\ali
{
\|f_{k}(y)\|\leq \|f_{k}(x)\|+\|f_{k}(x)-f_{k}(y)\|\leq (c_{\adr{nk-L_fk1}}(k)+c_{\adr{nk-L_fk2}}(k)\diam U)\|h_{k}(x)\|.
}
Thus we have
\ali
{
\|\prod_{k=0}^{n}f_{k}(x)-\prod_{k=0}^{n}f_{k}(y)\|\leq& \sum_{i=1}^{n}\left(\prod_{k=0}^{i-1}\|f_{k}(x)\|\right)\|f_{i}(x)-f_{i}(y)\|\left(\prod_{k=i+1}^{n}\|f_{k}(y)\|\right)\\
\leq&c_{\adr{nk-L_fk3}}\prod_{s=1}^{n}\|h_{s}(x)\|d(x,y).
}
\proe
For a convenience later, for $p\in \R$ with $p>0$ and $k,l\in \Z$ with $0\leq k\leq n$ and $0\leq l\leq k$, we define a function $G^{p}_{l,k}\,:\,E^{\infty}\to \R$ by
\ali
{
G^{p}_{l,k}(\om)=
\case
{
\di\frac{|g(\om)|^{p}}{g(\om)^{l}}\sum_{j_{1},\dots,j_{k}\geq 0\,:\atop{j_{1}+\cdots+j_{k}=l\atop{j_{1}+2j_{2}+\dots+k j_{k}=k}}}\frac{l!}{j_{1}!\cdots j_{k}!}g_{1}(\om)^{j_{1}}\cdots g_{k}(\om)^{j_{k}},&k\geq 1\\
|g(\om)|^{p},&k=0.\\
}
}
This function will be used in the expansion of $|g(\e,\cd)|^{p}$ (see (\ref{eq:gkp=})), To estimate this function, we will assert the following lemma:
\lems
\label{lem:prop_g^q}
Assume that $G=(V,E,i(\cd),t(\cd))$ be a directed multigraph. Assume that $g\,:\,E^{\infty}\to \R$ satisfies (g.2) and (g.3). Then for any integer $l\geq 0$, $q\in \R$,  and $\om,\up\in E^{\infty}$ with $\om_{0}=\up_{0}$, we have $||g(\om)|^{q}/g(\om)^{l}-|g(\up)|^{q}/g(\up)^{l}|\leq c_{\adr{cg^q}}|g(\om)|^{q-l}d_{\theta}(\om,\up)$ by putting $c_{\adl{cg^q}}=c_{\adr{cg^q}}(q,l)=(c_{\adr{g3}}(1+c_{\adr{g3}}\theta)^{|q|}+\theta^{-1}l)(1+c_{\adr{g3}}\theta)^{l}$.
\leme
\pros
Let $\om,\up\in E^{\infty}$ with $\om_{0}=\up_{0}$.
The condition (g.3) yields $|g(\om)^{-1}-g(\up)^{-1}|\leq |g(\om)^{-1}|c_{\adr{g3}}d_{\theta}(\om,\up)$. On the other hand, 
by virtue of the assumption (g.3) again, we have
\ali
{
(1+c_{\adr{g3}}\theta)^{-1}|g(\om)|\leq |g(\up)|\leq (1+c_{\adr{g3}}\theta)|g(\om)|.
}
Therefore Mean Valued Theorem and the condition (g.3) imply 
\ali
{
||g(\om)|^{q}-|g(\up)|^{q}|=&|\alpha|g(\om)|+(1-\alpha)|g(\up)||^{q-1}||g(\om)|-|g(\up)||\\
\leq&|\alpha|g(\om)|+(1-\alpha)|g(\up)||^{q}\frac{||g(\om)|-|g(\up)||}{\min(|g(\om)|,|g(\up)|)}\\
\leq&|\alpha|g(\om)|+(1-\alpha)|g(\up)||^{q}c_{\adr{g3}}d_{\theta}(\om,\up)\\
\leq&\case
{
\max(|g(\om)|,|g(\up)|)^{q}c_{\adr{g3}}d_{\theta}(\om,\up),&q\geq 0\\
\min(|g(\om)|,|g(\up)|)^{q}c_{\adr{g3}}d_{\theta}(\om,\up),&q< 0\\
}\\
\leq&(1+c_{\adr{g3}}\theta)^{|q|}c_{\adr{g3}}|g(\om)|^{q}d_{\theta}(\om,\up)
}
for some $\alpha\in [0,1]$. Choose any $e\in E$.
Proposition \ref{prop:nk-Lip_ex} regarding as $U=[e]$, $f_{0}=h_{0}:=|g|^{q}$, $f_{1}=\cdots=f_{l}=h_{1}=\cdots=h_{l}:=g^{-1}$, $c_{\adr{nk-L_fk1}}(\cd):= 1$, $c_{\adr{nk-L_fk2}}(0):=(1+c_{\adr{g3}}\theta)^{|q|}c_{\adr{g3}}$ and $c_{\adr{nk-L_fk2}}(1)=\cdots=c_{\adr{nk-L_fk2}}(l):=c_{\adr{g3}}$ implies that the assertion holds for the constant $c_{\adr{cg^q}}=\sum_{i=0}^{l}c_{\adr{nk-L_fk2}}(i)\prod_{j=i+1}^{l}(1+c_{\adr{nk-L_fk2}}(j)\theta)$. Hence the proof is complete.
\proe
\lems
\label{lem:bdd_Gtlk}
Assume that the incidence matrix of $E^{\infty}$ is finitely irreducible and the conditions $(g.1)$-$(g.5)$ with fixed nonnegative integer $n$ are satisfied. Then for any $p> p(n)$, $0\leq k\leq n$ and $0\leq l\leq k$, the function $G^{p}_{l,k}$ is a weakly H\"older continuous function. In particular,
\alil
{
|G^{p}_{l,k}(\om)|\leq& c_{\adr{cGtlk}}|g(\om)|^{p-\frac{k}{n}p(n)+\frac{k}{n}\underline{p}}\label{eq:Gtlk<=}\\
|G^{p}_{l,k}(\om)-G^{p}_{l,k}(\up)|\leq& c_{\adr{cGtlk2}}|g(\om)|^{p-\frac{k}{n}p(n)+\frac{k}{n}\underline{p}}d_{\theta}(\om,\up)\label{eq:|Gtlkom-Gtlkup|<=}
}
for any $\om,\up\in E^{\infty}$ with $\om_{0}=\up_{0}$ for some constants $c_{\adl{cGtlk}}=c_{\adr{cGtlk}}(k,l),c_{\adl{cGtlk2}}=c_{\adr{cGtlk2}}(k,l)>0$.
\leme
\pros
In the case $k=0$, the inequality (\ref{eq:|Gtlkom-Gtlkup|<=}) follows from Lemma \ref{lem:prop_g^q}.

In the case $k\geq 1$, let $\om,\up\in E^{\infty}$ with $\om_{0}=\up_{0}$. 
By the assumption (g.4), we have
$|g_{i}(\om)^{j}|\leq c_{\adr{g4}}^{j}|g(\om)|^{jt_{i}}$ and $|g_{i}(\up)^{j}|\leq (c_{\adr{g4}}+c_{\adr{g4_2}}\theta)^{j}|g(\om)|^{jt_{i}}$. Then it follows from Proposition \ref{prop:nk-Lip_ex} that for any positive integer $j$, $|g_{i}(\om)^{j}-g_{i}(\up)^{j}|\leq j (c_{\adr{g4}}+c_{\adr{g4_2}}\theta)^{j-1}c_{\adr{g4_2}}|g(\om)|^{jt_{i}}d_{\theta}(\om,\up)$. Thus we obtain
\ali
{
|g|^{p-l}|g_{1}^{j_{1}}\cdots g_{k}^{j_{k}}|\leq&c_{\adr{g4}}^{j_{1}+\cdots+j_{k}}|g|^{p-l+j_{1}t_{1}+\cdots+j_{k}t_{k}}\\
\Big|\frac{|g(\om)|^{p}}{g(\om)^{l}}\prod_{i=1}^{k}(g_{i}(\om))^{j_{i}}-\frac{|g(\up)|^{p}}{g(\up)^{l}}\prod_{i=1}^{k}(g_{i}(\up))^{j_{i}}\Big|\leq&
c_{\adr{cgkt2}}|g(\om)|^{p-l+j_{1}t_{1}+\cdots+j_{k}t_{k}}d_{\theta}(\om,\up)
}
from Proposition \ref{prop:nk-Lip_ex} again and Lemma \ref{lem:prop_g^q}, where
$c_{\adl{cgkt2}}=c_{\adr{cgkt2}}(p,l):=c_{\adr{cg^q}}(p,l)(c_{\adr{g4}}+c_{\adr{g4_2}}\theta)^{l}+l(c_{\adr{g4}}+c_{\adr{g4_2}}\theta)^{l-1}$ using $j_{1}+j_{2}+\cdots+j_{k}=l$.
When we put $p_{k}=\underline{p}+(n/k)(1-t_{k})$, then $p_{k}\leq p(n)$ and therefore $t_{k}\geq 1-(k/n)(p(n)-\underline{p})$ are satisfied. We also note that 
\ali
{
p-l+j_{1}t_{1}+\cdots+j_{k}t_{k}\geq &p-l+\sum_{i=1}^{k}(j_{i}-\frac{ij_{i}}{n}(p(n)-\underline{p}))
= p-\frac{k}{n}p(n)+\frac{k}{n}\underline{p}
}
is satisfied using $j_{1}+\cdots+j_{k}=l$ and $j_{1}+2j_{2}+\cdots+kj_{k}=k$. Hence the lemma is complete by putting $c_{\adr{cGtlk}}=\binom{k-1}{l-1}c_{\adr{g4}}^{l}$ and $c_{\adr{cGtlk2}}=\binom{k-1}{l-1}c_{\adr{cgkt2}}$.
\proe
\rem
{\label{rem:p-k/npn+k/np_>=}
When we take $p, \eta>0$ so that $p>p(n)+\eta$, 
we see $p-\frac{k}{n}p(n)+\frac{k}{n}\underline{p}\geq (1-\frac{k}{n})p+\frac{k}{n}\underline{p}+\frac{k}{n}\eta\geq \underline{p}+\eta$.
Namely, the series $\sum_{e\,:\,t(e)=i(\om_{0})}G^{p}_{l,k}(e\cd\om)f(e\cd\om)$ converges for bounded function $f$.
}
\lems
\label{lem:asymp_|gew|^t_ex2}
Assume that the incidence matrix of $E^{\infty}$ is finitely irreducible and the conditions $(g.1)$-$(g.5)$ and $(\psi.1)$-$(\psi.4)$ with fixed nonnegative integer $n$ are satisfied. Then for any $p> p(n)$, there exist weakly H\"older continuous functions $g_{0,p},g_{1,p},\dots, g_{n,p}, \ti{g}_{n,p}(\e,\cd)$ and positive constants $c_{\adl{tP1}}=c_{\adr{tP1}}(k),c_{\adl{tP2}}=c_{\adr{tP2}}(k)$ and $\e^\p$ such that for $0<\e<\e^\p$
\alil
{
|g(\e,\cd)|^{p}=g_{0,p}+g_{1,p}\e+\cdots+g_{n,p}\e^{n}+\ti{g}_{n,p}(\e,\cd)\e^{n}\label{eq:geom=_2}
}
with 
\alil
{
g_{k,p}(\om)=\sum_{l=0}^{k}\binom{p}{l}G^{p}_{l,k}(\om)\label{eq:gkp=}
}
satisfying that for any $k=0,1,\dots, n$ and for any $\om\in E^{\infty}$
\alil
{
|g_{k,p}(\om)|\leq &c_{\adr{tP1}}|g(\om)|^{p-\frac{k}{n}p(n)+\frac{k}{n}\underline{p}}\label{eq:gjt<=_2}\\
|g_{k,p}(\om)-g_{k,p}(\up)|\leq &c_{\adr{tP2}}|g(\om)|^{p-\frac{k}{n}p(n)+\frac{k}{n}\underline{p}}d_{\theta}(\om,\up) \text{ with }\om_{0}=\up_{0},\label{eq:|gjtom-gjtup|<=_2}
}
where $\binom{p}{l}$ is the binomial coefficient. Moreover for any nonempty compact subset $I$ of the interval $(p(n),+\infty)$, there exist constants $\ti{\eta}>0$ and $c_{\adl{tP3}}(\e)>0$ with $c_{\adr{tP3}}(\e)\to 0$ such that
\alil
{
&\sup_{p\in I}|\ti{g}_{n,p}(\e,\om)|\leq c_{\adr{tP3}}(\e)|g(\om)|^{\underline{p}+\ti{\eta}}\quad \text{ for any small }\e>0 \text{ and for any }\om\in E^{\infty}.\label{eq:supt_tgnt(ew)<=_2}
}
\leme
\pros
Lemma \ref{lem:bdd_Gtlk} guarantees the inequalities (\ref{eq:gjt<=_2}) and (\ref{eq:|gjtom-gjtup|<=_2}) by putting $c_{\adr{tP1}}=\sum_{l=0}^{k}|\binom{p}{l}|c_{\adr{cGtlk}}(k,l)$ and $c_{\adr{tP2}}=\sum_{l=0}^{k}|\binom{p}{l}|c_{\adr{cGtlk2}}(k,l)$.
It is remain to show that (\ref{eq:supt_tgnt(ew)<=_2}) holds. Fix $\om\in E^{\infty}$. For the convenience sake, we omit `$\om$' from the notation, i.e. we write $g=g(\om), g_{k}=g_{k}(\om), \ti{g}_{n}(\e)=\ti{g}_{n}(\e,\om)$, $g_{k,p}=g_{k,p}(\om)$ and $\ti{g}_{n,p}(\e)=\ti{g}_{n,p}(\e,\om)$.
Put $x(\e)=\sum_{k=1}^{n}g_{k}\e^{k}$ and $g(\e)=g+x(\e)+\ti{g}_{n}(\e)\e^{n}$.
We also assume we take $c_{\adr{g5}}(\e)$ in the assumption (g.5) so small with $c_{\adr{g5}}(\e)^{\ti{t}}<1/2$.
We take numbers $\underline{c}=\inf I$, $\overline{c}=\sup I$ and $\eta\in (0,\underline{c}-p(n))$. Then for any $p\in I$, we see $\overline{c}\geq p\geq \underline{c}>p(n)+\eta>\underline{p}+\eta>\underline{p}$. We start with the following claim:
\ncla{
The assertion (\ref{eq:supt_tgnt(ew)<=_2}) holds for $n=0$.
}
Indeed, we will consider the two cases: $|g|\leq c_{\adr{g5}}(\e)$ and $|g|> c_{\adr{g5}}(\e)$.
\smallskip
\\
Case: $|g|\leq c_{\adr{g5}}(\e)$. In this case, we have $|g|^{p-\underline{p}-\eta}\leq c_{\adr{g5}}(\e)^{p-\underline{p}-\eta}$ and so $|g|^{p}\leq c_{\adr{g5}}(\e)^{p-\underline{p}-\eta}|g|^{\underline{p}+\eta}$. Therefore
\ali
{
||g(\e)|^{p}-|g|^{p}|\leq& (c_{\adr{1or2pp-1_n}}-1)|g|^{p}+c_{\adr{1or2pp-1_n}}|\ti{g}_{0}(\e)|^{p}\\
\leq&(c_{\adr{1or2pp-1_n}}-1)c_{\adr{g5}}(\e)^{p-\underline{p}-\eta}|g|^{\underline{p}+\eta}+c_{\adr{1or2pp-1_n}}c_{\adr{g5}}(\e)^{p}|g|^{p\ti{t}}\qquad (\because (g.5))\\
\leq& ((c_{\adr{1or2pp-1_n}}-1)c_{\adr{g5}}(\e)^{p-\underline{p}-\eta}+c_{\adr{1or2pp-1_n}}c_{\adr{g5}}(\e)^{p})|g|^{\underline{p}+\ti{t}\eta}\quad (\because p>p(n)+\eta\geq \underline{p}/\ti{t}+\eta)\\
\leq& (c_{\adr{1or2pp-1_n}}-1+c_{\adr{1or2pp-1_n}})c_{\adr{g5}}(\e)^{\underline{c}-\underline{p}-\eta}|g|^{\underline{p}+\ti{t}\eta}=c_{\adr{g5}}(\e)^{\underline{c}-\underline{p}-\eta}|g|^{\underline{p}+\ti{t}\eta}
}
follows with $c_{\adl{1or2pp-1_n}}=\max(1,2^{p-1})$, where the first inequality holds by the basic facts of inequalities (e.g. \cite[Corollary 8.1.4.]{Kuczma3}).
\smallskip
\\
Case: $|g|> c_{\adr{g5}}(\e)$. It follows from $|g|^{1-\ti{t}}\geq c_{\adr{g5}}(\e)^{1-\ti{t}}$ that
$|\ti{g}_{0}(\e)|\leq c_{\adr{g5}}(\e)|g|^{\ti{t}}\leq c_{\adr{g5}}(\e)^{\ti{t}}|g|$.
Therefore, since $\e$ satisfies $c_{\adr{g5}}(\e)^{\ti{t}}<1/2$, the inequality
\ali
{
\frac{1}{2}|g|\leq |g+\alpha \ti{g}_{0}(\e)|\leq \frac{3}{2}|g|
}
is fulfilled for any $\alpha\in [-1,1]$. In particular, from $\sign(g)(g+\ti{g}_{0}(\e))=|g|+\sign(g)\ti{g}_{0}(\e)\geq \frac{1}{2}|g|>0$, we have $\sign(g)=\sign(g(\e))$ and therefore the equation $|g+\ti{g}_{0}(\e)|=|g|+\sign(g)\ti{g}_{0}(\e)$. By using Mean Valued Theorem for the function $a\mapsto (|g|+a)^{p}$, we get the following estimate:
\ali
{
||g(\e)|^{p}-|g|^{p}|=&p(|g|+\alpha \sign(g)\ti{g}_{0}(\e))^{p-1}|\sign(g)\ti{g}_{0}(\e)|\\
\leq &
\case
{
p(1/2)^{p-1} |g|^{p-1}|\ti{g}_{0}(\e)|,& p<1\\
p(3/2)^{p-1}|g|^{p-1}|\ti{g}_{0}(\e)|,& p\geq 1\\
}\\
\leq &\max(1, \overline{c}(3/2)^{\overline{c}-1})|g|^{p-1}|\ti{g}_{0}(\e)|\\
\leq &\max(1, \overline{c}(3/2)^{\overline{c}-1})c_{\adr{g5}}(\e)^{\ti{t}}|g|^{p}\qquad (\because c_{\adr{g5}}(\e)|g|^{\ti{t}}\leq c_{\adr{g5}}(\e)^{\ti{t}}|g|)\\
\leq &\max(1, \overline{c}(3/2)^{\overline{c}-1})c_{\adr{g5}}(\e)^{\ti{t}}|g|^{\overline{p}+\eta}\quad (\because p>\overline{p}+\eta)
}
by using the fact that $p(1/2)^{p-1}$ $(p<1)$ is bounded for $1$ and $[1,\infty)\ni p\mapsto p(3/2)^{p-1}$ is an increasing function.
Thus, the remainder of the form $|g(\e)|^{p}=|g|^{p}+\ti{g}_{0,p}(\e)$ is estimated as (\ref{eq:supt_tgnt(ew)<=_2}) by putting 
$c_{\adr{tP3}}(\e)=\max(1, \overline{c}(3/2)^{\overline{c}-1})c_{\adr{g5}}(\e)^{\min(\ti{t}, \underline{c}-\underline{p}-\eta)}$ and $\eta^\p=\ti{t}\eta$.
\cla
{\label{cla:toyasymp_|f|^p-p_-eta/2<=e^n}
The assertion (\ref{eq:supt_tgnt(ew)<=_2}) holds in the cases $n\geq 1$ and $|g|^{p-\underline{p}-\eta/2}\leq \e^{n}$, i.e. $|\ti{g}_{n,p}(\e)|\leq c_{\adr{ctFnex_n}}(\e)|g|^{\underline{p}+c_{\adr{ceta8_n}}\eta}$ for some positive numbers $c_{\adl{ctFnex_n}}(\e)$ and $c_{\adl{ceta8_n}}$ with $c_{\adr{ctFnex_n}}(\e)\to 0$.
}
Note the form $\ti{g}_{n,p}(\e)=(-|g|^{p}-\sum_{k=1}^{n}g_{k,p}\e^{k}+|g(\e)|^{p})/\e^{n}$. First we estimate the functions $|g|^{p}$ in this expression. Put $c_{\adl{1/c-p_-eta/2}}=1/(\overline{c}-\underline{p}-\eta/2)$. Then we see the inequality 
$|g|\leq \e^{\frac{n}{p-\underline{p}-\eta/2}}\leq \e^{c_{\adr{1/c-p_-eta/2}}n}$.
Since $|g|^{p-\eta/4}\leq |g|^{\underline{p}+\eta/4}\e^{n}$ and $|g|^{\eta/4}\leq \e^{c_{\adr{1/c-p_-eta/2}}n\eta/4}=:c_{\adl{ceta1_n}}(\e)$ are satisfied, we get the inequality
\alil
{
|g|^{p}=|g|^{\eta/4}|g|^{p-\eta/4}\leq c_{\adr{ceta1_n}}(\e)|g|^{\underline{p}+\eta/4}\e^{n}.\label{eq:|g|^p<=}
}
Next we will consider $g_{k,p}$ in the expression of $\ti{g}_{n,p}(\e)$. We have
\alil
{
|g_{k,p}|\leq& c_{\adr{tP1}}(k)|g|^{\frac{n-k}{n}p+\frac{k}{n}\underline{p}+\frac{k}{n}\eta}\qquad (\because (\ref{eq:gjt<=_2}) \text{ and }|g|^{\frac{k}{n}p}\leq |g|^{\frac{k}{n}p(n)+\frac{k}{n}\eta})\nonumber\\
\leq&c_{\adr{tP1}}(k)|g|^{\frac{n-k}{n}\frac{\eta}{2}+\underline{p}+\frac{k}{n}\eta}\e^{n-k}\quad (\because |g|^{p}\leq |g|^{\underline{p}+\eta/2}\e^{n})\nonumber\\
=&c_{\adr{tP1}}(k)|g|^{(1+\frac{k}{n})\frac{\eta}{2}+\underline{p}}\e^{n-k}\nonumber\\
\leq&c_{\adr{tP1}}(k)c_{\adr{ceta7_n}}(\e)|g|^{\underline{p}+\frac{\eta}{2}}\e^{n-k}\label{eq:|gkp|<=}
}
with $c_{\adl{ceta7_n}}(\e)=\e^{c_{\adr{1/c-p_-eta/2}}\eta/2}$. Finally we consider the inequality
\ali
{
|g(\e)|^{p}&\leq c_{\adr{max13p-1_n}}(|g|^{p}+|x(\e)|^{p}+|\ti{g}_{n}(\e)|^{p})
}
with $c_{\adl{max13p-1_n}}=\max(1,3^{\overline{c}-1})$. Let $p_{k}=\underline{p}+(n/k)(1-t_{k})$. It follows from $p>p(n)+\eta\geq p_{k}+\eta$ and $t_{k}=1-(k/n)(p_{k}-\underline{p})$ that
\ali
{
|x(\e)|^{p}\leq&c_{\adr{max1np-1_n}}c_{\adr{g4}}^{p}\sum_{k=1}^{n}|g|^{t_{k}p}\e^{kp}\\
\leq&
c_{\adr{max1np-1_n}}c_{\adr{g4}}^{p}\sum_{k=1}^{n}\case
{
|g|^{\underline{p}+t_{k}\eta}\e^{n},&(kp\geq n)\quad (\because p>\underline{p}/t_{k}+\eta)\\
|g|^{t_{k}p-\frac{1}{n}(p-\underline{p}-\eta/2)(n-kp)}\e^{n},&(kp<n)\quad (\because |g|^{\frac{1}{n}(p-\underline{p}-\eta/2)(n-kp)}\leq \e^{n-kp})\\
}\\
\leq&
c_{\adr{max1np-1_n}}c_{\adr{g4}}^{p}\sum_{k=1}^{n}\case
{
|g|^{\underline{p}+t_{k}\eta}\e^{n},&(kp\geq n)\\
|g|^{p-\frac{kp}{n}(p-\underline{p}-\eta)-(p-\underline{p}-\frac{\eta}{2})(1-\frac{kp}{n})}\e^{n},&(kp< n)\quad (\because t_{k}>1-(k/n)(p-\underline{p}-\eta))\\
}\\
=&
c_{\adr{max1np-1_n}}c_{\adr{g4}}^{p}\sum_{k=1}^{n}\case
{
|g|^{\underline{p}+t_{k}\eta}\e^{n},&(kp\geq n)\\
|g|^{\underline{p}+\frac{\eta}{2}+\frac{\eta}{2}\frac{kp}{n}}\e^{n},&(kp< n)\\
}\\
\leq&
c_{\adr{ceta3_n}}(\e)|g|^{\underline{p}+c_{\adr{ceta2_n}}\eta},
}
where $c_{\adl{max1np-1_n}}=\max\{1,n^{\overline{c}-1}\}$, $c_{\adl{ceta3_n}}(\e)=c_{\adr{max1np-1_n}}c_{\adr{g4}}^{p}n\e^{n c_{\adr{ceta2_n}}c_{\adr{1/c-p_-eta/2}}}$ and $c_{\adl{ceta2_n}}=\min(t_{1},\dots,t_{n},1/2)/2$. Furthermore, we obtain
\ali
{
|\ti{g}_{n}(\e)|^{p}\leq c_{\adr{g5}}(\e)^{p}|g|^{\ti{t}p}\leq c_{\adr{g5}}(\e)^{\underline{c}}|g|^{\underline{p}+\ti{t}\eta}
}
from $p>\underline{p}/\ti{t}+\eta$.
Thus $|g(\e)|^{p}$ is estimated by
\alil
{
|g(\e)|^{p}\leq c_{\adr{fptsn3_n}}(\e)|g|^{\underline{p}+c_{\adr{ceta4_n}}\eta}\e^{n}\label{eq:|ge|^p<=}
}
with
$c_{\adl{fptsn3_n}}(\e)=c_{\adr{max13p-1_n}}\max\{c_{\adr{ceta1_n}}(\e), c_{\adr{ceta3_n}}(\e), c_{\adr{g5}}(\e)^{\underline{c}}\}$ and 
$c_{\adl{ceta4_n}}=\min\{1/4,c_{\adr{ceta2_n}},\ti{t}\}$.
Consequently, the inequalities (\ref{eq:|g|^p<=}), (\ref{eq:|gkp|<=}) and (\ref{eq:|ge|^p<=}) imply
\ali
{
|\ti{g}_{n,p}(\e)|&\leq c_{\adr{ceta1_n}}(\e)|g|^{\underline{p}+\eta/4}+\sum_{k=1}^{n}c_{\adr{tP1}}(k)c_{\adr{ceta7_n}}(\e)|g|^{\underline{p}+\eta/2}+c_{\adr{fptsn3_n}}(\e)|g|^{\underline{p}+c_{\adr{ceta4_n}}\eta}\\
&\leq c_{\adr{ctFnex_n}}(\e)|g|^{\underline{p}+c_{\adr{ceta4_n}}\eta}
}
with $c_{\adr{ctFnex_n}}(\e)=c_{\adr{ceta1_n}}(\e)+\sum_{k=1}^{n}c_{\adr{tP1}}(k)c_{\adr{ceta7_n}}(\e)+c_{\adr{fptsn3_n}}(\e)$.
Hence the assertion is valid.
\cla
{\label{cla:|xe|+|tgnee^n|<=}
If $n\geq 1$, $|g|^{p-\underline{p}-\eta/2}>\e^{n}$ and $c_{\adr{ceta9_n}}(\e)< 1/2$ are satisfied, then the inequality
\alil
{
|x(\e)|+|\ti{g}_{n}(\e)\e^{n}|\leq c_{\adr{ceta9_n}}(\e)|g|\label{eq:|xex|+|tfexe^n|<=|fx|/2_2}
}
holds with $c_{\adl{ceta9_n}}(\e)=nc_{\adr{g4}}\e^{c_{\adr{ceta10_n}}}+c_{\adr{g5}}(\e)$ and $c_{\adl{ceta10_n}}= (\eta/2)/(\overline{c}-\underline{p}-\eta/2)>0$. In this case, we have $\sign(g(\e))=\sign(g+x(\e))=\sign(g)$.
}
Indeed, note that $c_{\adr{ceta10_n}}=1-(\overline{c}-\underline{p}-\eta)/(\overline{c}-\underline{p}-\eta/2)$ is less than $1$. The assumption implies
\ali
{ \e^{k-c_{\adr{ceta10_n}}}<|g|^{\frac{k-c_{\adr{ceta10_n}}}{n}(p-\underline{p}-\eta/2)}\leq |g|^{\frac{k-c_{\adr{ceta10_n}}}{n}(\underline{c}-\underline{p}-\eta/2)}=|g|^{\frac{k}{n}(\underline{c}-\underline{p}-\eta/2)-\frac{1}{n}(\eta/2)}.
}
We see
\ali
{
|x(\e)|+|\ti{g}_{n}(\e)\e^{n}|&\leq \sum_{k=1}^{n}c_{\adr{g4}}|g|^{t_{k}}\e^{k}+c_{\adr{g5}}(\e)|g|^{\ti{t}}\e^{n}\\
&\leq \sum_{k=1}^{n}c_{\adr{g4}}|g|^{t_{k}+\frac{k}{n}(\underline{c}-\underline{p}-\frac{\eta}{2})-\frac{\eta}{2n}}\e^{c_{\adr{ceta10_n}}}+c_{\adr{g5}}(\e)|g|^{\ti{t}+(p-\underline{p})-\frac{\eta}{2}}.
}
We now show $t_{k}+\frac{k}{n}(\underline{c}-\underline{p}-\frac{\eta}{2})-\frac{\eta}{2n}> 1$ and $\ti{t}+(p-\underline{p})-\frac{\eta}{2}> 1$. Since the number $p_{k}=\underline{p}+(n/k)(1-t_{k})$ satisfies $p_{k}+\eta\leq p(n)+\eta<\underline{c}$ by the choices of $p(n)$ and $\eta$, we obtain
\ali
{
t_{k}+\frac{k}{n}(\underline{c}-\underline{p}-\frac{\eta}{2})-\frac{\eta}{2n}=&1-\frac{k}{n}(p_{k}-\underline{p})+\frac{k}{n}(\underline{c}-\underline{p}-\frac{\eta}{2})-\frac{\eta}{2n}\\
=&1+\frac{k}{n}(\underline{c}-p_{k})-\frac{k\eta}{2n}-\frac{\eta}{2n}\\
=&1+\frac{k}{n}(\underline{c}-p_{k}-\eta)+\frac{k-1}{2n}\eta> 1
}
and
\ali
{
\ti{t}+(p-\underline{p})-\frac{\eta}{2}&=1-(\ti{p}-\underline{p})+p-\underline{p}-\frac{\eta}{2}=1+p-\ti{p}-\frac{\eta}{2}>1
}
with $\ti{p}:=\underline{p}+1-\ti{t}$ and with $\ti{p}+\eta<p$. Thus we have the estimate (\ref{eq:|xex|+|tfexe^n|<=|fx|/2_2}).

Finally, we prove the last assertion in the claim. Choose any $\e>0$ so that $c_{\adr{ceta9_n}}(\e)< 1/2$. We have
\ali
{
\frac{g}{g(\e)}=\frac{|g|}{|g|+\sign(g)x(\e)+\sign(g)\ti{g}_{n}(\e)}\geq \frac{|g|}{|g|+|x(\e)|+|\ti{g}_{n}(\e)|}\geq \frac{2}{3}>0.
}
This means that $\sign(g)=\sign(g(\e))$. Similarity,
\alil
{
\frac{g}{g+x(\e)}=\frac{|g|}{|g|+\sign(g)x(\e)}\geq \frac{|g|}{|g|+|x(\e)|}\geq \frac{2}{3}>0.\label{eq:g/g+xe=}
}
Hence we see $\sign(g)=\sign(g+x(\e))$.
\cla
{\label{cla:|ge|<=|g+sumgk|^p+Y}
If $n\geq 1$, $0<p<1$, $|g|^{p-\underline{p}-\eta/2}>\e^{n}$ and $c_{\adr{ceta9_n}}(\e)< 1/2$ are satisfied, then
\alil
{
|g+\sum_{k=1}^{n}g_{k}\e^{k}+\ti{g}_{n}(\e)\e^{n}|^{p}=|g+\sum_{k=1}^{n}g_{k}\e^{k}|^{p}+Y(\e)\label{eq:|f+x+tf|^p=|f+x|^p+r_2}
}
and $|Y(\e)|\leq c_{\adr{ceta11_n}}(\e)|g|^{p-1+\ti{t}}\e^{n}$ with $c_{\adl{ceta11_n}}(\e)\to 0$ are valid.
}
Note the form $|g+x(\e)+\ti{g}_{n}(\e)\e^{n}|=|g+x(\e)|+\sign(g)\ti{g}_{n}(\e)\e^{n}$ from the above claim for any small $\e>0$ so that $c_{\adr{ceta9_n}}(\e)<1/2$. By virtue of Mean Valued Theorem, we have
\ali
{
&|g+x(\e)+\ti{g}_{n}(\e)\e^{n}|^{p}\\
=&|g+x(\e)|^{p}+p(|g+x(\e)|+\alpha\sign(g)\ti{g}_{n}(\e)\e^{n})^{p-1}\sign(g)\ti{g}_{n}(\e)\e^{n}
}
for some $\alpha\in [0,1]$. Let $Y(\e)=p(|g+x(\e)|+\alpha\sign(g)\ti{g}_{n}(\e)\e^{n})^{p-1}\sign(g)\ti{g}_{n}(\e)\e^{n}$. We have the estimate
\ali
{
|Y(\e)|\leq& p(|g|-|x(\e)|-|\ti{g}_{n}(\e)\e^{n}|)^{p-1}|\ti{g}_{n}(\e)|\e^{n}\qquad (\because 0<p<1)\\
\leq&2^{1-p}p|g|^{p-1}c_{\adr{g5}}(\e)|g|^{\ti{t}}\e^{n}\\
\leq&2^{1-p}pc_{\adr{g5}}(\e)|g|^{p-1+1+\underline{p}-p+\eta}\e^{n}\quad (\because 1+\underline{p}-p+\eta<\ti{t})\\
\leq&c_{\adr{ceta11_n}}(\e)|g|^{\underline{p}+\eta}\e^{n}.
}
Thus we see the claim by putting $c_{\adr{ceta11_n}}(\e)=2c_{\adr{g5}}(\e)$.
\cla
{\label{cla:|g+...|^p=|g+...|^p+Z}
If $n\geq 1$, $p\geq 1$, $|g|^{p-\underline{p}-\eta/2}>\e^{n}$ and $c_{\adr{ceta9_n}}(\e)< 1/2$ are satisfied, then the equation
\alil
{
|g+\sum_{k=1}^{n}g_{k}\e^{k}+\ti{g}_{n}(\e)\e^{n}|^{p}=|g+\sum_{k=1}^{n}g_{k}\e^{k}|^{p}+Z(\e)\label{eq:|f+x+tf|^p=|f+x|^p+r_ex}
}
and $|Z(\e)|\leq c_{\adr{ceta12_n}}(\e)|g|^{\underline{p}+\ti{t}\eta}\e^{n}$ with $c_{\adl{ceta12_n}}(\e)\to 0$ hold.
}
Indeed, we notice $|g(\e)|=|g+x(\e)|+\sign(g)\ti{g}_{n}(\e)\e^{n}$.
Denoted by $m$ the largest integer that is not greater than $p$. Put $t=p-m$ and $Y(\e)=|g(\e)|^{t}-|g+x(\e)|^{t}$. Then we have
\ali
{
|g(\e)|^{p}=&|g(\e)|^{m}|g(\e)|^{t}\\
=&(|g+x(\e)|^{m}+\sum_{l=0}^{m-1}\binom{m}{l}|g+x(\e)|^{l}(\sign(g)\ti{g}_{n}(\e)\e^{n})^{m-l})(|g+x(\e)|^{t}+Y(\e))\\
=&|g+x(\e)|^{p}+\sum_{l=0}^{m-1}\binom{m}{l}|g+x(\e)|^{l}(\sign(g)\ti{g}_{n}(\e)\e^{n})^{m-l})(|g+x(\e)|^{t}+Y(\e))\\
&+|g+x(\e)|^{m}Y(\e)\\
=&|g+x(\e)|^{p}+Z(\e).
}
We have
\ali
{
|g+x(\e)|^{l+t}|\ti{g}_{n}(\e)|^{m-l}\e^{(m-l)n}\leq&\left(\frac{3}{2}\right)^{l+t}|g|^{l+t}c_{\adr{g5}}(\e)^{m-l}|g|^{\ti{t}(m-l)}\e^{(m-l)n}\qquad (\because (\ref{eq:g/g+xe=}) \text{ and } (g.5))\\
\leq&\e^{(m-l-1)n}\left(\frac{3}{2}\right)^{p}c_{\adr{g5}}(\e)|g|^{l+t\ti{t}+\ti{t}(m-l)}\e^{n}\\
=&\e^{(m-l-1)n}c_{\adr{fptsn*c_n}}(\e)|g|^{l+\ti{t}(p-l)}\e^{n}\\
\leq&\e^{(m-l-1)n}c_{\adr{fptsn*c_n}}(\e)|g|^{\underline{p}+\ti{t}\eta}\e^{n}
}
by putting $c_{\adl{fptsn*c_n}}(\e)=(3/2)^{\overline{c}}c_{\adr{g5}}(\e)$.
Here the last inequality follows from $\ti{t}>(\underline{p}-l+\eta\ti{t})/(p-l)$ for $l=0,1,\dots, m-1$ by using the fact $\ti{t}>(\underline{p}+\eta\ti{t})/p$. Furthermore,
\ali
{
|g+x(\e)|^{l}|\ti{g}_{n}(\e)|^{m-l}\e^{(m-l)n}|Y(\e)|\leq& \left(\frac{3}{2}\right)^{\overline{c}}|g|^{l}c_{\adr{g5}}(\e)^{m-l}|g|^{\ti{t}(m-l)}\e^{(m-l)n}c_{\adr{g5}}(\e)^{t}|g|^{\ti{t}t}\e^{n t}\\
\leq&c_{\adr{fptsn*c_n}}(\e)|g|^{\ti{t}(m-l)+l+\ti{t}t}\e^{(m-l)n+nt}\\
\leq&c_{\adr{fptsn*c_n}}(\e)|g|^{\ti{t}(p-l)+l}\e^{(p-l)n}\\
\leq&\e^{(m-l-1)n}c_{\adr{fptsn*c_n}}(\e)|g|^{\underline{p}+\ti{t}\eta}\e^{n},
}
where the first inequality uses the estimate $|Y(\e)|\leq |\ti{g}_{n}(\e)|^{t}\leq c_{\adr{g5}}(\e)^{t}|g|^{\ti{t}t}\e^{n t}$. 
On the other hand, the remainder $Y(\e)$ has also the estimate $|Y(\e)|\leq c_{\adr{ceta11_n}}(\e)|g|^{t-1+\ti{t}}\e^{n}$ by the claim \ref{cla:|ge|<=|g+sumgk|^p+Y}. Therefore
\ali
{
|g+x(\e)|^{m}|Y(\e)|\leq& \left(\frac{3}{2}\right)^{m}c_{\adr{ceta11_n}}(\e)|g|^{m+t-1+\ti{t}}\e^{n}\\
\leq& 2c_{\adr{fptsn*c_n}}(\e)|g|^{\underline{p}+\eta}\e^{n}.\qquad (\because \ti{t}>1-p+\underline{p}+\eta)
}
Consequently, we obtain
\ali
{
|Z(\e)|\leq& 2\sum_{l=0}^{m-1}\binom{m}{l}\e^{(m-l-1)n}c_{\adr{fptsn*c_n}}(\e)|g|^{\underline{p}+\ti{t}\eta}\e^{n}+2c_{\adr{fptsn*c_n}}(\e)|g|^{\underline{p}+\eta}\e^{n}\\
=&2\sum_{l=0}^{m-1}\frac{m}{m-l}\binom{m-1}{l}\e^{(m-l-1)n}c_{\adr{fptsn*c_n}}(\e)|g|^{\underline{p}\e^{n}+\ti{t}\eta}+2c_{\adr{fptsn*c_n}}(\e)|g|^{\underline{p}+\eta}\e^{n}\\
\leq&(2m(1+\e)^{m-1}+2)c_{\adr{fptsn*c_n}}(\e)|g|^{\underline{p}+\ti{t}\eta}\e^{n}.
}
Hence the assertion is valid by putting $c_{\adr{ceta12_n}}(\e)=(2m(1+\e)^{m-1}+2)c_{\adr{fptsn*c_n}}(\e)$.
\cla
{\label{cla:toyasymp_|f|^p-p_-eta/2>e^n}
If $n\geq 1$, $|g|^{p-\underline{p}-\eta/2}> \e^{n}$ and $c_{\adr{ceta9_n}}(\e)< 1/2$ are satisfied, then $|\ti{g}_{n,p}(\e)|\leq c_{\adr{ce4_n}}(\e)|g|^{\underline{p}+\ti{t}\eta}$ holds for some constant $c_{\adl{ce4_n}}(\e)>0$ with $c_{\adr{ce4_n}}(\e)\to 0$.
}
Indeed, note the form $|g+x(\e)|=|g|+\sign(g)x(\e)$ by Claim \ref{cla:|xe|+|tgnee^n|<=}.
We apply the Taylor expansion to the function $F\,:\,\e\mapsto (|g|+\sign(g)x(\e))^{p}$:
\ali
{
|g+x(\e)|^{p}=&F(\e)=F(0)+\sum_{k=1}^{n}\frac{F^{(k)}(0)}{k!}\e^{k}+\left(\frac{F^{(n)}(\alpha \e)}{n!}-\frac{F^{(n)}(0)}{n!}\right)\e^{n}
}
for some $\alpha\in [0,1]$. By virtue of Fa\`a di Bruno formula \cite{Constantine_Savits}, we obtain the equation $F^{(k)}(0)/k!=g_{k,p}$. We will show that the remainder $|F^{(n)}(\alpha \e)/n!-F^{(n)}(0)/n!|$ is bounded by $c_{\adr{ce1_n}}(\e)|g|^{\underline{p}+\eta/2}$ with some constant $c_{\adl{ce1_n}}(\e)\to 0$.
Since the function $g$ is the composition of the two functions $G\,:\,y\mapsto y^p$ and $H\,:\,\e \mapsto |g|+\sign(g)x(\e)$, we have
\alil
{
\frac{F^{(n)}(\alpha \e)}{n!}-&\frac{F^{(n)}(0)}{n!}=\frac{(G\circ H)^{(n)}(\alpha \e)}{n!}-\frac{(G\circ H)^{(n)}(0)}{n!}\nonumber\\
=&\sum_{j=1}^{n}\sum_{\lam_{1},\dots,\lam_{n}\geq 0\,:\,\atop{\lam_{1}+\cdots+\lam_{n}=j\atop{\lam_{1}+2\lam_{2}+\cdots+n\lam_{n}=n}}}\Bigg(G^{(j)}(H(\alpha \e))\prod_{i=1}^{n}\frac{H^{(i)}(\alpha \e)^{\lam_{i}}}{\lam_{i}!(i!)^{\lam_{i}}}-G^{(j)}(H(0))\prod_{i=1}^{n}\frac{H^{(i)}(0)^{\lam_{i}}}{\lam_{i}!(i!)^{\lam_{i}}}\Bigg)\nonumber\\
=&\sum_{j=1}^{n}\sum_{\lam_{1},\dots,\lam_{n}\geq 0\,:\,\atop{\lam_{1}+\cdots+\lam_{n}=j\atop{\lam_{1}+2\lam_{2}+\cdots+n\lam_{n}=n}}}c_{\adr{cgn_n}}\Big((A(\e)-A(0))B_{1}(0)\cdots B_{n}(0)+\nonumber\\
&\quad +\sum_{i=1}^{n}A(\e)B_{1}(\e)\cdots B_{i-1}(\e)(B_{i}(\e)-B_{i}(0))B_{i+1}(0)\cdots B_{n}(0)\Big)\label{eq:gne-gn}
}
putting $c_{\adl{cgn_n}}=\binom{p}{j}j!(\prod_{i=1}^{n}\lam_{i}!(i!)^{\lam_{i}})^{-1}$, 
$A(\e)=H(\alpha\e)^{p-j}$ and $B_{i}(\e)=H^{(i)}(\alpha\e)^{\lam_{i}}$ for $i=1,\dots, n$. We begin with the estimate $A(\e)-A(0)$:
\ali
{
|A(\e)-A(0)|=&|(|g|+\sign(g)x(\alpha\e))^{p-j}-|g|^{p-j}|\\
\leq& |g|^{p-j}||(1+\frac{\sign(g)x(\alpha\e)}{|g|})^{p-j}-1|\qquad (\because \text{Claim }\ref{cla:|xe|+|tgnee^n|<=})\\
\leq&2^{1+|p-j|}|g|^{p-j}c_{\adr{ceta9_n}}(\e)|p-j|\\
\leq&c_{\adr{ce|p-j|_n}}(\e)|g|^{p-j}
}
with $c_{\adl{ce|p-j|_n}}(\e)=2^{1+\max(|\overline{c}-1|,|\underline{c}-n|)}c_{\adr{ceta9_n}}(\e)\max(|\overline{c}-1|,|\underline{c}-n|)$,
where we use the basic inequality $|((1+b)^{p}-1|\leq 2^{1+|p|}|bp|$ whenever $b,p\in \R$ with $|b|<1/2$. Moreover, it follows from $H^{(i)}(\alpha\e)=\sign(g)\sum_{l=0}^{n-i}g_{l+i}\alpha^{l}\e^{l}(l+i)!/(i!)$ and the assumption $\e^{l}<|g|^{\frac{l}{n}(p-\underline{p}-\eta/2)}$
that for $\e\geq 0$,
\ali
{
|H^{(i)}(\alpha\e)|\leq& c_{\adr{g4}}\sum_{l=0}^{n-i}\frac{(l+i)!}{l!}|g|^{t_{l+i}}\e^{l}\\
\leq&c_{\adr{g4}}\sum_{l=0}^{n-i}\frac{(l+i)!}{l!}|g|^{1-\frac{l+i}{n}(p_{l+i}-\underline{p})+\frac{l}{n}(p-\underline{p}-\eta/2)}\\
\leq &c_{\adr{g4}}\sum_{l=0}^{n-i}\frac{(l+i)!}{l!}|g|^{1-\frac{l+i}{n}(p-\underline{p}-\eta)+\frac{l}{n}(p-\underline{p}-\eta/2)}\\
\leq &c_{\adr{csl+i/i_n}}(i)|g|^{1-\frac{i}{n}(p-\underline{p}-\eta)}
}
by putting $c_{\adl{csl+i/i_n}}(i)=c_{\adr{g4}}\sum_{l=0}^{n-i}(l+i)!/l!$ and by using the fact $p_{l+i}=\underline{p}+\frac{n}{l+i}(1-t_{l+i})\geq p(n)>p+\eta$. Therefore
\ali
{
|(A(\e)-A(0))B_{1}(0)\cdots B_{n}(0)|\leq&|H^{(i)}(\alpha\e)-H^{(i)}(0)||H^{(1)}(0)|^{\lam_{1}}\cdots |H^{(n)}(0)|^{\lam_{n}}\\
 \leq&c_{\adr{ce|p-j|_n}}(\e)c_{\adr{csl+i/i_n}}(1)^{\lam_{1}}\cdots c_{\adr{csl+i/i_n}}(n)^{\lam_{n}}|g|^{p-j+\lam_{1}+\cdots+\lam_{n}-\frac{\lam_{1}+2\lam_{2}+\cdots+n\lam_{n}}{n}(p-\underline{p}-\eta)}\\
=&c_{\adr{ce|p-j|_n}}(\e)c_{\adr{csl+i/i_n}}(1)^{\lam_{1}}\cdots c_{\adr{csl+i/i_n}}(n)^{\lam_{n}}|g|^{\underline{p}+\eta}
}
by choosing $\lam_{1},\dots, \lam_{n}$. On the other hand, we note the inequality $\e^{l-c_{\adr{ceta19_n}}}\leq |g|^{\frac{i}{n}(p-\underline{p}-\eta/2)-\frac{\eta}{2n}}\leq |g|^{\frac{i}{n}(p-\underline{p}-\eta)}$
by using the number $c_{\adl{ceta19_n}}=(\eta/2)/(\underline{c}-\underline{p}-\eta/2)$. We get the estimate
\ali
{
|H^{(i)}(\alpha\e)-H^{(i)}(0)|\leq& \sum_{l=1}^{n-i}\frac{(l+i)!}{i!}|g_{l+i}|(\alpha \e)^{l}\\
\leq&c_{\adr{g4}}\sum_{l=1}^{n-i}\frac{(l+i)!}{i!}|g|^{t_{l+i}}|g|^{\frac{i}{n}(p-\underline{p}-\eta)}\e^{c_{\adr{ceta19_n}}}\\
\leq&c_{\adr{csl+i/i_n}}(i)|g|^{1-\frac{l+i}{n}(p_{l+i}-\underline{p})}|g|^{\frac{i}{n}(p-\underline{p}-\eta)}\e^{c_{\adr{ceta19_n}}}\\
\leq&c_{\adr{csl+i/i_n}}(i)|g|^{1-\frac{l+i}{n}(p-\underline{p}-\eta)+\frac{i}{n}(p-\underline{p}-\eta)}\e^{c_{\adr{ceta19_n}}}\\
\leq&c_{\adr{csl+i/i_n}}(i)|g|^{1-\frac{i}{n}(p-\underline{p}-\eta)}\e^{c_{\adr{ceta19_n}}}.
}
Thus
\ali
{
|B_{i}(\e)-B_{i}(0)|=&|H^{(i)}(\alpha\e)^{\lam_{i}}-H^{(i)}(0)^{\lam_{i}}|\\
=&|H^{(i)}(\alpha\e)-H^{(i)}(0)||H^{(i)}(\alpha\e)^{\lam_{i}-1}+H^{(i)}(\alpha\e)^{\lam_{i}-2}H^{(i)}(0)+\cdots+H^{(i)}(0)^{\lam_{i}-1}|\\
\leq&c_{\adr{csl+i/i_n}}(i)^{\lam_{i}}|g|^{\lam_{i}-\frac{i\lam_{i}}{n}(p-\underline{p}-\eta)}\e^{c_{\adr{ceta19_n}}}.
}
Moreover,
\ali
{
|A(\e)|=|(|g|+\sign(g)x(\alpha\e))^{p-j}|\leq &
\case
{
(|g|+|x(\alpha \e)|)^{p-j},&p-j\geq 0\\
(|g|-|x(\alpha \e)|)^{p-j},&p-j< 0\\
}\\
\leq &
\case
{
(\frac{3}{2})^{p-j}|g|^{p-j},&p-j\geq 0\\
(\frac{1}{2})^{p-j}|g|^{p-j},&p-j< 0\\
}\\
\leq&c_{\adr{max3/21/2_n}}|g|^{p-j}
}
with $c_{\adl{max3/21/2_n}}=\max\{(\frac{3}{2})^{\overline{c}-1},(\frac{1}{2})^{\underline{c}-n}\}$. Consequently we obtain
\ali
{
&|A(\e)B_{1}(\e)\cdots B_{i-1}(\e)(B_{i}(\e)-B_{i}(0))B_{i+1}(0)\cdots B_{n}(0)|\\
\leq&c_{\adr{max3/21/2_n}}c_{\adr{csl+i/i_n}}(1)^{\lam_{1}}\cdots c_{\adr{csl+i/i_n}}(n)^{\lam_{n}}\e^{c_{\adr{ceta19_n}}}|g|^{p-j+\lam_{1}+\cdots+\lam_{n}-\frac{\lam_{1}+2\lam_{2}+\cdots+n\lam_{n}}{n}(p-\underline{p}-\eta)}\\
=&c_{\adr{ce2_n}}(\e)|g|^{\underline{p}+\eta}
}
by putting $c_{\adl{ce2_n}}(\e)=c_{\adr{max3/21/2_n}}c_{\adr{csl+i/i_n}}(1)^{\lam_{1}}\cdots c_{\adr{csl+i/i_n}}(n)^{\lam_{n}}\e^{c_{\adr{ceta19_n}}}$. As result, (\ref{eq:gne-gn}) has the following estimate:
\ali
{
\left|\frac{F^{(n)}(\alpha \e)}{n!}-\frac{F^{(n)}(0)}{n!}\right|
\leq &\sum_{j=1}^{n}\sum_{\lam_{1},\dots,\lam_{n}\geq 0\,:\,\atop{\lam_{1}+\cdots+\lam_{n}=j\atop{\lam_{1}+2\lam_{2}+\cdots+n\lam_{n}=n}}}c_{\adr{csl+i/i_n}}(1)^{\lam_{1}}\cdots c_{\adr{csl+i/i_n}}(n)^{\lam_{n}}c_{\adr{cgn_n}}(c_{\adr{ce|p-j|_n}}(\e)+c_{\adr{max3/21/2_n}}\e^{c_{\adr{ceta19_n}}})|g|^{\underline{p}+\eta}\\
=&c_{\adr{ce1_n}}(\e)|g|^{\underline{p}+\eta}.
}
In addition to Claim \ref{cla:|g+...|^p=|g+...|^p+Z}, this estimate implies that the remainder
$\ti{g}_{n,p}(\e)=Z(\e)+(F^{(n)}(\alpha\e)-F^{(n)}(0))/n!$ is bounded by $|\ti{g}_{n,p}(\e)|\leq c_{\adr{ce4_n}}(\e)|g|^{\underline{p}+\ti{t}\eta}$
with $c_{\adr{ce4_n}}(\e)=\max(c_{\adr{ceta11_n}}(\e),c_{\adr{ceta12_n}}(\e))+c_{\adr{ce1_n}}(\e)$. Thus the assertion is valid.
\cla
{
When $n\geq 1$, the main assertion (\ref{eq:supt_tgnt(ew)<=_2}) holds by putting $\e_{1}>0$ so that $c_{\adr{ceta9_n}}(\e)< 1/2$ for $\e<\e_{1}$.
}
Indeed, by virtue of Claim \ref{cla:toyasymp_|f|^p-p_-eta/2<=e^n} and Claim \ref{cla:toyasymp_|f|^p-p_-eta/2>e^n}, this assertion is given with $c_{\adr{tP3}}(\e)=\max(c_{\adr{ctFnex_n}}(\e), c_{\adr{ce4_n}}(\e))$. Hence the proof is complete.
\proe
\lems
\label{lem:asymp_|gew|^tpsi_ex2}
Assume that the incidence matrix of $E^{\infty}$ is finitely irreducible and the conditions $(g.1)$-$(g.5)$ and $(\psi.1)$-$(\psi.4)$ with fixed nonnegative integer $n$ are satisfied. Then for any $p> p(n)$ and there exist functions $\zeta_{1,p},\dots, \zeta_{n,p}\in F_{\theta,b}(E^{\infty},\R)$, $\ti{\zeta}_{n,p}(\e,\cd)\in C_{b}(E^{\infty},\R)$ and constants $\eta\in (0,p-p(n))$, $c_{\adl{tPp1}},c_{\adl{tPp2}},c_{\adl{tPp3}}>0$ such that
\alil
{
|g(\e,\cd)|^{p}\psi(\e,\cd)=|g|^{p}\psi+\zeta_{1,p}\e+\cdots+\zeta_{n,p}\e^{n}+\ti{\zeta}_{n,p}(\e,\cd)\e^{n}\label{eq:geom=_3_2}
}
satisfying with
\ali
{
\zeta_{k,p}=&\sum_{i=0}^{k}g_{i,p}\psi_{k-i}\ (k=0,\dots,n),\quad
\ti{\zeta}_{n,p}(\e,\cdot)=\ti{g}_{n,p}(\e,\om)\psi(\e,\cdot)+\sum_{i=0}^{n}g_{i,p}\ti{\psi}_{n-i}(\e,\cd)
}
that for any $\om,\up\in E^{\infty}$ with $\om_{0}=\up_{0}$ and $k=1,\dots, n$
\alil
{
|\zeta_{k,p}(\om)|\leq &c_{\adr{tPp1}}|g(\om)|^{p-\frac{k}{n}p(n)+\frac{k}{n}\underline{p}}|\psi(\om)|\label{eq:gjt<=_3_2}\\
|\zeta_{k,p}(\om)-\zeta_{k,p}(\up)|\leq &c_{\adr{tPp2}}|g(\om)|^{p-\frac{k}{n}p(n)+\frac{k}{n}\underline{p}}|\psi(\om)|d_{\theta}(\om,\up)\label{eq:|gjtom-gjtup|<=_3_2}
}
and for any nonempty compact subset $I$ of the interval $(p(n),+\infty)$, there exist constants $\eta^\p>0, c_{\adr{tPp3}}(\e)>0$ with $c_{\adr{tPp3}}(\e)\to 0$ such that
\alil
{
&\sup_{p\in I}|\ti{\zeta}_{n,p}(\e,\om)|\leq c_{\adr{tPp3}}(\e)|g(\om)|^{\underline{p}+\eta^\p}|\psi(\om)|\label{eq:tgnte<=_3_2}
}
for any small $\e>0$ and for any $\om\in E^{\infty}$.
\leme
\pros
By using Lemma \ref{lem:asymp_|gew|^t_ex2} and the conditions $(\psi.1)$-$(\psi.4)$ in addition to Proposition \ref{prop:nk-Lip_ex}, we get the assertion.
\proe
Denoted by $\LR(\XR)$ the set of all bounded linear operators acting on a norm space $\XR$.
\lems
\label{lem:asymp_Rop_ge^t_Mgene}
Assume that the incidence matrix of $E^{\infty}$ is finitely irreducible and the conditions $(g.1)$-$(g.5)$ and $(\psi.1)$-$(\psi.4)$ with fixed nonnegative integer $n$ are satisfied. Then for any nonempty compact subset $I\subset (p(n),\infty)$, there exist operators $\LR_{1,p},\dots, \LR_{n,p}\in \LR(F_{\theta,b}(E^{\infty}))$ and $\tLR_{n,p}(\e,\cd)\in \LR(C_{b}(E^{\infty}))$ $(p\in I)$ such that
\alil
{
&\LR_{p\log |g(\e,\cd)|+\log\psi(\e,\cd)}=\LR_{p\log|g|+\log\psi}+\LR_{1,p}\e+\cdots+\LR_{n,p}\e^{n}+\tLR_{n,p}(\e,\cd)\e^{n},\label{eq:asymp_Rop_ge^t-1_2}\\
&\sup_{p\in I}\|\tLR_{n,p}(\e,\cd)\|_{\infty}\to 0\label{eq:asymp_Rop_ge^t-2_2}
}
and $\sup_{p\in I}\|\LR_{k,p}\|_{\theta}<\infty$, where we put
\alil
{
\LR_{k,p}f(\om)=\sum_{e\in E\,:\,t(e)=i(\om_{0})}\zeta_{k,p}(e\cd\om)f(e\cd\om),\quad \ti{\LR}_{n,p}(\e,f)(\om)=\sum_{e\in E\,:\,t(e)=i(\om_{0})}\ti{\zeta}_{n,p}(\e,e\cd\om)f(e\cd\om).\label{eq:Lkt=tLnt=}
}
\leme
\pros
Note that the Ruelle operator of $p\log |g(\e,\cd)|+\log \psi(\e,\cd)$ is given by
\ali
{
\LR_{p\log |g(\e,\cd)|+\log\psi(\e,\cd)}f(\om)
=&\sum_{e\in E\,:\,t(e)=i(\om_{0})}|g(\e,e\cd\om)|^{p}\psi(\e,e\cd\om)f(e\cd\om).
}
By virtue of Lemma \ref{lem:asymp_|gew|^tpsi_ex2}, we get the expansion (\ref{eq:asymp_Rop_ge^t-1_2}) and convergence (\ref{eq:asymp_Rop_ge^t-2_2}).
In remain to check $\sup_{p\in I}\|\LR_{k,p}\|_{\theta}<\infty$ for $k=1,2,\dots, n$.

First we show $\sup_{p\in I}\|\LR_{k,p}\|_{\infty}<\infty$. Indeed, choose any $\eta>0$ so that $p(n)+\eta<\min I$.
From (\ref{eq:gjt<=_3_2}) in addition to Remark \ref{rem:p-k/npn+k/np_>=}, we notice that for any $p\in I$
\ali
{
\|\LR_{k,p}f\|_{\infty}\leq& c_{\adr{tPp1}}\|\LR_{(p-\frac{k}{n}p(n)+\frac{k}{n}\underline{p})\log |g|+\log\psi}1\|_{\infty}\|f\|_{\infty}\\
\leq&c_{\adr{tPp1}}\|\LR_{(\underline{p}+\eta)\log |g|+\log\psi}1\|_{\infty}\|f\|_{\infty}.
}
Furthermore $\|\LR_{(\underline{p}+\eta)\log |g|+\log\psi}1\|_{\infty}<\infty$ by $P((\underline{p}+\eta)\log |g|+\log\psi)<+\infty$ (see Proposition \ref{prop:finitepres_finiteRuelle}). Therefore, we see $\|\LR_{k,p}\|_{\infty}<\infty$ uniformly in $p\in I$.

Next we check the boundedness of $[\LR_{k,p}f]_{\theta}$ for $f\in F_{\theta,b}(E^{\infty})$. It follows from (\ref{eq:gjt<=_3_2}) and (\ref{eq:|gjtom-gjtup|<=_3_2}) that for $\om,\up\in E^{\infty}$ with $\om_{0}=\up_{0}$
\ali
{
&|\LR_{k,p}f(\om)-\LR_{k,p}f(\up)|\\
\leq& \sum_{e\in E\,:\,t(e)=i(\om_{0})}\left(|\zeta_{k,p}(e\cd\om)-\zeta_{k,p}(e\cd\up)||f(e\cd\om)|+|\zeta_{k,p}(e\cd\up)| |f(e\cd\om)-f(e\cd\up)|\right)\\
\leq&c_{\adr{tPp2}}\sum_{e\in E\,:\,t(e)=i(\om_{0})}|g(e\cd\om)|^{p-\frac{k}{n}p(n)+\frac{k}{n}\underline{p}}\psi(e\cd\om)\|f\|_{\infty}d_{\theta}(e\cd\om,e\cd\up)\\
&+c_{\adr{tPp1}}\sum_{e\in E\,:\,t(e)=i(\om_{0})}|g(e\cd\up)|^{p-\frac{k}{n}p(n)+\frac{k}{n}\underline{p}}\psi(e\cd\up)[f]_{\theta}d_{\theta}(e\cd\om,e\cd\up)\\
\leq&\theta (c_{\adr{tPp2}}+ c_{\adr{tPp1}})\|\LR_{(\underline{p}+\eta)\log|g|+\log\psi}1\|_{\infty}\|f\|_{\theta}d_{\theta}(\om,\up).
}
Thus $\sup_{p\in I}[\LR_{k,p}f]_{\theta}\leq \theta(c_{\adr{tPp2}}+c_{\adr{tPp1}})\|\LR_{(\underline{p}+\eta)\log|g|+\log\psi}1\|_{\infty}\|f\|_{\theta}<\infty$. Hence we obtain the boundedness of $\sup_{p\in I}\|\LR_{k,p}\|_{\theta}$.
\proe
\lem
{\label{lem:bdd_e^aphph}
Let $G=(V,E,i(\cd),t(\cd))$ be a directed multigraph such that the incidence matrix of $E^{\infty}$ is finitely irreducible. Let $\ph\in F_{\theta}(E^{\infty},\R)$ satisfying (g.2) and (g.3). Then for any number $\eta>0$ and integer $k\geq 1$, the function $\om\mapsto |g(\om)|^{\eta}(\log|g(\om)|)^{k}$ belongs to $F_{\theta,b}(E^{\infty})$. In particular, $\||g|^{\eta}(\log|g|)^{k}\|_{\infty}\leq c_{\adr{cg^etalogg^k}}$ and $[|g|^{\eta}(\log|g|)^{k}]_{\theta}\leq c_{\adr{cg^etalogg^k2}}$ for some constants $c_{\adl{cg^etalogg^k}}=c_{\adr{cg^etalogg^k}}(\eta,k), c_{\adl{cg^etalogg^k2}}=c_{\adr{cg^etalogg^k2}}(\eta,k)>0$.
}
\pros
Let $\alpha=\eta/k$. We have that for $\om\in E^{\infty}$
\ali
{
||g(\om)|^{\alpha}\log|g(\om)||=\frac{1}{\alpha}\left||g(\om)|^{\alpha}\log(|g(\om)|^{\alpha})\right|\leq \frac{1}{\alpha e}.
}
Thus $\||g|^{\eta}(\log|g|)^{k}\|_{\infty}\leq c_{\adr{cg^etalogg^k}}:=k^{k}/(e^{k}\eta^{k})$. Moreover, for $\om,\up\in E^{\infty}$ with $\om_{0}=\up_{0}$
\ali
{
||g(\om)|^{\alpha}\log|g(\om)|-|g(\up)|^{\alpha }\log|g(\up)||\leq& |g(\om)|^{\alpha}|\log|g(\om)|-\log|g(\up)||\\
&+|\log(|g(\up)|)|g(\up)|^{\alpha}|e^{\alpha\log|g(\om)|-\alpha\log|g(\up)|}-1|\\
\leq& ([\log|g|]_{\theta}+(\alpha e)^{-1}e^{\alpha[\log|g|]_{\theta}\theta}\alpha[\log|g|]_{\theta})d_{\theta}(\om,\up).
}
Thus $|g|^{\alpha}\log|g|$ is in $F_{\theta,b}(E^{\infty},\R)$. Hence we obtain 
\ali
{[e^{\eta \ph}\ph^{k}]_{\theta}\leq& k([\log|g|]_{\theta}+(\alpha e)^{-1}e^{\alpha[\log|g|]_{\theta}\theta}\alpha[\log|g|]_{\theta})(e \eta/k)^{-k+1}\\
\leq&(kc_{\adr{g3}}+(\eta e/k)^{-1}e^{\eta c_{\adr{g3}}\theta/k}\eta c_{\adr{g3}})(e \eta/k)^{-k+1}
}
by using $[\log|g|]_{\theta}\leq c_{\adr{g3}}$ and by putting $c_{\adr{cg^etalogg^k2}}:=k^{k}c_{\adr{g3}}(e\eta)^{-k+1}(1+e^{\eta c_{\adr{g3}}\theta/k-1})$.
\proe
The function $(p,\e)\mapsto \LR_{p\log g(\e,\cd)+\log\psi(\e,\cd)}$ has also an asymptotic expansion in the sense of the following lemma.
\lems
\label{lem:asymp_op_Mgene}
Assume that the incidence matrix of $E^{\infty}$ is finitely irreducible and the conditions $(g.1)$-$(g.5)$ and $(\psi.1)$-$(\psi.4)$ with fixed nonnegative integer $n$ are satisfied. Choose any nonempty compact subset $I\subset (p(n),+\infty)$. Then for any $s,p\in I$, there exist operators $\ZR_{v,q,s}\in \LR(F_{\theta,b}(E^{\infty}))$ $(0\leq v,q\leq n)$ and $\hat{\ZR}_{n,s,p}\in \LR(C_{b}(E^{\infty}))$ such that the Ruelle operator of $p\log|g(\e,\cd)|+\log\psi(\e,\cd)$ is expanded as
\alil
{
\LR_{p\log|g(\e,\cd)|+\log\psi(\e,\cd)}=&\sum_{v=0}^{n}\sum_{q=0}^{n}\ZR_{v,q,s}\e^{v}(p-s)^{q}+(p-s)^{n+1}\hat{\ZR}_{n,s,p}+\ti{\LR}_{n,p}(\e,\cd)\e^{n}\label{eq:asympsol_Beq_1_2}
}
and $\sup_{s,p\in I}\|\hat{\ZR}_{n,s,p}\|_{\infty}<+\infty$, where $\ZR_{0,0,s}$ equals $\LR_{s\log|g|+\log\psi}$.
Here operators $\ZR_{v,q,s}$ and $\hat{\ZR}_{n,s,p}$ are given by
\ali
{
\ZR_{v,q,s}f:=&\LR_{0}(h_{v,q,s}f) \text{ for }f\in C_{b}(E^{\infty}) \text{ with}\\
h_{v,q,s}:=&\sum_{k=0}^{v}\sum_{l=0}^{k}\sum_{j=0}^{\min(l,q)}\frac{a_{l,j,s}}{(q-j)!}(\log |g|)^{q-j}G^{s}_{k,l}\psi_{v-k}
}
and
\alil
{
\hat{\ZR}_{n,s,p}f:=&\sum_{v=0}^{n}\LR_{0}(\hat{h}_{v,s,p} f) \e^{v}\text{ with}\\
\hat{h}_{v,s,p}:=&\sum_{k=0}^{v}\sum_{l=0}^{k}\sum_{j=0}^{l}a_{l,j,s}\hat{\varGamma}_{s,p,n-j+1}G^{s}_{k,l}\psi_{v-k}\label{eq:hh_vsp=}\\
\hat{\varGamma}_{s,p,i}(\om):=&\int_{0}^{1}\frac{(1-u)^{i-1}}{(i-1)!}|g(\om)|^{u(p-s)}(\log|g(\om)|)^{i}\,du,
}
where $\LR_{0}f$ means $\sum_{e\in E\,:\,t(e)=i(\om_{0})}f(e\cd\om)$, and $a_{l,j,s}$ are numbers defined in (\ref{eq:aljs=}) below.
\leme
\pros
First we show the expansion (\ref{eq:asympsol_Beq_1_2}).
We take the notation $\underline{c}=\inf I$, $\overline{c}=\sup I$, $\eta\in (0,\underline{c}-p(n))$ and $p_{k}=\underline{p}+\frac{n}{k}(1-t_{k})$.
By Theorem \ref{lem:asymp_Rop_ge^t_Mgene}, we have the expansion
\ali
{
\LR_{p\log|g(\e,\cd)|+\log\psi(\e,\cd)}=\sum_{v=0}^{n}\LR_{v,p}\e^{v}+\ti{\LR}_{n,p}(\e,\cd)\e^{n},
}
where $\LR_{v,p}$ and $\ti{\LR}_{n,p}$ are defined by (\ref{eq:Lkt=tLnt=}).
Now we will extend each $\LR_{v,p}$ for $p$. To do this, we remark the following expansions $\binom{p}{l}=\sum_{j=0}^{l}a_{l,j,s}(p-s)^{j}$ by putting
\alil
{
a_{l,j,s}=&
\case
{
\displaystyle \binom{s}{l}& (j=0)\\
\displaystyle
\sum_{0\leq  i_{1},\cdots,i_{l-j}\leq l-1:\atop{i_{1}<\cdots <i_{l-j}}}\frac{1}{l!}\prod_{q=1}^{l-j}(s-i_{q})& (l\geq 1 \text{ and }0\leq j<l)
\medskip
\\
\displaystyle 1/l!&(l\geq 1 \text{ and }j=l)\\
0&(l<j).
}\label{eq:aljs=}
}
Moreover, Taylor expansion for $p\mapsto |g|^{p}$ implies that for any $i\geq 0$,
\ali
{
|g|^{p}=&e^{p\log|g|}=\sum_{q=0}^{i}\frac{(\log|g|)^{q}}{q!}|g|^{s}(p-s)^{q}+(p-s)^{i+1}|g|^{s}\hat{\varGamma}_{s,p,i+1}.
}
Therefore the function $\zeta_{v,p}=\sum_{k=0}^{v}\psi_{v-k}g_{k,p}$ has the form 
\ali
{
\zeta_{v,p}
=&\sum_{k=0}^{v}\psi_{v-k}\sum_{l=0}^{k}\sum_{j=0}^{l}a_{l,j,s}(p-s)^{j}G^{p}_{k,l}\\
=&\sum_{k=0}^{v}\sum_{l=0}^{k}\sum_{j=0}^{l}\psi_{v-k}a_{l,j,s}\Big(\sum_{q=0}^{n-j}\frac{(\log|g|)^{q}}{q!}G^{s}_{k,l}(p-s)^{q+j}+\hat{\varGamma}_{s,p,n-j+1}G^{s}_{k,l}(p-s)^{n+1}\Big)\\
=&\sum_{q=0}^{n}\sum_{k=0}^{v}\sum_{l=0}^{k}\sum_{j=0}^{\min(l,q)}a_{l,j,s}\frac{(\log|g|)^{q-j}}{(q-j)!}G^{s}_{k,l}\psi_{v-k}(p-s)^{q}+\hat{h}_{v,s,p}(p-s)^{n+1}\\
=&\sum_{q=0}^{n}h_{v,q,s}(p-s)^{q}+\hat{h}_{v,s,p}(p-s)^{n+1}.
}
Thus the equation (\ref{eq:asympsol_Beq_1_2}) is valid.

Next we will prove $\ZR_{v,q,s}\in \LR(F_{\theta,b}(E^{\infty}))$. In the expression of $h_{v,q,s}$, we rewrite $(\log|g|)^{q-j}G^{s}_{k,l}\psi_{v-k}=((\log|g|)^{q-j}|g|^{\eta/2})|g|^{-\eta/2}G^{s}_{k,l}\psi_{v-k}$. Then the function $(\log|g|)^{q-j}|g|^{\eta/2}$ is in $F_{\theta,b}(E^{\infty},\R)$ from Lemma \ref{lem:bdd_e^aphph}. 
By Proposition \ref{prop:nk-Lip_ex} with $f_{0}=(\log|g|)^{q-j}|g|^{\eta/2}$, $h_{0}=1$, $f_{1}=|g|^{-\eta/2}$, $h_{1}=|g|^{-\eta/2}$, $f_{2}=G^{s}_{k,l}$, $h_{2}=|g|^{s+\eta}$, $f_{3}=\psi_{v-k}$, $h_{3}=\psi$, $c_{\adr{nk-L_fk1}}(0)=c_{\adr{cg^etalogg^k}}(\eta/2,q-j)$, $c_{\adr{nk-L_fk1}}(1)=1$, 
$c_{\adr{nk-L_fk1}}(2)=c_{\adr{cGtlk}}(k,l)$, $c_{\adr{nk-L_fk1}}(3)=c_{\adr{psi3}}$, $c_{\adr{nk-L_fk2}}(0)=c_{\adr{cg^etalogg^k2}}(\eta/2,q-j)$, $c_{\adr{nk-L_fk2}}(1)=c_{\adr{g3}}(1+c_{\adr{g3}}\theta)^{\eta/2}$, $c_{\adr{nk-L_fk2}}(2)=c_{\adr{cGtlk2}}(k,l)$, and $c_{\adr{nk-L_fk2}}(3)=c_{\adr{psi3_2}}$, we get the estimate
\ali
{
&|(\log|g|)^{q-j}G^{s}_{k,l}\psi_{v-k}|\leq c_{\adr{cgGpsi}}(k,l,j)|g|^{s+\eta/2}\psi\\
&|(\log|g(\om)|)^{q-j}G^{s}_{k,l}(\om)\psi_{v-k}(\om)-(\log|g(\up)|)^{q-j}G^{s}_{k,l}(\up)\psi_{v-k}(\up)|\leq c_{\adr{cgGpsi2}}|g(\om)|^{s+\eta/2}d_{\theta}(\om,\up)
}
for the constants $c_{\adl{cgGpsi}}=c_{\adr{cgGpsi}}(k,l,j)=c_{\adr{cg^etalogg^k}}(\eta/2,q-j)c_{\adr{cGtlk}}(k,l)c_{\adr{psi3}}$ and $c_{\adl{cgGpsi2}}=c_{\adr{cgGpsi2}}(k,l,j)=c_{\adr{nk-L_fk3}}$ with $n=3$.
Thus we obtain
\ali
{
\|\ZR_{v,q,s}\|_{\infty}&\leq c_{\adr{cHvqs}}\|\LR_{(\underline{p}+\eta/2)\log|g|+\log|\psi|}1\|_{\infty}\\
[\ZR_{v,q,s}f]_{\theta}&\leq (c_{\adr{cHvqs2}}\theta\|f\|_{\infty}+c_{\adr{cHvqs}}\theta[f]_{\theta})\|\LR_{(\underline{p}+\eta/2)\log|g|+\log|\psi|}1\|_{\infty},
}
where we define constants $c_{\adl{cHvqs}}$ and $c_{\adl{cHvqs2}}$ by
\ali
{
c_{\adr{cHvqs}}=\sum_{k=0}^{v}\sum_{l=0}^{k}\sum_{j=0}^{\min(l,q)}\frac{a_{l,j,s}}{(q-j)!}c_{\adr{cgGpsi}}(k,l,j),\qquad c_{\adr{cHvqs2}}=\sum_{k=0}^{v}\sum_{l=0}^{k}\sum_{j=0}^{\min(l,q)}\frac{a_{l,j,s}}{(q-j)!}c_{\adr{cgGpsi2}}(k,l,j).
}
Hence $\ZR_{v,q,s}\in \LR(F_{\theta,b}(E^{\infty},\R))$ is guaranteed.

It remain to check the boundedness of $\hat{\ZR}_{n,s,p}$. In the expression of $\hat{h}_{v,s,p}$, we have
\ali
{
&|\hat{\varGamma}_{s,p,n-j+1}(\om)G^{s}_{k,l}(\om)\psi_{v-k}(\om)|\\
\leq&c_{\adr{cGtlk}}(k,l)c_{\adr{psi3}}|g(\om)|^{s-\frac{k}{n}p(n)+\frac{k}{n}\underline{p}}\psi(\om)|\hat{\varGamma}_{s,p,n-j+1}(\om)|\\
\leq&c_{\adr{cGtlk}}(k,l)c_{\adr{psi3}}|g(\om)|^{\underline{p}+\eta}\psi(\om)\int_{0}^{1}\frac{(1-u)^{n-j}}{(n-j)!}|\log(|g(\om)|)|^{n-j+1}\,du\quad (\because \text{Remark }\ref{rem:p-k/npn+k/np_>=})\\
\leq&\frac{c_{\adr{cGtlk}}(k,l)c_{\adr{psi3}}c_{\adr{cg^etalogg^k}}(\eta/2,n-j+1)}{(n-j+1)!}\psi(\om)|g(\om)|^{\underline{p}+\eta/2}\quad (\because \text{ Lemma }\ref{lem:bdd_e^aphph}).
}
As a result, we get
$|\hat{h}_{v,s,p}(\om)|\leq c_{\adr{chHvst}}\psi(\om)|g(\om)|^{\underline{p}+\eta/2}$ for some constant $c_{\adl{chHvst}}>0$.
Thus $\hat{\ZR}_{n,s,p}$ is bounded by $c_{\adr{chHvst}}\|\LR_{(\underline{p}+\eta/2)\log|g|+\log\psi}1\|_{\infty}<+\infty$ uniformly in $s,p\in I$.
\proe
Now we describe a remark concerning the coefficient of the solution $s(\e)$ before the proof of Theorem \ref{th:asympsol_Beq_Mgene}.
\rem
{\label{rem:coeff_sk}
Assume that the incidence matrix of $E^{\infty}$ is finitely irreducible and the conditions $(g.1)$-$(g.5)$ and $(\psi.1)$-$(\psi.4)$ with fixed nonnegative integer $n$ are satisfied. The coefficient $s_{k}$ $(k=1,\dots, n)$ and remainder $\ti{s}_{n}(\e)$ of the solution $s(\e)$ are given by
\alil
{
s_{k}=&\frac{-1}{\nu(h\log|g|)}\Big(\sum_{i=1}^{k-1}\nu_{i}(\LR_{s(0)\log|g|+\log\psi|}(h\log|g|))s_{k-i}+\sum_{i=0}^{k-1}\nu_{i}(\NR_{k-i}h)\Big)\label{eq:tn=...}\\
\ti{s}_{0}(\e)=&-\frac{\nu(\e,\tLR_{0,s(\e)}(\e,h))}{\nu(\e,\hat{\ZR}_{0,s(0),s(\e)}h)} \text{ and if }n\geq 1 \text{ then}\label{eq:ttne=...n=0}\\
\ti{s}_{n}(\e)=&\frac{-1}{\nu(h\log|g|)}\Big(\sum_{u=1}^{n}\ti{\nu}_{n-u}(\e,\NR_{u}h)+\sum_{i=1}^{n-1}\ti{s}_{n-i}(\e)\nu_{i}(\ZR_{0,1,s(0)}h)\label{eq:ttne=...n>=1}\\
&+\nu(\e,\ti{\LR}_{n,s(\e)}(\e,h))+\frac{s(\e)-s(0)}{\e}\ti{\nu}_{n-1}(\e,\ZR_{0,1,s(0)}h)+\nu(\e,\hat{\NR}_{n+1}(\e,h))\e\Big).\nonumber
}
Here
$\nu, \nu_{k}$, $\nu(\e,\cd)$ and $\ti{\nu}_{n-u}(\e,\cd)$ appear in the asymptotic expansions $\nu(\e,f)=\nu(f)+\sum_{k=1}^{m}\nu_{k}(f)\e^{k}+\ti{\nu}_{m}(\e,f)\e^{m}$ $(0\leq m\leq n-1)$ of the Perron eigenvector $\nu(\e,\cd)$ of $\LR_{s(\e)\log|g(\e,\cd)|+\log\psi}$ given by Corollary \ref{cor:asymp_e.vec}. We define operators $\NR_{u}, \hat{\NR}_{n+1}(\e,\cd)\in \LR(C_{b}(E^{\infty}))$ by
\alil
{
\NR_{u}=&\sum_{0\leq v\leq u,\ 0\leq q\leq u-v\,:\,(v,q)\neq (0,1)}s_{q,u-v}\ZR_{v,q,s(0)}\ \ (u=1,\dots, n)\label{eq:NRu=}\\
\hat{\NR}_{n+1}(\e,f)=&\sum_{0\leq v,q\leq n\,:\,(v,q)\neq (0,1)}\ti{s}_{q,n-v}(\e)\ZR_{v,q,s(0)}f+\hat{\ZR}_{n,s(0),s(\e)}f\bigg(\frac{s(\e)-s(0)}{\e}\bigg)^{n+1},\label{eq:hNRn+1=}
}
where $s_{q,i}$ and $\ti{s}_{q,i}(\e)$ are the coefficient and the remainder of the expansion $(s(\e)-s(0))^{k}=\sum_{i=0}^{j}s_{q,i}\e^{i}+\ti{s}_{q,j}(\e)\e^{j}$, respectively (see (\ref{eq:tki=}) and (\ref{eq:ttkne=}) for detail).
}
\noindent
(Proof of Theorem \ref{th:asympsol_Beq_Mgene}).
\pros
Put $\Ph(\e,s,\cd)=s\log|g(\e,\cd)|+\log\psi(\e,\cd)$ and $\Ph(s,\cd)=s\log|g|+\log\psi$ for convenience. For $s\in I$ and for $\e>0$, let $(\lam_{s}(\e),h_{s}(\e,\cd),\nu_{s}(\e,\cd))$ be the triplet of the operator $\LR_{\Ph(\e,s,\cd)}$ and $(\lam_{s},h_{s},\nu_{s})$ the triplet of $\LR_{\Ph(s,\cd)}$ given by Theorem \ref{th:exGibbs}. For convenience, we may write $\nu(\e,\cd):=\nu_{s(\e)}(\e,\cd)$, $\nu:=\nu_{s(0)}$ and $h:=h_{s(0)}$.
\ncla
{
The solution $s(\e)$ exists in $I$ for any small $\e>0$, and converges to $s(0)$ as $\e\to 0$.
}
Indeed, it follows from (\ref{eq:asympsol_Beq_1_2}) with $n=0$, $a_{0,0,s}=1$, $G_{0,0}^{s}=|g|^{s}$ and $\psi_{0}=\psi$ that
\ali
{
\LR_{\Ph(\e,s,\cd)}h=&\LR_{\Ph(s(0),\cd)}h+(s-s(0))\hat{\ZR}_{0,s(0),s}h+\tLR_{0,s}(\e,h)\\
\hat{\ZR}_{0,s(0),s}h=&\LR_{0}(\int_{0}^{1}|g|^{u(s-s(0))+s(0)}\,du\log|g|\psi h)
}
By using the equations $\LR_{\Ph(\e,s,\cd)}^{*}\nu_{s}(\e,\cd)=\lam_{s}(\e)\nu_{s}(\e,\cd)$ and $\LR_{\Ph(s(0),\cd)}h=\lam_{s(0)}h$, we obtain
\alil
{
\nu_{s}(\e,h)(\lam_{s}(\e)-\lam_{s(0)})=\nu_{s}(\e,\hat{\ZR}_{0,s(0),s}h)(s-s(0))+\nu_{s}(\e,\tLR_{0,s}(\e,h)).\label{eq:nue(LRtelogge...)=}
}
Now choose any small $\eta>0$ so that $s(0)+\eta,s(0)-\eta\in I$. Take $\overline{c}=\sup I$. For any $s\in [s(0)-\eta,s(0)+\eta]$, we have the estimate
\alil
{
-\nu_{s}(\e,\hat{\ZR}_{0,s(0),s}h_{s(0)})
\geq &
\case
{
\nu_{s}(\e,\LR_{0}(|g|^{s}\psi (-\log|g|) h_{s(0)}),&s\geq s(0)\\
\nu_{s}(\e,\LR_{0}(|g|^{s(0)}\psi (-\log|g|) h_{s(0)}),&s<s(0)\\
}\nonumber\\
\geq&\case
{
(-\log\|g\|_{\infty})\nu_{s}(\e,\LR_{\Ph(\overline{c},\cd)}h_{\overline{c}}\frac{h_{s(0)}}{h_{\overline{c}}}),&s\geq s(0)\\
(-\log\|g\|_{\infty})\lam_{s(0)}\nu_{s}(\e,h_{s(0)}),&s<s(0)\\
}\quad (\because s\leq \overline{c})\nonumber\\
\geq&-c_{\adr{cvteZ}}\log\|g\|_{\infty}>0\label{eq:nuthZ0th<0}
}
with $c_{\adl{cvteZ}}=\min(\lam_{\overline{c}}\inf_{\om}h_{\overline{c}}(\om)/\|h_{\overline{c}}\|_{\infty},\lam_{s(0)})\inf_{\om}h_{s(0)}(\om)$. Now fix $s\in [s(0)-\eta,s(0)+\eta]$ with $s\neq s(0)$. For any small $\e>0$ with $\|\tLR_{0,s}(\e,h_{s(0)})\|_{\infty}/|s-s(0)|\leq -c_{\adr{cvteZ}}\log(\|g\|_{\infty})/2$, the equation (\ref{eq:nue(LRtelogge...)=}) yields
\ali
{
\nu_{s}(\e,h_{s(0)})\frac{\lam_{s}(\e)-\lam_{s(0)}}{s-s(0)}=&\nu_{s}(\e,\hat{\ZR}_{0,s(0),s(\e)}h_{s(0)})+\frac{\nu_{s}(\e,\tLR_{0,s(\e)}(\e,h_{s(0)}))}{s-s(0)}\leq\frac{c_{\adr{cvteZ}}}{2}\log\|g\|_{\infty}<0.
}
In addition to $\nu_{s}(\e,h_{s(0)})>0$, this implies that $\lam_{s}(\e)<\lam_{s(0)}$ if $s>s(0)$ and $\lam_{s}(\e)>\lam_{s(0)}$ if $s<s(0)$ for any small $\e>0$. Since $\log\lam_{s}(\e)=P(\Ph(\e,s,\cd))$ and $\log\lam_{s(0)}=P(\Ph(s(0),\cd))=p_{0}$ are satisfied, we obtain
\ali
{
P(\Ph(\e,s(0)+\eta,\cd))<p_{0}<P(\Ph(\e,s(0)-\eta,\cd))
}
for a fix number $\eta>0$ with $s(0)+\eta,s(0)-\eta\in I$ and for any small $\e>0$. It follows from this inequality and the strictly monotone decreasingly of the map $s\mapsto P(\Ph(\e,s,\cd))$ that there exists a unique $s(\e)\in [s(0)-\eta,s(0)+\eta]$ so that $P(\Ph(\e,s(\e),\cd))=p_{0}$. By arbitrary choosing $\eta>0$, $s(\e)$ converges to $s(0)$ as $\e\to 0$. In particular, $\ti{s}_{0}(\e)=s(\e)-s(0)$ estimates as (\ref{eq:ttne=...n=0}) from (\ref{eq:nue(LRtelogge...)=}).
\cla
{\label{cl:te=t0+Oe}
$s(\e)=s(0)+O(\e)$ as $\e\to 0$ in the case $n\geq 1$.
}
Indeed, since the form $\tLR_{0,s(\e)}(\e,\cd)=\LR_{1,s(\e)}\e+\tLR_{1,s(\e)}(\e,\cd)\e$ is satisfied from (\ref{eq:asymp_Rop_ge^t-1_2}), (\ref{eq:ttne=...n=0}) implies 
\ali
{
\frac{s(\e)-s(0)}{\e}=\frac{-\nu(\e,\LR_{1,s(\e)}+\tLR_{1,s(\e)}(\e,\cd))}{\nu(\e,\hat{\ZR}_{0,s(0),s(\e)}h)}.
}
By virtue of Theorem \ref{th:asymp_e.vec} by regarding $(\XR_{0},\|\cd\|_{0})=(C_{b}(E^{\infty}),\|\cd\|_{\infty})$ and $(\XR_{1},\|\cd\|_{1})=(F_{\theta,b}(E^{\infty}),\|\cd\|_{\theta})$, $\LR(\e,\cd)=\LR_{\Ph(\e,s(\e),\cd)}$ and $\LR=\LR_{\Ph(s(0),\cd)}$,  the measure $\nu(\e,\cd)$ converges to $\nu$ weakly. In fact, the conditions (L.1) and (L.2) are satisfied by Theorem \ref{th:exGibbs}, the condition (L.3) is yielded by $\|\nu(\e,\cd)\|=1$ and $0<\inf h\leq \sup h<\infty$, and the condition (L.4) follows from $\LR_{\Ph(\e,s(\e),\cd)} \to \LR_{\Ph(s(0),\cd)}$ in $C_{b}(E^{\infty})$. 
Moreover, the operator $\LR_{1,s(\e)}$ is bounded by Lemma \ref{lem:asymp_Rop_ge^t_Mgene} for any small $\e>0$.
In addition to (\ref{eq:nuthZ0th<0}) and (\ref{eq:asymp_Rop_ge^t-2_2}), it follows from the boundedness $\LR_{1,s(\e)}$ and convergence $\nu(\e,\cd)$ that $(s(\e)-s(0))/\e$ is bounded.
\cla
{
The asymptotic expansion (\ref{eq:t(e)=t0+t1e+...}) of $s(\e)$ is satisfied in the case $n\geq 1$.
}
Let $n\geq 1$. Assume that $s(\e)$ has the form $(n-1)$-asymptotic expansion $s(\e)=s_{0}+s_{1}\e+\cdots+s_{n-1}\e^{n-1}+o(\e^{n-1})$ with $s_{0}=s(0)$. 
To see the $n$-asymptotic behavior of $s(\e)$, we need the expansion of $(s(\e)-s(0))^{k}$ for $k\geq 0$:
\ali
{
(s(\e)-s(0))^{k}=
\case
{
s_{1,0}+s_{1,1}\e+\dots+s_{1,n-1}\e^{n-1}+\ti{s}_{n-1}(\e)\e^{n-1}&(k=1)\\
s_{k,0}+s_{k,1}\e+\dots+s_{k,n-1}\e^{n-1}+s_{k,n}\e^{n}+\ti{s}_{k,n}(\e)\e^{n}&(k\geq 2)\\
}
}
with 
\alil
{
s_{k,i}=&
\case
{
1&(k=i=0)\\
0&(k=0 \text{ and }i\geq 1)\\
0&(k\geq 1 \text{ and }0\leq i\leq k-1)\\
s_{i}&(k= 1 \text{ and }1\leq i\leq n-1)\\
\di\sum_{j_{1},\cdots,j_{i-1}\geq 0\,:\,\atop{j_{1}+\cdots+j_{i-1}=k\atop{j_{1}+2j_{2}+\cdots+(i-1)j_{i-1}=i}}}\frac{s_{1}^{j_{1}}\cdots s_{i-1}^{j_{i-1}}}{j_{1}!\cdots j_{i-1}!}& (k\geq 2 \text{ and }k\leq i\leq n)\\
}\label{eq:tki=}\\
\ti{s}_{k,n}(\e)=&\sum_{i=n+1}^{k(n-1)}\sum_{j_{1},\dots,j_{n}\geq 0\,:\,\atop{j_{1}+\cdots+j_{n}=k\atop{j_{1}+2j_{2}+\cdots+(n-1)j_{n-1}+(n-1)j_{n}=i}}}\frac{s_{1}^{j_{1}}\cdots s_{n-1}^{j_{n-1}}\ti{s}_{n-1}(\e)^{j_{n}}}{j_{1}!\cdots j_{n}!}\e^{i-n}.\label{eq:ttkne=}
}
Note that $s_{1,i}=s_{i}$ holds for $1\leq i\leq n-1$. Thus the expansion (\ref{eq:asympsol_Beq_1_2}) implies
\alil
{
&\LR_{\Ph(\e,s(\e),\cd)}f\label{eq:Ltelogg=_asymp}\\
=&\sum_{v=0}^{n}\sum_{q=0}^{n}(s(\e)-s(0))^{q}\ZR_{v,q,s(0)}f\e^{v}+(s(\e)-s(0))^{n+1}\hat{\ZR}_{n,s(0),s(\e)}f+\ti{\LR}_{n,s(\e)}(\e,f)\e^{n}\nonumber\\
=&(s(\e)-s(0))\ZR_{0,1,s(0)}f +\sum_{0\leq v,q\leq n\,:\,\atop{(v,q)\neq (0,1)}}(s(\e)-s(0))^{q}\ZR_{v,q,s(0)}f \e^{v}\nonumber\\
&+(s(\e)-s(0))^{n+1}\hat{\ZR}_{n,s(0),s(\e)}f+\ti{\LR}_{n,s(\e)}(\e,f)\e^{n}\nonumber\\
=&\LR_{\Ph(s(0),\cd)}f+\ZR_{0,1,s(0)}f (s(\e)-s(0))+\sum_{u=1}^{n}\NR_{u}f\e^{u}+\hat{\NR}_{n+1}(\e,f)\e^{n+1}+\ti{\LR}_{n,s(\e)}(\e,f)\e^{n}\nonumber
}
for $f\in C_{b}(E^{\infty})$ using (\ref{eq:NRu=}) and (\ref{eq:hNRn+1=}).
By the definition of $\hat{\NR}_{n+1}(\e,\cd)$, this operator is bounded uniformly in any small $\e>0$.
Since $s(\e)$ has $(n-1)$-asymptotic expansion, we see by (\ref{eq:Ltelogg=_asymp}) that $\LR_{\Ph(\e,s(\e),\cd)}$ at least has $(n-1)$-asymptotic expansion in $\LR(C_{b}(E^{\infty}))$. Thus it follows from Corollary \ref{cor:asymp_e.vec} that $\nu(\e,\cd)$ has the form $\nu(\e,\cd)=\nu+\sum_{k=1}^{n-1}\nu_{k} \e^{k}+\ti{\nu}_{n-1}(\e,\cd)\e^{n-1}$ and $|\ti{\nu}_{n-1}(\e,f)|\to 0$ for each $f\in F_{\theta,b}(E^{\infty})$. We have
\ali
{
0=&\nu(\e,(e^{p_{0}}-e^{p_{0}})h)=\nu(\e,(\LR_{\Ph(\e,s(\e),\cd)}-\LR_{\Ph(s(0),\cd)})h)\\
=&(s(\e)-s(0))\nu(\e,\ZR_{0,1,s(0)}h)+\sum_{u=1}^{n}\nu(\e,\NR_{u}h)\e^{u}+\nu(\e,\hat{\NR}_{n+1}(\e,h))\e^{n+1}+\nu(\e,\ti{\LR}_{n,s(\e)}(\e,h))\e^{n}.
}
Consequently, 
we get the form $\ti{s}_{n-1}(\e)=s_{n}\e+\ti{s}_{n}(\e)\e$ by putting (\ref{eq:tn=...}) with $k=n$ and (\ref{eq:ttne=...n>=1}) and $\ti{s}_{n}(\e)$ vanishes. Thus this claim is satisfied.
\cla
{
The estimate (\ref{eq:tsne=}) of the remainder $\ti{s}_{n}(\e)$ is valid.
}
First assume $n=0$. Recall the form (\ref{eq:ttne=...n=0}) of $\ti{s}_{0}(\e)$.
Since $s(\e)$ converges to $s(0)$, then $\LR_{s(\e)\log|g(\e,\cd)|+\log\psi(\e,\cd)}$ converges to $\LR_{s(0)\log|g|+\log\psi}$ in $C_{b}(E^{\infty},\C)$ as $\e\to 0$ by Lemma \ref{lem:asymp_op_Mgene}. Therefore the measure $\nu(\e,\cd)$ converges to $\nu$ weakly from Corollary \ref{cor:asymp_e.vec}. Moreover, it follows from the Mean Valued Theorem that $\hat{\varGamma}_{s(0),s(\e),1}(\e,\om)=\int_{0}^{1}|g(\om)|^{u(s(\e)-s(0))}\log|g(\om)|\,du$ converges to $\log|g(\om)|$ uniformly in $\om\in E^{\infty}$. Thus $\nu(\e,\hat{\ZR}_{0,s(0),s(\e)}h)\to \nu(h\log |g|)$. The assertion is valid in the case $n=0$.

Next assume $n\geq 1$. Since $s(\e)$ has an $n$-order asymptotic expansion, $\LR_{s(\e)\log|g(\e,\cd)|+\log\psi(\e,\cd)}$ also has an $n$-order asymptotic expansion for $\e$ in $\LR(C_{b}(E^{\infty}))$ from the expansion (\ref{eq:Ltelogg=_asymp}). Thus Corollary \ref{cor:asymp_e.vec} says that $\nu(\e,\cd)$ has an $n$-order asymptotic expansion. 
In the expression (\ref{eq:ttne=...n>=1}) of $\ti{s}_{n}(\e)$, we notice the form $\ti{s}_{n-1}(\e)=(s_{n}+\ti{s}_{n}(\e))\e$ and $\ti{\nu}_{n-1}(\e,f)=(\nu_{n}(f)+\ti{\nu}_{n}(\e,f))\e$. Hence the proof is complete.
\proe
\subsection{Proof of Proposition \ref{prop:asympsol_Beq_est_rem2}}
\pros
Choose any compact neighborhood $I$ of $s_{0}$ so that $I\subset (p(n),\underline{p}+(n+1)(1-t_{0}))\setminus \Z$. Put $\underline{c}=\inf I$, $\overline{c}=\sup I$ and 
$x(\e,\om)=|g(\e,\om)|-|g(\om)|$. 
We begin with the estimate of $x(\e,\om)$.
\ncla
{\label{cl:xeom>=xg^t0}
There exists $c_{\adl{AS3}}>0$ such that for any $\om \in E^{\infty}$, $c_{\adr{AS3}}|g(\om)|^{t_{0}}\e\leq x(\e,\om)$ for any small $\e>0$.
}
Indeed, let $s_{\om}=\sign(g(\om))$. It follows from (g.9) that for any $2\leq k\leq n$, $|g_{k}(\om)|\leq c_{\adr{g4}}|g(\om)|^{t_{k}}\leq c_{\adr{g4}}|g(\om)|^{t_{0}}$. Then we have
\ali
{
\frac{\sign(g(\om))g(\e,\om)-|g(\om)|}{\e|g(\om)|^{t_{0}}}=&\sum_{k=1}^{n}\frac{\sign(g(\om))g_{k}(\om)}{|g(\om)|^{t_{0}}}\e^{k-1}\\
\geq &c_{\adr{g9}}-\sum_{k=2}^{n}c_{\adr{g4}}\e^{k-1}>c_{\adr{g9}}-c_{\adr{g4}}\frac{\e}{1-\e}>\frac{c_{\adr{g9}}}{2}>0
}
for any small $\e>0$. This implies that the signature of $\sign(g(\om))g(\e,\om)$ is plus for any small $\e>0$ and therefore the signature of $g(\e,\om)$ equals the signature of $g(\om)$. This also yields $x(\e,\om)=\sign(g(\om))g(\e,\om)-|g(\om)|>0$ and thus the assertion is valid by putting $c_{\adr{AS3}}=c_{\adr{g9}}/2$.
\cla
{
\label{cla:LR_t+(n+1)(1-s1)}
$\|\LR_{(\underline{c}-n(1-t_{1}))\log|g|}1\|_{\infty}<\infty$ and $\|\LR_{(\overline{c}-(n+1)(1-t_{0}))\log|g|}1\|_{\infty}=\infty$.
}
Since $\underline{c}-(1-t_{1})n$ is greater than $\underline{p}$, $P((\underline{c}-(1-t_{1})n)\log |g|)$ is finite and so is $\|\LR_{(\underline{c}-(1-t_{1})n)\log |g|}1\|_{\infty}$. On the other hand, from $\overline{c}-(1-t_{0})(n+1)$ is less than $\underline{p}$, $P((\overline{c}-(1-t_{1})n)\log |g|)$ is infinite and it yields $\|\LR_{(\overline{c}-(1-t_{1})n)\log |g|}1\|_{\infty}=\infty$. Therefore this claim is valid.
\smallskip
\par
We let
\ali
{
E(\e)=\{e\in E\,:\,\inf_{\om\in [e]}|g(\om)|\geq 2c_{\adr{g4}}\e\}.
}
Then we see that $E(\e)$ is an including finite set and $\lim_{\e\to 0}E(\e)=E$. We will use the fact that for any $e\in E(\e)$, $\om\in [e]$ and $0<\e<1/2$
\alil
{
|x(\e,\om)|\leq c_{\adr{g4}}\sum_{k=1}^{n}|g(\om)|^{t_{k}}\e^{k}\leq 2c_{\adr{g4}}\e\leq |g(\om)|\label{eq:xeom<=|gw|}
}
by $|g(\om)|\leq 1$.
\cla
{
There exists $\ti{b}\in E$ such that $\inf_{\om\in [\ti{b}]}\sum_{e\in E(\e)\,:\,t(e)=i(\ti{b})}|g(e\cd\om)|^{\overline{c}-(n+1)(1-t_{0})}\to \infty$ as $\e\to 0$.
}
Choose any large number $M>0$. Since incidence matrix is finitely irreducible, there exists a finite subset $\{b_{1},\dots, b_{N}\}$ of $E$ such that for any $e\in E$, $t(e)=i(b_{k})$ for some $k$. Namely when we put
\ali
{
E_{k}=&\{e\in E\,:\,t(e)=i(b_{k})\}
}
for each $k=1,2,\dots, N$, then $E=\bigcup_{k=1}^{N}E_{k}$ is satisfied.
From $\|\LR_{(\overline{c}-(1-t_{0})(n+1))\log |g|}1\|_{\infty}=+\infty$, there is $\om\in E^{\infty}$ satisfying $\LR_{(\overline{c}-(1-t_{0})(n+1))\log |g|}1(\om)>M$. Moreover, there exists $\e_{0}>0$ such that for any $0<\e<\e_{0}$, $\sum_{e\in E(\e)\,:\,e\cd\om\in E^{\infty}}|g(e\cd\om)|^{\overline{c}-(1-t_{0})(n+1)}>M$. 
We notice that for any $\om^{k}\in [b_{k}]$ $(k=1,2,\dots, N)$
\ali
{
\sum_{e\in E(\e)\,:\,t(e)=i(\om_{0})}|g(e\cd\om)|^{\overline{c}-(1-t_{0})(n+1)}&\leq \sum_{k=1}^{N}\sum_{e\in E_{k}\cap E(\e)\,:\,t(e)=i(\om_{0})}|g(e\cd \om)|^{\overline{c}-(1-t_{0})(n+1)}\\
&\leq \sum_{k=1}^{N}\sum_{e\in E_{k}\cap E(\e)\,:\,t(e)=i(b_{k})}|g(e\cd\om^{k})|^{\overline{c}-(1-t_{0})(n+1)}(1+c_{\adr{cg^q}}\theta)\\
&\qqqqqqqqqquad (\because \text{Lemma } \ref{lem:prop_g^q})
}
Thus
\ali
{
M\leq& \sum_{k=1}^{N}\inf_{\up\in [b_{k}]}\sum_{e\in E(\e)\,:\,t(e)=i(b_{k})}|g(e\cd\up)|^{\overline{c}-(1-t_{0})(n+1)}(1+c_{\adr{cg^q}}\theta).
}
Since $M$ is an arbitrary large number, the right hand side tends to $+\infty$ as $\e\to 0$. Hence the assertion of the claim is valid for some $\ti{b}\in \{b_{k}\}_{k=1}^{N}$.
\cla
{
There exist a sign $\ti{s}\in \{+1,-1\}$ and a constant $c_{\adl{tLRnse>=}}>0$ such that
$\ti{s}\tLR_{n,s}(\e,h)(\om)\geq -c_{\adr{tLRnse>=}}$ uniformly in $s\in I$ and $\om\in E^{\infty}$, and $\inf_{s\in I}\inf_{\om\in [\ti{b}]}\ti{s}\tLR_{n,s}(\e,h)(\om)/\e\to +\infty$.
}
Let $\om\in E^{\infty}$. By virtue of the Taylor expansion for the function $F\,:\,x\mapsto (|g(\om)|+x)^{s}$, we obtain the form
\alil
{
|g(\e,\om)|^{s}=F(x(\e,\om))&=\sum_{l=0}^{n}\binom{s}{l}|g(\om)|^{s-l}x(\e,\om)^{l}+\binom{s}{n+1}(|g(\om)|+\alpha x(\e,\om))^{s-n-1}x(\e,\om)^{n+1}\nonumber\\
&=g_{0,s}+g_{1,s}\e+\cdots+g_{n,s}\e^{n}+\ti{g}_{n,s}(\e,\cd)\e^{n}\label{eq:|geom|^s=}
}
with
\ali
{
\ti{g}_{n,s}(\e,\om)=&\sum_{l=0}^{n}\binom{s}{l}|g(\om)|^{s-l}\sign(g(\om))^{l}\sum_{j=n+1}^{nl}\sum_{j_{1},\cdots ,j_{n}\geq 0:\atop{j_{1}+\cdots+j_{n}=l\atop{j_{1}+2j_{2}+\dots+nj_{n}=j}}}\frac{l!}{j_{1}!\cdots j_{n}!}g_{1}(\om)^{j_{1}}\cdots g_{n}(\om)^{j_{n}}\e^{j-n}\\
&+\binom{s}{n+1}(|g(\om)|+\alpha x(\e,\om))^{s-n-1}\left(\frac{x(\e,\om)}{\e}\right)^{n+1}\e=I_{1}(\e,\om)+I_{2}(\e,\om)
}
for some $\alpha=\alpha(s,\e,\om)\in [0,1]$.
By noting $t_{1}\leq t_{0}\leq t_{k}$ for $k\geq 2$ and by the estimate $|g|^{s-l}|g_{1}|^{j_{1}}\cdots |g_{n}|^{j_{n}}\leq c_{\adr{g4}}^{l}|g|^{s-l+lt_{1}}\leq \max(c_{\adr{g4}},1)^{n}|g|^{s-(1-t_{1})n}$, we get the inequality
\alil{
|I_{1}(\e,\om)|\leq c_{\adr{I_gew^s}}|g(\om)|^{s-(1-t_{1})n}\e\label{eq:Ieom<=}
}
for any $\om\in E^{\infty}$ for some constant $c_{\adl{I_gew^s}}>0$. To estimate $I_{2}(\e,\cd)$, we note that since the integers $1,2,\dots, n$ are not in $I$, the signature of $\binom{s}{n+1}$ for $s\in I$ does not depend on choice $s\in I$.
Then we put $\ti{s}=\sign(\binom{s}{n+1})$ for a $s\in I$. By the same reason as above, $\inf_{s\in I}|\binom{s}{n+1}|=:c_{\adl{inf_binsn+1}}$ is positive.
Remark also that for $e\in E(\e)$
\ali
{
(|g(e\cd\om)|+\alpha x(\e,e\cd\om))^{s-n-1}&\geq 
\case
{
|g(e\cd\om)|^{s-n-1},&s-n-1\geq 0\\
(|g(e\cd\om)|+x(\e,e\cd\om))^{s-n-1},&s-n-1<0\\
}\\
&\geq \case
{
|g(e\cd\om)|^{s-n-1},&s-n-1\geq 0\\
2^{s-n-1}|g(e\cd\om)|^{s-n-1},&s-n-1<0\\
}\\
&\geq 2^{-|s-n-1|}|g(e\cd\om)|^{s-n-1}
}
by (\ref{eq:xeom<=|gw|}). Thus we obtain that for any $s\in I$ and $\om\in E^{\infty}$,
\alil
{
\frac{\ti{s}I_{2}(\e,e\cd\om)}{\e}=&\Big|\binom{s}{n+1}\Big|(|g(e\cd\om)|+\alpha x(\e,e\cd\om))^{s-n-1}\left(\frac{x(\e,e\cd\om)}{\e}\right)^{n+1}\nonumber\\
\geq &
\case
{
c_{\adr{tsIIeew>=}}|g(e\cd\om)|^{\overline{c}-(1-t_{0})(n+1)},&\om\in [e] \text{ for some }e\in E(\e)\\
0,&\text{otherwise}
}\label{eq:IIeom<=}
}
with $c_{\adl{tsIIeew>=}}=c_{\adr{inf_binsn+1}}2^{-\sup_{s\in I}|s-n-1|}c_{\adr{AS3}}^{n+1}>0$.
Consequently, (\ref{eq:Ieom<=}) and (\ref{eq:IIeom<=}) imply that for any $\om\in E^{\infty}$
\ali
{
\frac{\ti{s}\tLR_{n,s}(\e,h)(\om)}{\e}= &\sum_{e\in E\,:\,t(e)=i(\om_{0})}(\ti{s}I_{1}(\e,e\cd\om)+\ti{s}I_{2}(\e,e\cd\om))h(e\cd\om)\frac{1}{\e}\\
\geq
&\case
{ \di-c_{\adr{tLRnse>=}}+c_{\adr{tsIIeew>=}}(\inf_{\up}h(\up))\sum_{e\in E(\e)\,:\,t(e)=i(\om_{0})}|g(e\cd\om)|^{\overline{c}-(n+1)(1-t_{0})},&\om\in[\ti{b}]\\
-c_{\adr{tLRnse>=}},&\text{otherwise}
}
}
by putting $c_{\adr{tLRnse>=}}=c_{\adr{I_gew^s}}\|\LR_{(\underline{c}-(1-t_{1})n)\log|g|}1\|_{\infty}\|h\|_{\infty}$. Thus the assertion is valid.
\cla
{
The assertion of this corollary is valid.
}
Recall the form $\ti{s}_{n}(\e)$ ($n\geq 1)$ of (\ref{eq:tsne=}) in Theorem \ref{th:asympsol_Beq_Mgene}. Then we have
\ali
{
\frac{\ti{s}\ti{s}_{n}(\e)}{\e}=&\frac{1}{-\nu(h\log|g|)}\nu(\e,\frac{\ti{s}\tLR_{n,s(\e)}(\e,h)}{\e})+O(1)\\
\geq&-\frac{-c_{\adr{tLRnse>=}}}{-\nu(h\log|g|)}+
\frac{\nu(\e,[\ti{b}])}{-\nu(h\log|g|)}\frac{\inf_{s\in I}\inf_{\om\in [\ti{b}]}\ti{s}\tLR_{n,s}(\e,h)(\om)}{\e}+O(1)\\
\to&+\infty
}
as $\e\to 0$ by using the above claim in addition to the fact $\nu(\e,[\ti{b}])\to \nu([\ti{b}])>0$. Hence the assertion $|\ti{s}_{n}(\e)|/\e\to +\infty$ is guarantied.
\proe
\subsection{Proof of Theorem \ref{th:asymp_dim}}
For a sake of convenience, we may write the composite map $T_{\om_{0}}(\e,\cd) T_{\om_{1}}(\e,\cd) \cdots T_{\om_{n}}(\e,\cd)$ by $T_{\om_{0}\om_{1}\cdots \om_{n}}(\e,\cd)$. Similarity, $T_{\om_{0}\om_{1}\cdots \om_{n}}$ means $T_{\om_{0}} T_{\om_{1}} \cdots T_{\om_{n}}$. 
Assume that the condition $(G.1)_{n}$ is satisfied. We take open and relative compact subsets $(U_{v})$ of $\R^{D}$ and numbers $r\in (0,1)$ and $r_{0}>0$ such that $U_{v}=\bigcup_{x\in J_{v}}B(x,r_{0})$, $J_{v}\subset U_{v}\subset \overline{U_{v}}\subset O_{v}$ for any $v\in V$, and $\sup_{e\in E}\sup_{x\in U_{t(e)}}\|T_{e}^\p(x)\|\leq r$, where $B(x,r_{0})$ is the open ball with center $x$ and radius $r_{0}$.
We begin with the following fact:
\lem
{\label{lem:induc_asympGDMS}
For any $n\geq 1$, if the conditions $(G.1)_{n}$ and $(G.2)_{n}$ are satisfied, then so are the conditions $(G.1)_{n-1}$ and $(G.2)_{n-1}$.
}
\pros
Assume that the conditions $(G.1)_{n}$ and $(G.2)_{n}$ are satisfied. It suffices to prove that the conditions (iii) and (iv) in $(G.2)_{n-1}$ are fulfilled for $n\geq 2$. Since $\ti{T}_{e,n-1}(\e,\cd)$ has the form $T_{e,n}\e+\ti{T}_{e,n}(\e,x)\e$, convergence $\sup_{e\in E}\sup_{x\in J_{t(e)}}(\|\frac{\partial}{\partial x}\ti{T}_{e,n-1}(\e,x)\|/\|T_{e}^\p(x)\|^{\ti{t}_{1}})\to 0$ is yielded by putting $\ti{t}_{1}=\min(t(n,1),\ti{t}_{0})$. Therefore (iii) is valid for $(G.2)_{n-1}$. To check (iv) in $(G.2)_{n-1}$, we note the forms
\ali
{
p(n-1)=&\max\{\underline{p}-\frac{n-1}{1}(1-t_{1}), \cdots, \underline{p}-\frac{n-1}{n-1}(1-t_{n-1}),\\
&\qqqqqqqquad\underline{p}/t_{1},\dots, \underline{p}/t_{n-1},\underline{p}+1-\ti{s},\underline{p}/\ti{s}\}\\
\ti{s}=&\min\{t_{n-1},\ti{t}_{1},\frac{\ti{t}_{1}}{D}+\frac{D-1}{D}t(1,1),\dots, \frac{\ti{t}_{1}}{D}+\frac{D-1}{D}t(n-1,1)\}.
}
By the definition of $t_{k}$ in (\ref{eq:tk=}), the inequality $t_{n-1} \geq t_{n}$ holds. By the same reason, we have $t(n,1),\cdots, t(1,1)\geq t_{n}$ and therefore $\ti{s}\geq \ti{t}$, where $\ti{t}$ is defined in (\ref{eq:tt=}). Thus we see $p(n-1)\leq p(n)$. Hence $\dim_{H}K/D>p(n)$ implies $\dim_{H}K/D>p(n-1)$.
\proe
\lem
{\label{lem:bdd_Telkx-...}
If the condition $(G.1)_{n}$-(ii) holds, then there exists a constant $c_{\adl{Gbd2}}>0$ such that for any $e\in E$, $x\in J_{t(e)}$, $y\in O_{t(e)}$ with $|x-y|<r_{0}$, $0\leq l\leq n$, and $0\leq k\leq 1+n-l$
\alil
{
\|T_{e,l}^{(k)}(x)-T_{e,l}^{(k)}(y)\|\leq c_{\adr{Gbd2}} \|T_{e}^\p(x)\|^{t(l,k)}|x-y|^{\beta},\quad \|T_{e,l}^{(k)}(y)\|\leq c_{\adr{Gbd2}}\|T_{e}^\p(x)\|^{t(l,k)}.
}
}
\pros
We have
\ali
{
\Big|\frac{T_{e,l}^{(k)}(x)}{\|T_{e}^\p(x)\|^{t(l,k)}}-\frac{T_{e,l}^{(k)}(y)}{\|T_{e}^\p(x)\|^{t(l,k)}}\Big|\leq& \Big|\frac{T_{e,l}^{(k)}(x)}{\|T_{e}^\p(x)\|^{t(l,k)}}-\frac{T_{e,l}^{(k)}(y)}{\|T_{e}^\p(y)\|^{t(l,k)}}\Big|+\Big|\frac{T_{e,l}^{(k)}(y)}{\|T_{e}^\p(y)\|^{t(l,k)}}-\frac{T_{e,l}^{(k)}(y)}{\|T_{e}^\p(x)\|^{t(l,k)}}\Big|\\
\leq&c_{\adr{G2ii_c1}}|x-y|^{\beta}+c_{\adr{G2ii_c2}}\|T_{e}^\p(y)\|^{t(l,k)}\frac{|\|T_{e}^\p(x)\|^{t(l,k)}-\|T_{e}^\p(y)\|^{t(l,k)}|}{\|T_{e}^\p(x)\|^{t(l,k)}\|T_{e}^\p(y)\|^{t(l,k)}}\\
\leq&(c_{\adr{G2ii_c1}}+c_{\adr{G2ii_c2}}(1+c_{\adr{Gbd}}r_{0}^\beta))|x-y|^{\beta}
}
for some $c_{\adl{G2ii_c1}},c_{\adl{G2ii_c2}}>0$ by using the condition $(G.2)_{n}$-(ii) and the condition (v) in GDMS. Therefore the former assertion is fulfilled. The letter assertion follows from the above inequality.
\proe
Let $\pi(\e,\cd)$ be the coding map of $K(\e)$ for $\e>0$.
\lem
{\label{lem:asymp_pi}
Assume that the conditions $(G.1)_{n}$ and $(G.2)_{n}$ are satisfied. Choose any $r_{1}\in (r,1)$. Then there exist functions $\pi_{1},\pi_{2},\dots, \pi_{n}\in F_{r_{1},b}(E^{\infty},\R^{D})$ and $\ti{\pi}(\e,\cd)\in C_{b}(E^{\infty},\R^{D})$ such that
$\pi(\e,\cd)=\pi+\pi_{1}\e+\cdots+\pi_{n}\e^{n}+\ti{\pi}_{n}(\e,\cd)\e^{n}$ and $\|\ti{\pi}_{n}(\e,\cd)\|_{\infty}:=\sup_{\om\in E^{\infty}}|\ti{\pi}_{n}(\e,\om)|\to 0$ as $\e\to 0$.
}
\pros
This assertion mostly follows from the proof of \cite[Lemma 3.1]{T2016}. When we use this proof, we need the boundedness of $|T_{e,k}^{(i)}(y)|$ uniformly in $e\in E$ and $y\in U_{t(e)}$ for each $k$ and $i$. This fact is satisfied by the condition $(G.2)_{n}$-(ii) in particular. Therefore the proof of \cite[Lemma 3.1]{T2016} implies
\ali
{
\pi_{j}(\om)=&\sum_{k=0}^{\infty}T_{\om_{0}\cdots \om_{k-1}}^\p(\pi\si^{k}\om)(R_{j}(\pi\si^{k}\om)),\ \ 
\ti{\pi}_{n}(\e,\om)=\sum_{k=0}^{\infty}T_{\om_{0}\cdots \om_{k-1}}^\p(\pi\si^{k}\om)(\ti{R}_{n}(\e,\pi\si^{k}\om)),
}
where $R_{j}$ and $\ti{R}_{n}(\e,\cd)$ are defined inductively
\ali
{
R_{j}(\om)&=T_{\om_{0},j}(\pi\si\om)+\sum_{0\leq l\leq j-1,1\leq k\leq j-l:\atop{(l,k)\neq (0,1)}}\sum_{i_{1},\dots,i_{k}\geq 1:\atop{i_{1}+\cdots+i_{k}=j-l}}\frac{T_{\om_{0},l}^{(k)}(\pi\si\om)(\pi_{i_{1}}(\si\om),\dots, \pi_{i_{k}}(\si\om))}{k!}\\
\ti{R}_{n}(\e,\om)&=\sum_{0\leq l\leq j-1,1\leq k\leq j-l:\atop{(l,k)\neq (0,1)}}\sum_{i=n-l+1}^{kn}\sum_{1\leq i_{1},\dots, i_{k}\leq n-1:\atop{i_{1}+\cdots+i_{k}=i}}\frac{T_{\om_{0},l}^{(k)}(x)(\pi_{i_{1}}(\si\om),\dots, \pi_{i_{k}}(\si\om))}{k!}\e^{i-n+l}\\
+&\sum_{0\leq l\leq j-1,1\leq k\leq j-l:\atop{(l,k)\neq (0,1)}}\sum_{i=1}^{k}\frac{T_{\om_{0},l}^{(k)}(x)}{k!}(\underbrace{z(\e),\dots, z(\e)}_{i-1},\underbrace{\ti{\pi}_{n-1}(\e,\si\om)}_{i\text{-th}},\underbrace{x(\e)-x,\dots, x(\e)-x}_{n-i})\e^{l-1}\\
+&\sum_{l=0}^{n}L(n-l,T_{\om_{0},l},x(\e),x)\left(\frac{x(\e)-x}{\e}\right)^{n-l}+\ti{T}_{\om_{0},n}(\e,x(\e)),
}
where $z(\e)=\sum_{k=1}^{n-1}\pi_{k}(\si\om)\e^{k}$, $x(\e)=\pi(\e,\si\om)$, $x=\pi\si\om$, and $L(n-l,T_{\om_{0},l},x(\e),x)=\int_{0}^{1}\frac{(1-t)^{n-l-1}}{(n-l-1)!}(T_{\om_{0},l}^{(n-l)}(x+t(x(\e)-x))-T_{\om_{0},l}^{(n-l)}(x))dt$. The facts $\pi_{k}\in F_{r_{1},b}(E^{\infty},\R^{D})$ and $\|\ti{\pi}_{n}(\e,\cd)\|_{\infty}\to 0$ follow from \cite{T2016} again.
\proe
Now we will give the asymptotic expansion of the function $\om\mapsto \det \frac{\partial}{\partial x}T_{\om_{0}}(\e,\pi(\e,\si\om))$. Put
\ali
{
u(k,i)=\min\{t(i_{1},j_{1}+1)+\cdots+t(i_{D},j_{D}+1)\,:\,&0\leq i_{1},\dots, i_{D}\leq k,\ i_{1}+\dots+i_{D}=k\\
&0\leq j_{1},\dots, j_{D}\leq i,\ j_{1}+\cdots+j_{D}=i\}
}
for each $k=1,\dots, n$. Then we see the equation
\alil
{
t_{k}=&\frac{1}{D}\min(\{u(k,0)\}\cup\{u(l,i)\,:\,l=0,\dots, k-1,\ i=1,\dots,k-l\}).\label{eq:sk=_GDMS}
}
\lems
\label{lem:asmyp_e^phe_MOREgene}
Assume that the conditions $(G.1)_{n}$ and $(G.2)_{n}$ are satisfied. Then the functions $g(\e,\om):=\det \frac{\partial}{\partial x}T_{\om_{0}}(\e,\pi(\e,\si\om))$ and $g(\om):=\det T_{\om_{0}}^{\prime}(\pi\si\om)$ satisfy the conditions (g.1)-(g.5).
\leme
\pros
For $e\in E$, $x\in O_{t(e)}$, $0\leq k\leq n$ and $\e>0$, we write $T_{e}(\e,x)=(t_{e,1}(\e,x),\dots, t_{e,D}(\e,x))$, $T_{e,k}(x)=(t_{e,k,1}(x),\dots, t_{e,k,D}(x))$ and $\ti{T}_{e,n}(\e,x)=(\ti{t}_{e,n,1}(\e,x),\dots, \ti{t}_{e,n,D}(\e,x))$, where $T_{e,0}=T_{e}$.
Note the form
\ali
{
\det \frac{\partial}{\partial x}T_{e}(\e,x)=\sum_{\eta}\sgn(\eta)\frac{\partial t_{e,1}(\e,x)}{\partial x_{\eta(1)}}\cdots \frac{\partial t_{e,D}(\e,x)}{\partial x_{\eta(D)}}
}
for $x=(x_{1},x_{2},\dots, x_{D})\in J_{t(e)}$, where $\eta$ is taken over all permutations on $\{1,2,\dots, D\}$ and $\sgn(\eta)$ denotes the sign of $\eta$. We also recall the form
\ali
{
\det \frac{\partial}{\partial x}T_{e}(\e,x)=\det T_{e}^\p(x)+\kappa_{e,1}(x)\e+\cdots+\kappa_{e,n}(x)\e^{n}+\ti{\kappa}_{e,n}(\e,x)\e^{n},
}
where we let
\alil
{
\kappa_{e,k}(x)=&\sum_{\eta}\sgn(\eta)\sum_{0\leq i_{1},\dots, i_{D}\leq k\,:\,\atop{i_{1}+\cdots+i_{D}=k}}\prod_{p=1}^{D}\frac{\partial t_{e,i_{p},p}(x)}{\partial x_{\eta(p)}}\nonumber\\
\ti{\kappa}_{e,n}(\e,x)=&\sum_{\eta}\sgn(\eta)\sum_{i=n+1}^{Dn}\sum_{0\leq i_{1},\dots, i_{D}\leq n\,:\,\atop{i_{1}+\cdots+i_{D}=i}}\prod_{p=1}^{D}\frac{\partial t_{e,i_{p},p}(x)}{\partial x_{\eta(p)}}\e^{i-n}\label{eq:tkenex=}\\
&+\sum_{\eta}\sgn(\eta)\sum_{j=1}^{D}\bigg\{\left(\prod_{q=1}^{j-1}\sum_{l=0}^{n}\frac{\partial t_{e,l,q}(x)}{\partial x_{\eta(q)}}\e^{l}\right)\frac{\partial \ti{t}_{e,n,j}(\e,x)}{\partial x_{\eta(j)}}\left(\prod_{p=j+1}^{D}\frac{\partial t_{e,p}(\e,x)}{\partial x_{\eta(p)}}\right)\bigg\}\nonumber
}
for each $e\in E$, $x\in O_{t(e)}$, $k=1,2,\dots, n$ and $\e>0$ (see \cite[Lemma 3.2]{T2016}). Note that $\kappa_{e,k}$ is of class $C^{n-k+\beta}$ and has the form
\ali
{
\kappa_{e,k}^{(i)}(x)=&\sum_{j_{1},\dots, j_{D}\geq 0\,:\atop{j_{1}+\cdots+j_{D}=i}}\sum_{\eta}\sgn(\eta)\sum_{0\leq i_{1},\dots, i_{D}\leq n\,:\,\atop{i_{1}+\cdots+i_{D}=k}}\prod_{p=1}^{D}\left(\frac{\partial t_{e,i_{p},p}}{\partial x_{\eta(p)}}\right)^{(j_{p})}(x)
}
for each $i=0,1,\dots, n-k$ and $x\in O_{t(e)}$. For $1\leq j\leq D$, we let $z(0)\in \R^{D}$ as $z(0)_{j}=1$ and $z(0)_{j_{0}}=0$ for $j_{0}\neq j$. By using Lemma \ref{lem:bdd_Telkx-...}, we have
\ali
{
\left\|\left(\frac{\partial t_{e,i,q}}{\partial x_{j}}\right)^{(p)}(y)\right\|=&\sup_{z(1),\dots, z_{p}\in \R^{D}\,:\,\atop{|(z(0),z(1),\dots, z(p))|\leq 1}}|\sum_{1\leq i_{0}\cdots i_{p}\leq D}\frac{\partial^{p+1}t_{e,i,q}(y)}{\partial x_{i_{0}}\cdots x_{i_{p}}}(z(0),z(1),\dots, z(p))|\\
\leq&\|t_{e,i,q}^{(p+1)}(y)\|\leq \|T_{e,i}^{(p+1)}(y)\|\leq c_{\adr{Gbd2}}\|T_{e}^\p(x)\|^{t(i,p+1)}
}
and
\ali
{
&\left\|\left(\frac{\partial t_{e,i,q}}{\partial x_{j}}\right)^{(p)}(x)-\left(\frac{\partial t_{e,i,q}}{\partial x_{j}}\right)^{(p)}(y)\right\|\\
=&\sup_{z(1),\dots, z_{p}\in \R^{D}\,:\,\atop{|(z(0),\dots, z(p))|\leq 1}}\left|\sum_{1\leq i_{1},\dots, i_{p}\leq D}\left\{\frac{\partial^{p+1}t_{e,i,q}(x)}{\partial x_{i_{0}}\cdots x_{i_{p}}}(z(0),\dots, z(p))-\frac{\partial^{p+1}t_{e,i,q}(y)}{\partial x_{i_{0}}\cdots x_{i_{p}}}(z(0),\dots, z(p))\right\}\right|\\
\leq&\|t_{e,i,q}^{(p+1)}(x)-t_{e,i,q}^{(p+1)}(y)\|\\
\leq&\|T_{e,i}^{(p+1)}(x)-T_{e,i}^{(p+1)}(y)\|\leq \max\{c_{\adr{Gbd}},c_{\adr{Gbd2}}\}\|T_{e}^\p(x)\|^{t(i,p+1)}|x-y|^{\beta}
}
for $x\in J_{t(e)}$ and $y\in U_{t(e)}$ with $|x-y|<r_{0}$, and for $0\leq p\leq n-k$. Therefore
\alil
{
|\kappa_{e,k}^{(i)}(y)|\leq&\sum_{j_{1},\dots, j_{D}\geq 0\,:\atop{j_{1}+\cdots+j_{D}=i}}\sum_{\eta}\sum_{0\leq i_{1},\dots, i_{D}\leq n\,:\,\atop{i_{1}+\cdots+i_{D}=k}}\prod_{p=1}^{D}\left\|\left(\frac{\partial t_{e,i_{p},p}}{\partial x_{\eta(p)}}\right)^{(j_{p})}(y)\right\|\nonumber\\
\leq&\sum_{j_{1},\dots, j_{D}\geq 0\,:\atop{j_{1}+\cdots+j_{D}=i}}\sum_{\eta}\sum_{0\leq i_{1},\dots, i_{D}\leq k\,:\,\atop{i_{1}+\cdots+i_{D}=k}}(c_{\adr{Gbd2}})^{D}\|T_{e}^\p(x)\|^{t(i_{1},j_{1}+1)+\cdots+t(i_{D},j_{D}+1)}\nonumber\\
\leq&c_{\adr{Gbd6}}\|T_{e}^\p(x)\|^{u(k,i)}\label{eq:k_ek^i(x)<=...}
}
for some constant $c_{\adl{Gbd6}}$. Moreover, we have
\ali
{
&\|\kappa_{e,k}^{(i)}(x)-\kappa_{e,k}^{(i)}(y)\|
\leq c_{\adr{Gbd8}}\|T_{e}^\p(x)\|^{u(k,i)}|x-y|^{\beta}
}
by using Proposition \ref{prop:nk-Lip_ex} for each $e\in E$, $x\in J_{t(e)}$, $y\in U_{t(e)}$ with $|x-y|<r_{0}$ and $i=0,1,\dots, n-k$ for some constant $c_{\adl{Gbd8}}>0$.
On the other hand, the form (\ref{eq:tkenex=}) implies
\ali
{
\|\ti{\kappa}_{e,n}(\e,x)\|\leq& D!(n+1)^{D}c_{\adr{Gbd2}}^{D}\e\|T_{e}^\p(x)\|^{\min_{n+1\leq i\leq Dn}u(i,0)}\\
&+D!D(2c_{\adr{Gbd2}})^{D-1} c_{\adr{Gbd3}}(\e)\sum_{j=1}^{D}\|T_{e}^\p(x)\|^{(j-1)\min\{t(1,1),\cdots,t(n,1)\}+\ti{t}+(D-j)\min\{t(1,1),\cdots,t(n,1),\ti{t}\}}\\
\leq&c_{\adr{Gbd095}}(\e)\|T_{e}^\p(x)\|^{\ti{u}}
}
with
\ali
{
\ti{u}=\min\{u(n+1,0),\dots, u(Dn,0),\ti{t}+(D-1)\min\{t(1,1),\cdots,t(n,1),\ti{t}\}\}
}
for any $e\in E$, $x\in J_{t(e)}$ and small $\e>0$ by putting 
$c_{\adl{Gbd095}}(\e)=D!\max((n+1)^{D}c_{\adr{Gbd2}}^{D}\e, D(2c_{\adr{Gbd2}})^{D-1} c_{\adr{Gbd3}}(\e))$ with $c_{\adr{Gbd095}}(\e)\to 0$.
Consequently, by the asymptotic expansion of the composite functions (see \cite[Proposition 2.3]{T2016}), we obtain
\ali
{
g(\e,\om)=g(\om)+g_{1}(\om)\e+\cdots +g_{n}(\om)\e^{n}+\ti{g}_{n}(\e,\om)\e^{n}
}
with
\alil
{
g_{j}(\om)=&\kappa_{\om_{0},j}(\pi\si\om)+\sum_{l=0}^{j-1}\sum_{k=1}^{j-l}\sum_{i_{1},\dots,i_{k}\geq 1\,:\,\atop{i_{1}+\cdots+i_{k}=j-l}}\frac{\kappa_{\om_{0},l}^{(k)}(x)(\pi_{i_{1}}(\si\om),\dots,\pi_{i_{k}}(\si\om))}{k!}\label{eq:gj=}\\
\ti{g}_{n}(\e,\om)=&\sum_{l=0}^{n-1}\sum_{k=1}^{n-l}\sum_{i=n-l+1}^{kn}\sum_{1\leq i_{1},\dots,i_{k}\leq n\,:\atop{i_{1}+\cdots+i_{k}=i}}\frac{\kappa_{\om_{0},l}^{(k)}(x)(\pi_{i_{1}}(\si\om),\dots,\pi_{i_{k}}(\si\om))}{k!}\e^{i-n+l}\nonumber\\
&+\sum_{l=0}^{n}\sum_{k=1}^{n-l}\sum_{i=1}^{k}\kappa_{\om_{0},l}^{(k)}(x)(\underbrace{z(\e),\dots,z(\e)}_{i-1},\underbrace{\ti{\pi}_{n}(\e,\si\om)}_{i\text{-th}},\underbrace{x(\e)-x,\dots,x(\e)-x}_{k-i})\frac{\e^{l}}{k!}\nonumber\\
&+\sum_{l=0}^{n-1}\int_{0}^{1}\frac{(1-t)^{n-l-1}}{(n-l-1)!}(\kappa_{\om_{0},l}^{(n-l)}(x(\e)+t(x-x(\e)))-\kappa_{\om_{0},l}^{(n-l)}(x))\,dt\left(\frac{x(\e)-x}{\e}\right)^{n-l}\nonumber\\
&+\kappa_{\om_{0},n}(x(\e))-\kappa_{\om_{0},n}(x)+\ti{\kappa}_{\om_{0},n}(\e,x(\e)),\label{eq:tgne=}
}
where $z(\e)=\sum_{k=1}^{n}\pi_{k}(\si\om)\e^{k}$, $x(\e)=\pi(\e,\si\om)$ and $x=\pi\si\om$. Then we see
\ali
{
|g_{j}(\om)|\leq &|\kappa_{\om_{0},j}(\pi\si\om)|+\sum_{l=0}^{j-1}\sum_{k=1}^{j-l}\sum_{i_{1},\dots,i_{k}\geq 1\,:\,\atop{i_{1}+\cdots+i_{k}=j-l}}\frac{\|\kappa_{\om_{0},l}^{(k)}(\pi\si\om)\|\|\pi_{i_{1}}\|_{\infty}\cdots\|\pi_{i_{k}}\|_{\infty}}{k!}\\
\leq&c_{\adr{Gbd6}}\left(\|T_{\om_{0}}^\p(\pi\si\om)\|^{u(j,0)}+\sum_{l=0}^{j-1}\sum_{k=1}^{j-l}\sum_{i_{1},\dots,i_{k}\geq 1\,:\,\atop{i_{1}+\cdots+i_{k}=j-l}}\|T_{\om_{0}}^\p(\pi\si\om)\|^{u(l,k)}\frac{\|\pi_{i_{1}}\|_{\infty}\cdots\|\pi_{i_{k}}\|_{\infty}}{k!}\right)\\
\leq&c_{\adr{Gbd12}}\|T_{\om_{0}}^\p(\pi\si\om)\|^{\min(\{u(j,0)\}\cup\{u(l,k)\,:\,l=0,\dots, j-1,\ k=1,\dots,j-l\})}=c_{\adr{Gbd12}}\|T_{\om_{0}}^\p(\pi\si\om)\|^{D t_{j}}
}
with constant $c_{\adl{Gbd12}}>0$, where $t_{j}$ is defined by (\ref{eq:tk=}).
Moreover, for $\om,\up\in E^{\infty}$ with $\om_{0}=\up_{0}$,
\ali
{
|g_{j}(\om)-g_{j}(\up)|\leq& |\kappa_{\om_{0},j}(\pi\si\om)-\kappa_{\om_{0},j}(\pi\si\up)|\\
&+\sum_{l=0}^{j-1}\sum_{k=1}^{j-l}\sum_{i_{1},\dots,i_{k}\geq 1\,:\,\atop{i_{1}+\cdots+i_{k}=j-l}}\frac{1}{k!}\Big\{\|\kappa_{\om_{0},l}^{(k)}(\pi\si\om)-\kappa_{\om_{0},l}^{(k)}(\pi\si\up)\|\|\pi_{i_{1}}\|_{\infty}\cdots\|\pi_{i_{k}}\|_{\infty}+\\
&+\sum_{q=1}^{k}\|\kappa_{\om_{0},l}^{(k)}(\pi\si\om)\|_{\infty}\||\pi_{i_{q}}(\si\om)-\pi_{i_{q}}(\si\up)|\prod_{1\leq u\leq k\,:\,k\neq q}\|\pi_{i_{u}}\|_{\infty}\Big\}\\
\leq& c_{\adr{Gbd8}}\|T_{\om_{0}}^\p(\pi\si\om)\|^{u(j,0)}[\pi]_{r^\beta}r^{-\beta}d_{r^\beta}(\om,\up)+\\
&+\sum_{l=0}^{j-1}\sum_{k=1}^{j-l}\sum_{i_{1},\dots,i_{k}\geq 1\,:\,\atop{i_{1}+\cdots+i_{k}=j-l}}\frac{1}{k!}\Big\{c_{\adr{Gbd8}}\|T_{\om_{0}}^\p(\pi\si\om)\|^{u(l,k)}[\pi]_{r^\beta}r^{-\beta}\|\pi_{i_{1}}\|_{\infty}\cdots\|\pi_{i_{k}}\|_{\infty}\\
&+\sum_{q=1}^{k}c_{\adr{Gbd6}}\|T_{\om_{0}}^\p(\pi\si\om)\|^{u(l,k)}[\pi_{i_{q}}]_{r_{1}}r_{1}^{-1}\prod_{1\leq u\leq k\,:\,k\neq q}\|\pi_{i_{u}}\|_{\infty}\Big\}d_{\theta}(\om,\up)\\
=&c_{\adr{Gbd9}}|g(\om)|^{t_{j}}d_{\theta}(\om,\up).
}
by putting $r_{1}\in (r^\beta,\theta)$ for some constant $c_{\adl{Gbd9}}>0$. On the other hand, it follows from the definition of the remainder $\ti{g}_{n}(\e,\om)$ addition to the condition (v)  in GDMS's definition that
\ali
{
&|\ti{g}_{n}(\e,\om)|\\
\leq&c_{\adr{Gbd6}}\sum_{l=0}^{n}\sum_{k=1}^{n-l}\left(\e\frac{(\sum_{j=1}^{n}\|\pi_{j}\|_{\infty})^{k}}{k!}+\|\ti{\pi}_{n}(\e,\cd)\|_{\infty}\sum_{i=1}^{k}\|z(\e)\|^{i-1}\|x(\e)-x\|^{k-i}\right)\|T^\p_{\om_{0}}(\pi\si\om)\|^{u(l,k)}\\
&+c_{\adr{Gbd6}}\sum_{l=0}^{n-1}\frac{n-l}{(n-l+1)!}\left\|\frac{\pi(\e,\cd)-\pi}{\e}\right\|_{\infty}^{n-l}\|\pi(\e,\cd)-\pi\|_{\infty}\|T^\p_{\om_{0}}(\pi\si\om)\|^{u(l,n-l)}\\
&+c_{\adr{Gbd8}}\|\pi(\e,\cd)-\pi\|^{\beta}\|T^\p_{\om_{0}}(\pi\si\om)\|^{u(n,0)}+c_{\adr{Gbd095}}(\e)\|T_{e}^\p(x(\e))\|^{\ti{u}}\\
\leq &c_{\adr{Gbd10}}(\e)\|T_{\om_{0}}^\p(\pi\si\om)\|^{\min\{D t_{n},\ti{u}\}}=c_{\adr{Gbd10}}(\e)|g(\om)|^{\ti{t}}
}
for any $\om\in E^{\infty}$ for some number $c_{\adl{Gbd10}}(\e)$ with $\lim_{\e\to 0}c_{\adr{Gbd10}}(\e)=0$, where the last inequality uses the fact $\|T_{e}^\p(x(\e))\|^{\ti{u}}\leq c_{\adr{Gbd095}}(\e)(1+c_{\adr{Gbd}}\|\pi(\e,\cd)-\pi\|_{\infty})^{\ti{u}}\|T_{\om_{0}}^\p(\pi\si\om)\|^{\ti{u}}$. 
Hence the proof is complete.
\proe
\noindent
({\it Proof of Theorem \ref{th:asymp_dim}}). This theorem follows immediately from Theorem \ref{th:asympsol_Beq_Mgene} and Lemma \ref{lem:asmyp_e^phe_MOREgene} by putting $\psi(\e,\cd)\equiv 1$.
\subsection{Proof of Theorem \ref{th:AP_NPL_MOREgene}}
\pros
Choose any $r_{1}\in (r,1)$.
We can give the asymptotic expansion of the coding map $\pi(\e,\cd)$ as well as Lemma \ref{lem:asymp_pi}:
\ali
{
\pi(\e,\cd)=\pi+\pi_{1}\e+\cdots+\pi_{n}\e^{n}+\ti{\pi}_{n}(\e,\cd)\e^{n}
}
with $\pi_{1},\dots, \pi_{n}\in F_{r_{1},b}(E^{\infty},\R^{D})$ and $\|\ti{\pi}_{n}(\e,\cd)\|_{\infty}\to 0$. By the form $g(\e,\om)=DT_{\om_{0}}(\e,\pi(\e,\si\om))$, we have
\ali
{
g(\e,\om)=g(\om)+g_{1}\e+\cdots+g_{n}\e^{n}+\ti{g}_{n}(\e,\cd)\e^{n}
}
with $g(\om)=DT_{\om_{0}}(\pi \si\om)$,
\ali
{
&g_{k}(\om)=S_{\om_{0},k}(\pi\si\om)+\sum_{l=0}^{k-1}\sum_{j=1}^{k-l}\sum_{i_{1}\cdots i_{j}\geq 1:\atop{i_{1}+\cdots+i_{j}=k-l}}\frac{1}{j!}S_{\om_{0},l}^{(j)}(\pi\si\om)(\pi_{i_{1}}(\si\om),\cdots,\pi_{i_{j}}(\si\om))
}
and the remainder $\ti{g}_{n}(\e,\cd)$ given by replacing $\kappa_{\om_{0},l}$ in (\ref{eq:tgne=}) with $S_{\om_{0},l}$.
By a similar argument of the proof of Lemma \ref{lem:asmyp_e^phe_MOREgene}, we obtain the assertion.
\proe
\subsection{Proof of Theorem \ref{th:ex_5}}
\pros
(1) When $a\geq 5$, the number $p(n)$ becomes zero. Therefore $\dim_{H}K>p(n)$ is satisfied whenever $n\geq 0$. Thus Theorem \ref{th:asymp_dim} implies $s(\e)$ has asymptotic expansion with any order $n\geq 0$. The coefficients of $s(\e)$ are calculated as follows. Recall the form (\ref{eq:tn=...}) of $s_{k}$. We see $h\equiv 1$ and $\nu([e])=1/2^{e}$ for $e\in E$. We have
\ali
{
\NR_{u}h=&\sum_{0\leq v\leq u,0\leq q\leq u-v\,:\,\atop{(v,q)\neq (0,1)}}s_{q,u-v}\ZR_{v,q,s(0)}\\
=&\sum_{0\leq v\leq u,0\leq q\leq u-v\,:\,\atop{(v,q)\neq (0,1)}}\sum_{j=0}^{\min(v,q)}s_{q,u-v}\frac{a_{v,j,s(0)}}{(q-j)!}\LR_{(s-(1-s_{0})v)\ph}(\left(\log \left(\frac{1}{5^{\om_{0}}}\right)\right)^{q-j}h)\\
=&\sum_{0\leq v\leq u,0\leq q\leq u-v\,:\,\atop{(v,q)\neq (0,1)}}\sum_{j=0}^{\min(v,q)}s_{q,u-v}\frac{a_{v,j,s(0)}}{(q-j)!}(-1)^{q-j}(\log 5)^{q-j}\sum_{e=1}^{\infty}e^{q-j}\left(\frac{5^{v}}{2 a^{v}}\right)^{e}.
}
In particular, this is a constant function. Moreover,
\ali
{
\nu(h\log|g|)=&\sum_{e=1}^{\infty}\frac{\log(1/5^e)}{2^{e}}=-\log 5\sum_{e=1}^{\infty}\frac{e}{2^{e}}=-2\log 5\\
\ZR_{0,1,s(0)}h(\om)=&\sum_{e=1}^{\infty}\frac{1}{2^{e}}\log\left(\frac{1}{5^e}\right)=-2\log 5.
}
By $\nu(\e,1)=\nu(1)=1$, we obtain $\nu_{i}(\NR_{k-i}h)=\nu_{i}(\ZR_{0,1,s(0)}h)=0$ for $1\leq i\leq k-1$. Thus we get the equation
\ali
{
s_{k}=&\frac{-1}{\nu(h\log|g|)}\left(\sum_{i=1}^{k-1}\nu_{i}(\ZR_{0,1,s(0)}h)s_{k-i}+\sum_{i=0}^{k-1}\nu_{i}(\NR_{k-i}h)\right)=\frac{1}{2\log 5}\nu(\NR_{k}1).
}
This yields the form (\ref{eq:sk=}) of $s_{k}$.
\smallskip
\\
(2) Assume $1<a<5$. Put $a_{0}=\log 5/\log(5/a)$. We show some claims below:
\ncla
{\label{cla:tseasympLRnse}
For each $n\geq 0$, if $\e=o(\ti{\LR}_{n,s(\e)}(\e,h))$ then we have
\alil
{
\ti{s}_{n}(\e)\asymp \ti{\LR}_{n,s(\e)}(\e,h).\label{eq:tsneasymp...}
}
}
Indeed, note the Ruelle operator of $\ph(\e,\om)=\log\|\frac{\partial}{\partial x}T_{\om_{0}}(\e,\pi(\e,\si\om))\|$ has the form $\LR_{s(\e)\phe}f(\om)=\sum_{e=1}^{\infty}\left(1/5^{e})+(\e/a^{e})\right)^{s(\e)}$.
Here, $\ti{\LR}_{n,s(\e)}(\e,\cd)$ is given as follows: By applying Taylor theorem to the function $x\mapsto (1/5^e+x)^{s(\e)}$, we have the expansion
\alil
{
\left(\frac{1}{5^{e}}+\frac{\e}{a^{e}}\right)^{s(\e)}=&\sum_{k=0}^{n}\binom{s(\e)}{k}\left(\frac{1}{5^{e}}\right)^{s(\e)-k}\left(\frac{\e}{a^{e}}\right)^{k}+\binom{s(\e)}{n+1}\left(\frac{1}{5^{e}}+\alpha \frac{\e}{a^{e}}\right)^{s(\e)-n-1}\left(\frac{\e}{a^{e}}\right)^{n+1}\label{eq:1/5e+e/ae)se=...}
}
for each $e\in E$ and $\e>0$ for some $\alpha=\alpha(e,n+1,\e,a)\in [0,1]$. Therefore
\ali
{
\ti{\LR}_{n,s(\e)}(\e,f)(\om)=&\binom{s(\e)}{n+1}\sum_{e\in E}\left(\frac{1}{5^{e}}+\alpha \frac{\e}{a^{e}}\right)^{s(\e)-n-1}\left(\frac{\e}{a^{e}}\right)^{n+1}\frac{1}{\e^{n}}f(e\cd\om).
}
This implies that $\tLR_{n,s(\e)}(\e,h)$ is a constant function. We obtain the claim by (\ref{eq:tsne=}) in Theorem \ref{th:asympsol_Beq_Mgene}.
\cla
{
The inequality $n<a_{0}s(0)\leq n+1$ is satisfied, where $s(0)=\log 2/\log 5$.
}
Indeed, the assumption $a\leq 5/2^{1/(n+1)}$ implies
\ali
{
\frac{\log a}{\log 5}&\leq 1-\frac{1}{n+1}s(0),\qquad
\frac{\log 5}{\log (5/a)}\leq \frac{n+1}{s(0)}
}
and therefore $s(0)a_{0}\leq n+1$.
By a similar argument above, if $n\geq 1$ then $5/2^{1/n}<a$ implies $s(0)a_{0}>n$.
Thus the claim is valid for any $n\geq 0$.
\medskip
\par
We write $\tLR_{n,s(\e)}(\e,1)(\om)\e^{n}=\binom{s(\e)}{n+1}\sum_{e\in E}R_{e}(\e)$ with $R_{e}(\e)=\left(1/5^{e}+\alpha \e/a^{e}\right)^{s(\e)-n-1}\left(\e/a^{e}\right)^{n+1}$. Now we will prove that $\sum_{e\in E}R_{e}(\e)\asymp \e^{a_{0}s(\e)}$.
Note that
\ali
{
\frac{1}{5^{e}}>\frac{\e}{a^{e}} \iff& e< \frac{\log \e}{\log (a/5)}=:a_{1}(\e) \iff e< \lceil a_{1}(\e) \rceil =:a_{2}(\e),
}
where the notation $\lceil\ \rceil$ means round up to the nearest integer. 
Recall the notation $M(n+1,s(\e))$ in Proposition \ref{prop:lower_int_biom} replacing $n:=n+1$, $a:=1/5^e$, $x:=\e/a^e$ and $s:=s(\e)$. Since $\e\mapsto M(n+1,s(\e))$ and $\e\mapsto L(n+1,s(\e))$ are continuous, there exists $\e_{0}>0$ such that $M(n+1,s(\e))\geq M(n+1,s(0))/2$ and $L(n+1,s(\e))\geq L(n+1,s(0))/2$ for any $0<\e<\e_{0}$. Put $a_{3}=\lceil\frac{\log (M(n+1,s(0))/2)}{\log(a/5)}\rceil$. 
We decompose $\sum_{e}R_{e}(\e)$ into
\ali
{
\sum_{e}R_{e}(\e)=\sum_{e=1}^{a_{2}(\e)-1}R_{e}(\e)+\sum_{e=a_{2}(\e)}^{a_{2}(\e)+a_{3}-1}R_{e}(\e)+\sum_{e=a_{2}(\e)+a_{3}}^{\infty}R_{e}(\e)=I_{1}(\e)+I_{2}(\e)+I_{3}(\e).
}
\cla
{\label{cla:limsupI+II/e^as}
Let $1<a<5$ and $n\geq 0$ the largest integer satisfying $a\leq 5/2^{1/(n+1)}$. Then $\limsup_{\e\to 0}(I_{1}(\e)+I_{2}(\e))/(-\e^{n+1}\log \e)<+\infty$ if $a=5/2^{1/(n+1)}$ for some $n\geq 0$, and $\limsup_{\e\to 0}(I_{1}(\e)+I_{2}(\e))/\e^{a_{0}s(0)}<+\infty$ otherwise.
}
Indeed, we have
\ali
{
I_{1}(\e)+I_{2}(\e)\leq& \sum_{e=1}^{a_{2}(\e)+a_{3}-1}\left(\frac{1}{5^{e}}\right)^{s(\e)-n-1}\left(\frac{\e}{a^{e}}\right)^{n+1}\\
\leq &\sum_{e=1}^{a_{2}(\e)+a_{3}-1}\left(\frac{1}{5^{e}}\right)^{s(0)-n-1}\left(\frac{\e}{a^{e}}\right)^{n+1}\quad (\because s(\e)\geq s(0))\\
=&\e^{n+1}\sum_{e=1}^{a_{2}(\e)+a_{3}-1}\left(\frac{5^{n+1}}{2 a^{n+1}}\right)^{e}\\
\leq&
\case
{
\frac{-\log \e}{-\log(a/5)}\e^{n+1}+a_{3}\e^{n+1},&\text{ if }a=5/2^{1/(n+1)} \text{ for some }n\geq 0\\
\e^{n+1}\frac{\left(\frac{5^{n+1}}{2a^{n+1}}\right)^{a_{2}(\e)+a_{3}}-\frac{5^{n+1}}{2a^{n+1}}}{\frac{5^{n+1}}{2a^{n+1}}-1},&\text{otherwise}.
}
}
To show the assertion, it is sufficient to check that $\e^{n+1}(\frac{5^{n+1}}{2a^{n+1}})^{a_{2}(\e)}\leq c\e^{a_{0}s(0)}$ for some $c>0$. This is implied by the facts $a_{2}(\e)\leq a_{1}(\e)-1$, $(1/2)^{a_{1}(\e)}=\e^{a_{0}s(0)}$ and $(5/2)^{(n+1)a_{1}(\e)}=\e^{-n-1}$, and by putting $c=((2a^{n+1})/5^{n+1})$.
This yields the assertion of the claim.
\cla
{
\label{cla:limsupI3/e^a0se}
Let $1<a<5$ and $n\geq 0$ the largest integer satisfying $a\leq 5/2^{1/(n+1)}$. Then $\limsup_{\e\to 0}I_{3}(\e)/\e^{a_{0}s(0)}<+\infty$.
}
To show this, we will apply (\ref{eq:1/5e+e/ae)se=...}) to Appendix \ref{sec:interme_binom} taking $a:=1/5^e$, $x:=\e/a^e$ and $s:=s(\e)$.
we note that $a_{3}$ satisfies $(a/5)^{a_{3}}\leq M(n+1,s(0))/2\leq M(n+1,s(\e))$ and $a_{2}(\e)$ fills $(a/5)^{a_{2}(\e)}/\e\leq 1$. Therefore, $(a/5)^{a_{2}(\e)+a_{3}}/\e\leq M(n+1,s(\e))$ is satisfied. Moreover, it follows from Proposition \ref{prop:lower_int_biom} that
\ali
{
e\geq a_{2}(\e)+a_{3} \Rightarrow \frac{1}{\e}\left(\frac{a}{5}\right)^{e}\leq M(n+1,s(\e)) \Rightarrow 
\case
{
\alpha \geq L(1,s(\e))& n=0\\
\alpha \geq L(n+1,s(\e))\left(\frac{1}{\e}\left(\frac{a}{5}\right)^{e}\right)^{\frac{n-s(\e)}{n+1-s(\e)}}& n\geq 1.\\
}
}
Note also that $L(n+1,s(\e))\geq L(n+1,s(0))/2=:L(n+1)>0$ for $0<\e<\e_{0}$.
In the case when $n=0$, we have
\ali
{
I_{3}(\e)\leq& \sum_{e=a_{2}(\e)+a_{3}}^{\infty}\left(\frac{1}{5^{e}}+L(1)\frac{\e}{a^{e}}\right)^{s(\e)-1}\frac{\e}{a^{e}}\\
\leq &\frac{1}{L(1)}\sum_{e=a_{2}(\e)+a_{3}}^{\infty}\left((1+L(1))\frac{\e}{a^{e}}\right)^{s(0)}\quad (\because \frac{1}{5^{e}}\leq \frac{\e}{a^{e}} \text{ and }s(\e)>s(0))\\
=&\frac{(1+L(1))^{s(0)}}{L(1)(1-1/a)}\e^{s(0)}\left(\frac{1}{a^{s(0)}}\right)^{a_{2}(\e)+a_{3}}\\
\leq&\frac{(1+L(1))^{s(0)}}{L(1)(1-1/a)}\left(\frac{1}{a^{s(0)}}\right)^{a_{3}}\e^{a_{0}s(0)}\quad (\because a^{a_{2}(\e)}\geq a^{a_{1}(\e)}=\e^{\log_{a/5}{a}}=\e^{-a_{0}+1}).
}
Thus the assertion of the claim holds under the case $n=0$.

In the case when $n\geq 1$, we obtain
\ali
{
I_{3}(\e)\leq& \sum_{e=a_{2}(\e)+a_{3}}^{\infty}\left(\frac{1}{5^{e}}+L(n+1)\left(\frac{1}{\e}\left(\frac{a}{5}\right)^{e}\right)^{\frac{n-s(\e)}{n+1-s(\e)}}\frac{\e}{a^{e}}\right)^{s(\e)-n-1}\left(\frac{\e}{a^{e}}\right)^{n+1}\\
\leq&\sum_{e=a_{2}(\e)+a_{3}}^{\infty}\left(\frac{\e}{a^{e}}\right)^{s(\e)}\left(\frac{1}{\e}\left(\frac{a}{5}\right)^{e}+L(n+1)\left(\frac{1}{\e}\left(\frac{a}{5}\right)^{e}\right)^{\frac{n-s(\e)}{n+1-s(\e)}}\right)^{s(\e)-n-1}\\
\leq&L(n+1)^{s(\e)-n-1}\sum_{e=a_{2}(\e)+a_{3}}^{\infty}\left(\frac{\e}{a^{e}}\right)^{s(\e)}\left(\frac{1}{\e}\left(\frac{a}{5}\right)^{e}\right)^{-n+s(\e)}\quad(\because s(\e)-n-1<0)\\
=&L(n+1)^{s(\e)-n-1}\e^{n}\sum_{e=a_{2}(\e)+a_{3}}^{\infty}\left(\frac{1}{5^{s(\e)}}\left(\frac{5}{a}\right)^{n}\right)^{e}\\
\leq& L(n+1)^{s(0)-n-1}\e^{n}\frac{\left(\frac{1}{2}\left(\frac{5}{a}\right)^{n}\right)^{a_{2}(\e)+a_{3}}}{1-\frac{1}{2}\left(\frac{5}{a}\right)^{n}},\quad (\because s(\e)>s(0) \text{ and }5^{s(0)}=2)
}
where in the last expression, we remark $(5/a)^{n}/2<1$ by the definition of $n$. In the last expression, we notice the estimate
$((5/a)^{n}/2)^{a_{2}(\e)}\leq ((5/a)^{n}/2)^{a_{1}(\e)}=\e^{a_{0}s(0)}\e^{-n}$.
Thus we obtain the assertion of Claim \ref{cla:limsupI3/e^a0se}.
\cla
{\label{cla:liminfI1/e^a0se}
Let $1<a<5$ and $n\geq 0$ the largest integer satisfying $a\leq 5/2^{1/(n+1)}$. Then $\liminf_{\e\to 0}I_{1}(\e)/(-\e^{n+1}\log \e)>0$ if $a=5/2^{1/(n+1)}$ for some $n\geq 0$, and $\liminf_{\e\to 0}I_{1}(\e)/\e^{a_{0}s(0)}>0$ otherwise.
}

By virtue of Claim \ref{cla:tseasympLRnse}-\ref{cla:limsupI3/e^a0se}, $s(\e)=s(0)+t(\e)\e^{a_{1}}$ and $t(\e)=O(1)$ are satisfied with $a_{1}:=a_{0}s(0)-\eta$ for any small $\eta>0$.
Then we have that for each $n\geq 0$,
\ali
{
I_{1}(\e)=&\sum_{e=1}^{a_{2}(\e)-1}\left(\frac{1}{5^{e}}+\alpha \frac{\e}{a^{e}}\right)^{s(\e)-n-1}\left(\frac{\e}{a^{e}}\right)^{n+1}\\
\geq& \sum_{e=1}^{a_{2}(\e)-1}\left(\frac{1}{5^{e}}+ \frac{\e}{a^{e}}\right)^{s(\e)-n-1}\left(\frac{\e}{a^{e}}\right)^{n+1}\\
=&\sum_{e=1}^{a_{2}(\e)-1}\frac{\left(\frac{1}{5^{e}}+ \frac{\e}{a^{e}}\right)^{s(0)+t(\e)\e^{a_{1}}}}{\left(\frac{1}{\e}\left(\frac{a}{5}\right)^{e}+ 1\right)^{n+1}}\\
\geq &\sum_{e=1}^{a_{2}(\e)-1}\frac{\left(\frac{1}{5^{e}}\right)^{s(0)}}{\left(\frac{2}{\e}\left(\frac{a}{5}\right)^{e}\right)^{n+1}}\left(\frac{\e}{a^{e}}\right)^{t(\e)\e^{a_{1}}}\quad (\because t(\e)>0)\\
\geq&
\left(\frac{\e}{a^{a_{1}(\e)}}\right)^{t(\e)\e^{a_{1}}}\frac{\e^{n+1}}{2^{n+1}}\case
{
a_{1}(\e)-1,&\text{ if }a=5/2^{1/(n+1)}\\
\frac{\left(\frac{5^{n+1}}{2a^{n+1}}\right)^{a_{1}(\e)}-\frac{5^{n+1}}{2a^{n+1}}}{\frac{5^{n+1}}{2a^{n+1}}-1},&\text{ if }a<5/2^{1/(n+1)}
}\quad (\because a_{2}(\e)\geq a_{1}(\e)).
}
Here we notice that $a_{2}(\e)\geq a_{1}(\e)=\log\e/\log(a/5)$ and $\left(5^{n+1}/(2a^{n})\right)^{a_{1}(\e)}=\e^{a_{0}s(0)}\e^{-n-1}$.
Remark that
\ali
{
\left(\frac{\e}{a^{a_{1}(\e)}}\right)^{t(\e)\e^{a_{1}}}
= &\exp(t(\e)\e^{a_{1}}\log \e)\exp(-t(\e)\frac{\log a}{\log(a/5)}\e^{a_{1}}\log \e)\to \exp(0)\exp(0)=1
}
as $\e\to 0$. 
Thus the assertion of Claim \ref{cla:liminfI1/e^a0se} is yielded.
\smallskip
\par
By the claims \ref{cla:limsupI+II/e^as}-\ref{cla:liminfI1/e^a0se} and the fact $\sum_{e\in E}R_{e}(\e)\geq I_{1}(\e)$, we obtain $\sum_{e\in E}R_{e}(\e)\asymp -\e^{n+1}\log \e$ if $a=5/2^{1/(n+1)}$, and 
$\asymp \e^{a_{0}s(0)}$ otherwise. Thus so be $\tLR_{n,s(\e)}(\e,h)$.
 Hence the proof follows from Claim \ref{cla:tseasympLRnse}.
\proe
\subsection{Proof of Proposition \ref{prop:IIFS_2-asymp}}
\pros
We take $t_{1}$ and $t_{2}$ so that
\ali
{
t_{1}=&\sup\{t\in (0,1]\,:\,\sup_{\om\in E^{\infty}}\frac{|g_{1}(\om)|}{|g(\om)|^{t}}=\sup_{e\in E}(\frac{5^{t}}{4})^{e}<+\infty\} \iff t_{1}=\frac{\log 4}{\log 5}\\
t_{2}=&\sup\{t\in (0,1]\,:\,\sup_{\om\in E^{\infty}}\frac{|g_{1}(\om)|}{|g(\om)|^{t}}=\sup_{e\in E}(\frac{5^{t}}{3})^{e}<+\infty\} \iff t_{2}=\frac{\log 3}{\log 5}.
}
In view of Theorem \ref{th:asympsol_Beq_Mgene}, the number $p(n)$ is defined by
$p(n)=\max(n(1-t_{1}),\ n(1-t_{2})/2)$
by noting $\underline{p}=0$ and $\ti{t}=0$. Therefore, $p(n)=n\max(1-\log 4/\log 5,\ (1-\log 3/\log 5)/2)=n(1-\log 3/\log 5)/2$ by $1-\log 4/\log 5=0.1386\cdots$ and $(1-\log 3/\log 5)/2=0.1586\cdots$. By virtue of Corollary \ref{cor:asympsol_Beq_Mgene}, when $s(0)>p(n)$, the dimension $\dim_{H}K(\e)$ has an asymptotic expansion with order $n$ at $\e=0$. We see
\ali
{
s(0)=\frac{\log 2}{\log 5}>n\frac{1-\log 3/\log 5}{2} &\iff n<\frac{\log 2/\log 5}{(1-\log 3/\log 5)/2}=\frac{2\log 2}{\log 5-\log 3}=2.713\cdots.\\
&\iff n\leq 2.
}
Hence $\dim_{H}K(\e)$ has at least the $2$-asymptotic expansion at $\e=0$.
\proe
\appendix
\section{Thermodynamic formalism and Ruelle operators}\label{sec:thermo}
In this section, we will present useful results for the proof of the main theorem. We will recall the notion of thermodynamic formalism and some facts of Ruelle transfer operators which were manly introduced by \cite{MU}.

We use the notation defined in Section \ref{sec:intro}.
The incidence matrix $A$ of the graph $G$ is called {\it finitely primitive} if there exist an integer $n\geq 1$ and a finite subset $F$ of $E^{n}$ such that for any $e,e^\p\in E$, $ewe^\p$ is a path on the graph $G$ for some $w\in F$. Note that $A$ is finitely primitive if and only if $(E^{\infty},\si)$ is topologically mixing and $A$ has the BIP property. Then it is stronger than finitely irreducible.
A function $\psi\,:\,E^{\infty}\to \R$ is {\it acceptable} if there exists a constant $c_{\adl{acce}}\geq 1$ such that for any $e\in E$ and $\om,\up \in [e]$, $e^{\psi(\om)-\psi(\up)}\leq c_{\adr{acce}}$ (\cite{MU} for the terminology).
For real-valued function $\psi$ on $E^{\infty}$, the topological pressure $P(\psi)$ of $\psi$ is given by
\alil
{\label{eq:toppres}
P(\psi)=\lim_{n\to \infty}\frac{1}{n}\log \sum_{w\in E^{n}\,:\,[w]\neq \emptyset}\exp(\sup_{\om\in [w]}\sum_{k=0}^{n-1}\psi(\si^{k}\om))
}
formally. If $\psi$ is acceptable, then $P(\psi)$ exists in $[-\infty,+\infty]$ (see \cite{MU}).
We mainly consider the pressure function $t\mapsto P(t\psi)\in [-\infty,+\infty]$ with $\sup_{\om\in E^{\infty}}\psi(\om)<0$. In this case, it is basic fact that the pressure function is strictly monotone decreasing,  convex as limit of convex function, and continuous. In particular, $\lim_{t\to +\infty}P(t\psi)=-\infty$ holds.

For a real-valued function $\psi$ on $E^{\infty}$, the Ruelle operator $\LR_{\psi}$ associated to $\psi$ is defined by
\ali
{
\LR_{\psi} f(\om)=\sum_{e\in E\,:\,t(e)=i(\om_{0})}e^{\psi(e\cd \om)}f(e\cd\om)
}
if this series converges in $\C$ for a complex-valued function $f$ on $E^{\infty}$ and for $\om\in E^{\infty}$. Here $e\cd \om$ is the concatenation of $e$ and $\om$, i.e. $e\cd\om=e\om_{0}\om_{1}\cdots$. It is known that if the incidence matrix is finitely irreducible and $\psi$ is in $F_{\theta}(E^{\infty},\R)$ with finite topological pressure, then $\LR_{\psi}$ becomes a bounded linear operator both on the Banach spaces $F_{\theta,b}(E^{\infty})$ and $C_{b}(E^{\infty})$. We begin with the following proposition.
\prop
{\label{prop:finitepres_finiteRuelle}
Let $G=(V,E,i(\cd),t(\cd))$ be a directed multigraph such that the incidence matrix of $E^{\infty}$ is finitely irreducible. Take an acceptable function $\psi\,:\,F_{\theta}(E^{\infty})\to \R$. Then $P(\psi)<\infty$ if and only if $\|\LR_{\psi}1\|_{\infty}<\infty$.
}
\pros
It is known in \cite[Proposition 2.1.9]{MU} that $Z:=\sum_{e\in E}\exp(\sup_{\om\in [e]}\psi(\om))<\infty$ if and only if $P(\psi)<\infty$. Then we will show $Z<\infty$ if and only if $\|\LR_{\psi}1\|_{\infty}<\infty$. Noting $\|\LR_{\psi}\|_{\infty}\leq Z$, it is sufficient to show that $Z$ is finite if $\|\LR_{\psi}1\|_{\infty}<\infty$. Since the BIP property is satisfied, there exists a finite set $\{b_{1},b_{2},\dots, b_{N}\}$ such that for any $e\in E$, there exists $1\leq i(e)\leq N$ such that $t(e)=i(b_{i(e)})$. Choose any $\up^{i(e)}\in [b_{i(e)}]$ for each $e\in E$. We have
\ali
{
Z\leq c_{\adr{acce}}\sum_{e\in E}e^{\psi(\up^{i(e)})}
= c_{\adr{acce}}\sum_{i=1}^{N}\sum_{e\in E\,:\,i(e)=i}e^{\psi(\up^{i})}
\leq c_{\adr{acce}}\sum_{i=1}^{N}\LR_{\psi}(\up^{i})\leq c_{\adr{acce}}N\|\LR_{\psi}\|_{\infty}<+\infty.
}
Hence the assertion is valid.
\proe
A Borel probability measure $\mu$ on $E^{\infty}$ is said to be a {\it Gibbs measure} of the potential $\psi$ if there exist constants $c\geq 1$ and $P\in \R$ such that for any $\om\in E^{\infty}$ and $n\geq 1$
\alil
{
c^{-1}\leq \frac{\mu(\{\up\in E^{\infty}\,:\,\up_{i}=\om_{i},\ 0\leq i<n\})}{\exp(-nP+\sum_{k=0}^{n-1}\psi(\si^{k}\om))}\leq c.\label{eq:Gibbs}
}
Recall the notation $\LR(\XR)$ which is the set of all bounded linear operators acting on a norm space $\XR$.
The following is a version of Ruelle-Perron-Frobenius Theorem:
\thm
{\label{th:exGibbs}
Let $G=(V,E,i(\cd),t(\cd))$ be a directed multigraph such that the incidence matrix of $E^{\infty}$ is finitely irreducible. Assume that $\psi\in F_{\theta}(E^{\infty},\R)$ with $P(\psi)<\infty$. Then there exists a unique triplet $(\lam,h,\nu)\in \R\times F_{\theta,b}(E^{\infty})\times C_{b}(E^{\infty})^{*}$ such that the following are satisfied:
\ite
{
\item The number $\lam$ is positive and a simple maximal eigenvalue of the operator $\LR_{\psi}\in \LR(F_{\theta,b}(E^{\infty}))$ and is equal to $\exp(P(\psi))$.
\item The operator $\LR_{\psi}\in \LR(F_{\theta,b}(E^{\infty}))$ has the decomposition
\ali
{
\LR_{\psi}=\lam \PR+\RR
}
with $\PR\RR=\RR\PR=O$. Here the operator $\PR$ is a projection onto the one-dimensional eigenspace of the eigenvalue $\lam$. Moreover, this has the form $\PR f=\int_{E^{\infty}}f h\,d\nu$ for $f\in C_{b}(E^{\infty})$, where $h\in F_{\theta,b}(E^{\infty},\R)$ is the corresponding eigenfunction of $\lam$ and $\nu$ is the corresponding eigenvector of $\lam$ of the dual $\LR_{\psi}^{*}$ with $\nu(h)=1$. Here $h$ is bounded uniformly away from zero and infinity, and $\nu$ is a Borel probability measure on $E^{\infty}$. In particular, $h\nu$ is the $\si$-invariant Gibbs measure of $\psi$.
\item The spectrum of $\RR\in \LR(F_{\theta,b}(E^{\infty}))$ is contained in $\{z\in \C\,:\,|z-\lam|\geq \rho\}$ for some $\rho>0$.
}
}
\pros
We will prove (1)(2). When the incidence matrix of $E^{\infty}$ is finitely primitive, these assertions (1)(2) follow from \cite{Sar99}. 
The existence $(\lam,h,\nu)$ can be extended to the finitely irreducible case. In fact, we decompose $E^{\infty}$ into $\Si_{0}, \Si_{1}, \dots, \Si_{p-1}$ such that each $(\Si_{i},\si^{p})$ is topologically mixing and $\si \Si_{i}\subset \Si_{(i+1) \mod p}$ for some integer $p\geq 1$. It can be also check that $(\Si_{i},\si^{p})$ satisfies the big images and pre-images property and therefore this is finitely primitive. We take $(\lam_{i},h_{i},\nu_{i})\in \R\times F_{\theta,b}(\Si_{i})\times M(\Si_{i})$ so that $\LR_{\psi}^{p}h_{i}=\lam_{i}^{p} h_{i}$, $(\LR_{\psi}^{*})^{p}\nu_{i}=\lam_{i}^{p} \nu_{i}$ and $\nu_{i}(h_{i})=\nu_{i}(1)=1$. When we put $\nu_{*}=\sum_{i=0}^{p-1}\lam_{0}^{-i}(\LR_{\psi}^{*})^{i}\nu_{0}$ and $h_{*}=\sum_{i=0}^{p-1}\lam_{0}^{i}\LR_{\psi}^{i}h_{0}$, the measure $\nu:=\nu_{*}/\nu_{*}(1)$ and the function $h:=h_{*}/\nu(h_{*})$ satisfy (1)(2). The boundedness of $h$ is yielded by the boundedness of $h_{i}$ for each $i$. The equality $\lam_{0}=\exp(P(\psi))$ follows from the equation $P(\psi)=\lim_{n\to \infty}(1/n)\log\|\LR_{\psi}^{n}1\|_{\infty}$ in addition to the finitely irreducibility and H\"older continuously of $\psi$. The equation $h\nu=\mu$ is guarantied from \cite[Corollary 2.7.5]{MU}.
Finally, the assertion $\PR\RR=\RR\PR=O$ directly follows from the definition of $\PR$. 
\smallskip
\par
In order to see $(3)$, it only has to be shown that the essential spectral radius of $\LR_{\psi}$ is less than the spectral radius of $\LR_{\psi}$. Indeed, the operator $\LR_{\psi}$ satisfies that there exist $k\geq 1$, $0<\rho<\lam$ and $c_{\adl{LYineq}}>0$ such that for any $f\in F_{\theta,b}(E^{\infty})$, $\|\LR_{\psi}^{k}f\|_{\theta}\leq \rho^{k}\nu(h |f|)+c_{\adr{LYineq}}\|f\|_{\theta}$ as the condition of Hennion's theorem (see \cite{AD,Sar09}).
\proe
\section{Asymptotic perturbation of eigenvectors of bounded linear operators}\label{sec:abstract_asymp}
In this section, we study asymptotic behaviour of the eigenvalues and eigenvectors of perturbed bounded linear operators under an abstract setting.

Put $\K=\R$ or $\K=\C$. Let $(\mathcal{X}_{0},\|\cdot\|_{0})$ be a normed space over $\K$ and $(\mathcal{X}_{1},\|\cdot\|_{1})$ a Banach space over $\K$ such that $\mathcal{X}_{1}\subset \mathcal{X}_{0}$ and $\|f\|_{0}\leq \|f\|_{1}$ for any $f\in \mathcal{X}_{1}$. 
We write $\mathcal{X}^{*}$ as the dual space of $\mathcal{X}$ and $\LR^{*}\in \LR(\XR^{*})$ as the dual operator of $\LR\in \LR(\XR)$.

Let $\mathcal{L}\in \mathcal{L}(\mathcal{X}_{0})\cap \mathcal{L}(\mathcal{X}_{1})$ and $\mathcal{L}(\epsilon,\cdot)\in \mathcal{L}(\mathcal{X}_{0})$. Take $(\lam,\nu), (\lam(\e),\nu(\e,\cd))\in \K\times \XR_{0}^{*}$ so that $\LR^{*}\nu=\lam \nu$ and $\LR(\e,\cd)^{*}\nu(\e,\cd)=\lam(\e)\nu(\e,\cd)$. We assume the following conditions:
\begin{itemize}
\item[(L.1)] There exists $h\in \mathcal{X}_{1}$ such that $\LR h=\lam h$ and $\nu(h)=1$.
\item[(L.2)] The operator $\LR$ has the decomposition $\LR=\lam \PR+\RR$ satisfying that (i) $\PR\in \LR(\XR_{0})\cap \LR(\XR_{1})$ has the form $\PR f=\nu(f)h$, (ii) $\PR\RR=\RR\PR=O$, and (iii) $\lambda$ is in the resolvent set of the operator $\RR\in \mathcal{L}(\mathcal{X}_{1})$.
\item[(L.3)] $\limsup_{\epsilon\to 0}\|\nu(\epsilon,\cdot)\|_{0}^{*}/\nu(\e,h)<\infty$, where $\|\nu(\epsilon,\cdot)\|_{0}^{*}=\sup_{f\in \mathcal{X}_{0}\,:\,\|f\|_{0}\leq 1}|\nu(\epsilon,f)|$.
\item[(L.4)] There exist operators $\mathcal{L}_{1},\dots, \mathcal{L}_{n}\in \mathcal{L}(\mathcal{X}_{0})\cap \mathcal{L}(\mathcal{X}_{1})$ and $\ti{\LR}_{n}(\e,\cd)\in \LR(\XR_{0})$ such that $\mathcal{L}(\epsilon,\cdot)=\mathcal{L}+\mathcal{L}_{1}\epsilon+\cdots+\mathcal{L}_{n}\epsilon^{n}+\tilde{\mathcal{L}}_{n}(\epsilon,\cdot)\epsilon^{n}$ and $\|\tilde{\mathcal{L}}_{n}(\epsilon,f)\|_{0}\to 0$ as $\epsilon\to 0$ for each $f\in \mathcal{X}_{1}$.
\end{itemize}
Let $\SR\in \LR(\XR_{1})$ be $\SR=(\RR-\lam \IR)^{-1}(\IR-\PR)$, where $\IR$ is the identity operator on $\XR_{1}$. 
Numbers $\lam_{k}\in \K$ and linear functionals $\kappa_{k}\,:\,\XR_{1}\to \K$ $(1\leq k\leq n)$ are defined by
\begin{align}
\lambda_{k}=&\sum_{j=1}^{k}\kappa_{k-j}(\mathcal{L}_{j}h),\label{eq:proof2}\\
\kappa_{k}(f)=&\sum_{j=1}^{k}\kappa_{k-j}((\lambda_{j}\mathcal{I}-\mathcal{L}_{j})\mathcal{S}f)\quad \text{ for each }f\in \mathcal{X}_{1}\label{eq:proof3}
\end{align}
inductively with $\kappa_{0}=\nu$ and $\LR_{0}=\LR$.
Then we have the following:
\begin{theorem}\label{th:asymp_e.vec}
Assume that the conditions (L.1)-(L.4) are satisfied for a fixed integer $n\geq 0$. We define $\kappa(\e,\cd)\in \XR_{0}^{*}$ by $\kappa(\e,f)=\nu(\e,f)/\nu(\e,h)$ for $f\in \XR_{0}$.
Then
\begin{itemize}
\item[(1)] 
$\lambda(\epsilon)=\lambda+\lambda_{1}\epsilon+\cdots+\lambda_{n}\epsilon^{n}+o(\epsilon^{n}) \text{ in }\K$;
\item[(2)] $\kappa(\epsilon,f)=\nu(f)+\kappa_{1}(f)\epsilon+\cdots+\kappa_{n}(f)\epsilon^{n}+o(\epsilon^{n}) \text{ in }\K \text{ for each }f\in \mathcal{X}_{1}$.
\end{itemize}
\end{theorem}
\cor
{\label{cor:asymp_e.vec}
In addition to the conditions (L.1)-(L.4), we assume that $\liminf_{\e\to 0}|\nu(\e,h)|>0$, and there exists $1_{\XR}\in \XR_{1}$ such that $\nu(\e,1_{\XR})=1$ for any $\e>0$. Then the eigenvector of $\lam(\e)$ has the expansion
\ali
{
\nu(\e,f)=\nu(f)/\nu(1_{\XR})+\nu_{1}(f)\e+\cdots+\nu_{n}(f)\e^{n}+\ti{\nu}_{n}(\e,f)\e^{n}
}
and $|\ti{\nu}_{n}(\e,f)|\to 0$ in $\K$ for each $f\in \XR_{1}$, where we put $\nu_{k}(f)=\sum_{0\leq i,j\leq k\,:\,i+j=k}b_{j}\kappa_{i}(f)$, $b_{0}=1/\nu(1_{\XR})$ and
\ali
{
b_{j}=\sum_{l=1}^{j}\frac{1}{\nu(1_{\XR})^{l+1}}\sum_{i_{1},\dots, i_{j}\geq 0\,:\,\atop{i_{1}+\cdots+i_{j}=l\atop{i_{1}+2i_{2}+\cdots+j\cd i_{j}=j}}}\frac{(-1)^{l}l!}{i_{1}!\cdots i_{j}!}\kappa_{1}(1_{\XR})^{i_{1}}\cdots \kappa_{l}(1_{\XR})^{i_{j}} \quad \text{ for }1\leq j\leq n.
}
}
\pros
When we put $f=1_{\XR}$ in Theorem \ref{th:asymp_e.vec}(2), we have $\nu(\e,h)^{-1}=\nu(1_{\XR})+\sum_{k=1}^{n}\kappa_{k}(1_{\XR})\e^{k}+o(\e^{n})$ and $\nu(1_{\XR})\neq 0$. The assertion follows from the asymptotic expansion of $\nu(\e,h)$ and the form $\nu(\e,f)=\kappa(\e,f)\nu(\e,h)$.
\proe
\rem
{
\item[(1)] The following results are generalizations of the results of \cite{T2011} which gave the asymptotic behavior of the maximal eigenvalue of the Ruelle operators with finite state and the corresponding eigenprojection. 
\item[(2)] When the reminder $\tLR_{n}(\e,\cd)$ satisfies $\|\tLR_{n}(\e,\cd)\|_{1}\to 0$, the above results are implied by the general asymptotic perturbation theory \cite{Kato}. If $\LR_{n}(\e,\cd)$ is a Ruelle operator with finite state and satisfies $\|\tLR_{n}(\e,\cd)\|_{0}\to 0$ with $\|\cd\|_{0}:=\|\cd\|_{\infty}$, then similar assertions followed from \cite{T2011}. Keller and Liverani in \cite{KL} considered convergence of eigenvalue and eigenprojection in an abstract setting under a uniform Lasota-Yorke type inequality such as $\|\LR(\e,f)^{n}\|_{1}\leq c\alpha^{n}\|f\|_{1}+c M^{n}\|f\|_{0}$ for any $\e>0$ and $f\in \XR_{1}$ for some constant $c>0$, $0<M\leq \sup_{\e>0}\|\LR(\e,\cd)\|_{0}$ and $0<\alpha<M$. Under such an inequality, Gou\"ezel and Liverani in \cite[Section 8]{GL} studied the asymptotic perturbation of bounded linear operators. We stress that our assertion does not need a uniform Lasota-Yorke type inequality.
}
\noindent
({\bf Proof of Theorem \ref{th:asymp_e.vec}}). 
We start with the equation
$(\mathcal{L}-\lambda \mathcal{I})\mathcal{S}=\mathcal{I}-\mathcal{P}$ on $\XR_{1}$. By the definition of the operator $\mathcal{P}$, this is projection, i.e. $\mathcal{P}^{2}=\mathcal{P}$. The equation follows from $(\LR-\lam\IR)(\IR-\PR)=\LR-\lam\IR$ and $(\mathcal{I}-\mathcal{P})(\RR-\lambda \mathcal{I})=\RR-\lambda\mathcal{I}+\lambda\mathcal{P}=\mathcal{L}-\lambda\mathcal{I}$.

We first prove the assertions (1)(2) in the case when $n=0$. Consider the equation
\begin{align}
(\lambda(\epsilon)-\lambda)\nu(\epsilon,h)=\nu(\epsilon,(\mathcal{L}(\epsilon,\cdot)-\mathcal{L})h)\label{eq:(la-l)PaP=Pa(La-L)P0}
\end{align}
by using $\LR(\e,\cd)^{*}\nu(\e,\cd)=\lam(\e)\nu(\e,\cd)$ and $\LR h=\lam h$. This yields
$|\lambda(\epsilon)-\lambda|\leq \|\kappa(\epsilon,\cdot)\|_{0}^{*}\|\ti{\LR}_{0}(\e,h)\|_{0}\to 0$
with the conditions (L.3) and (L.4). Therefore we have $\lambda(\epsilon)\to \lambda$. On the other hand, we obtain that for each $f\in \mathcal{X}_{1}$
\begin{align*}
|\kappa(\epsilon,(\mathcal{I}-\mathcal{P})f)|&=|\kappa(\epsilon,(\mathcal{L}-\lambda \mathcal{I})\mathcal{S}f)|\\
&=|\kappa(\epsilon,(\mathcal{L}-\mathcal{L}(\epsilon,\cdot)+(\lambda(\epsilon)-\lambda) \mathcal{I})\mathcal{S}f)|&\\
&\leq\|\kappa(\epsilon,\cdot)\|_{0}^{*}(\|(\mathcal{L}(\epsilon,\cdot)-\mathcal{L})\mathcal{S}f\|_{0}+|\lambda(\epsilon)-\lambda|\|\mathcal{S}f\|_{0})\to 0
\end{align*}
as $\epsilon\to 0$. This and the fact $\kappa(\e,h)\equiv 1$ imply $\kappa(\e,f)\to \nu(f)$ for $f\in \XR_{1}$.

Assume $n\geq 1$. 
To show the assertions (1)(2), we assume that the assertions (1)(2), (\ref{eq:proof2}) and (\ref{eq:proof3}) are valid for each $n^{\prime}=0,1,\dots,n-1$. We will check the case $n^{\prime}=n$. 
By (\ref{eq:proof2}) for each $n^{\prime}=1,2,\dots, n-1$, and the equation (\ref{eq:(la-l)PaP=Pa(La-L)P0}), we have the equation
\begin{align*}
&\frac{\lambda(\epsilon)-\lambda-\lambda_{1}\epsilon-\cdots-\lambda_{n-1}\epsilon^{n-1}}{\epsilon^{n}}
=\kappa(\e,\frac{\lambda(\epsilon)-\lambda-\lambda_{1}\epsilon-\cdots-\lambda_{n-1}\epsilon^{n-1}}{\epsilon^{n}}h)\nonumber\\
=&\kappa(\e,\frac{\mathcal{L}(\epsilon,\cdot)-\mathcal{L}-\sum_{l=1}^{n-1}\mathcal{L}_{l}\epsilon^{l}}{\epsilon^{n}}h)+\sum_{l=1}^{n-1}\frac{\kappa(\e,\mathcal{L}_{l}h)\epsilon^{l}-\sum_{j=1}^{l}\kappa_{l-j}(\mathcal{L}_{j}h)\epsilon^{j}}{\epsilon^{n}}\\
=&\kappa(\e,\frac{\mathcal{L}(\epsilon,\cdot)-\mathcal{L}-\sum_{l=1}^{n-1}\mathcal{L}_{l}\epsilon^{l}}{\epsilon^{n}}h)+\sum_{l=1}^{n-1}\frac{\kappa(\e,\mathcal{L}_{l}h)-\sum_{j=0}^{n-l-1}\kappa_{j}(\mathcal{L}_{l}h)\epsilon^{j}}{\epsilon^{n-l}}\nonumber\\
\to&\nu(\LR_{n}h)+\sum_{l=1}^{n-1}\kappa_{n-l}(\LR_{l}h)=\sum_{l=1}^{n}\kappa_{n-l}(\LR_{l}h)=:\lam_{n}.
\end{align*}
Thus, (1) and (\ref{eq:proof2}) are valid for $n$. Finally, we check (2) and (\ref{eq:proof3}). We obtain
\ali
{
&\frac{\kappa(\e,f)-\sum_{l=0}^{n-1}\kappa_{l}(f)\epsilon^{l}}{\epsilon^{n}}\nonumber
=\kappa(\e,\frac{\IR-\PR}{\epsilon^{n}}f)-\sum_{l=1}^{n-1}\frac{\kappa_{l}(f)\e^{l}}{\e^{n}}\nonumber\\
=&\kappa(\e,\frac{\mathcal{L}-\mathcal{L}(\epsilon,\cdot)+(\lambda(\epsilon)-\lambda)\mathcal{I}}{\epsilon^{n}}\mathcal{S}f)-\sum_{l=1}^{n-1}\frac{\sum_{i=1}^{l}\kappa_{l-i}((\lambda_{i}\mathcal{I}-\mathcal{L}_{i})\mathcal{S}f)}{\epsilon^{n-l}}\nonumber\\
=&-\kappa(\e,\frac{\mathcal{L}(\epsilon,\cdot)-\mathcal{L}-\sum_{l=1}^{n-1}\mathcal{L}_{l}\epsilon^{l}}{\epsilon^{n}}\mathcal{S}f)+\kappa(\e,\frac{\lambda(\epsilon)-\sum_{l=0}^{n-1}\lambda_{l}\epsilon^{l}}{\epsilon^{n}}\mathcal{S}f)\nonumber\\
&+\kappa(\e,\frac{\sum_{l=1}^{n-1}(\lambda_{l}\IR-\mathcal{L}_{l})\epsilon^{l}}{\epsilon^{n}}\mathcal{S}f)-\sum_{l=1}^{n-1}\frac{\sum_{i=1}^{l}\kappa_{l-i}((\lambda_{i}\mathcal{I}-\mathcal{L}_{i})\mathcal{S}f)}{\epsilon^{n-l}}\nonumber\\
=&-\kappa(\e,\frac{\mathcal{L}(\epsilon,\cdot)-\mathcal{L}-\sum_{l=1}^{n-1}\mathcal{L}_{l}\epsilon^{l}}{\epsilon^{n}}\mathcal{S}f)+\frac{\lambda(\epsilon)-\lambda-\sum_{l=1}^{n-1}\lambda_{l}\epsilon^{l}}{\epsilon^{n}}\kappa(\e,\mathcal{S}f)\nonumber\\
&+\sum_{l=1}^{n-1}\frac{\kappa(\e,\cd)-\nu(\cd)-\sum_{j=1}^{n-l-1}\kappa_{j}(\cd)\epsilon^{j}}{\epsilon^{n-l}}((\lambda_{l}\mathcal{I}-\mathcal{L}_{l})\mathcal{S}f)\nonumber\\
\to&-\nu(\LR_{n}\SR f)+\lam_{n}\nu(\SR f)+\sum_{l=1}^{n-1}\kappa_{n-l}((\lam_{l}\IR-\LR_{l})\SR f)=\sum_{l=1}^{n}\kappa_{n-l}((\lam_{l}\IR-\LR_{l})\SR f)=:\kappa_{n}(f).
}
Hence (2) and (\ref{eq:proof3}) are fulfilled for $n$.
\qed
\section{Estimate of intermediate point of the binomial expansion}\label{sec:interme_binom}
Taylor expansion implies that for any $n\geq 1$, $a>0$ and $s\in (0,1]$, the map $x\mapsto (a+x)^{s}$ has the form
\alil
{
(a+x)^{s}=a^s+\sum_{k=1}^{n-1}\binom{s}{k}a^{s-k}x^k+\binom{s}{n}(a+x\alpha)^{s-n}x^{n}\label{eq:(a+x)^s=}
}
for some constant $\alpha=\alpha(n,a,s,x)\in [0,1]$, where $\binom{s}{k}$ is the binomial coefficient $s(s-1)\dots (s-k+1)/k!$. In this section, we will estimate the lower bound of the intermediate point $\alpha$ which plays an important role in giving the asymptotic expansion of $\exp(t\phe)$ (see Proof of Theorem \ref{th:ex_5}). Note that the estimate of $\alpha$ was studied by \cite{Haber,Haber_Shisha}.
\prop
{\label{prop:lower_int_biom}
Assume that the map $x\mapsto (a+x)^{s}$ has the expansion (\ref{eq:(a+x)^s=}). Then there exist two positive continuous functions $(0,1)\ni s\mapsto L(n,s)$, $(0,1)\ni s\mapsto M(n,s)\in (0,1)$ such that for any $a,x>0$ with $0<a/x\leq M(n,s)$
\ite
{
\item[(1)] if $n=1$ then $\alpha(n,a,s,x)\geq L(n,s)$;
\item[(2)] if $n\geq 2$ then $\alpha(n,a,s,x)\geq (a/x)^{\frac{n-1-s}{n-s}}L(n,s)$.
}
}
\pros
(1) Assume $n=1$. We have
\ali
{
(a+x)^{s}&=a^{s}+s(a+\alpha x)^{s-1}x\\
(\frac{a}{x}+1)^{s}-\left(\frac{a}{x}\right)^{s}&=s(\frac{a}{x}+\alpha)^{s-1}.
}
If $x,a$ satisfy $a/x\leq 2^{-s/(1-s)-1}<1$, then
\ali
{
\alpha=&\frac{1}{\left(\left(\frac{a}{x}+1\right)^{s}-\left(\frac{a}{x}\right)^{s}\right)^{1/(1-s)}}-\frac{a}{x}\geq \frac{1}{2^{s/(1-s)}}-\frac{1}{2}\frac{1}{2^{s/(1-s)}}=\frac{1}{2\cd 2^{s/(1-s)}}=\frac{1}{2^{1/(1-s)}}.
}
Thus we obtain the assertion by putting $L(1,s)=2^{-1/(1-s)}$ and $M(1,s)=2^{-s/(1-s)-1}$.
\smallskip
\\
(2) Assume $n\geq 2$. We will solve the equation (\ref{eq:(a+x)^s=}) for $\alpha$. This equation implies
\ali
{
\left(\frac{a}{x}+1\right)^{s}=\sum_{k=0}^{n-1}\binom{s}{k}\left(\frac{a}{x}\right)^{s-k}+\binom{s}{n}\left(\frac{a}{x}+\alpha\right)^{s-n}.
}
Noting the fact $\sign \binom{s}{n}=(-1)^{n-1}$, we have
\ali
{
&\left|\binom{s}{n}\right|\left(\frac{a}{x}+\alpha\right)^{s-n}=(-1)^{n-1}\left(\frac{a}{x}+1\right)^{s}+(-1)^{n}\sum_{k=0}^{n-1}\binom{s}{k}\left(\frac{a}{x}\right)^{s-k}\\
=&\left(\frac{a}{x}\right)^{s-n+1}\left(b(s,a,x)+\left|\binom{s}{n-1}\right|\right)
}
by putting $b(s,a,x)=(-1)^{n-1}\left(\frac{a}{x}\right)^{-s+n-1}\left(\frac{a}{x}+1\right)^{s}+(-1)^{n}\sum_{k=0}^{n-2}\binom{s}{k}\left(\frac{a}{x}\right)^{n-1-k}$.
Therefore
\ali
{
\alpha=\left(\frac{a}{x}\right)^{\frac{n-s-1}{n-s}}\left\{\left(\frac{\left|\binom{s}{n}\right|}{\left|\binom{s}{n-1}\right|+b(s,a,x)}\right)^{1/(n-s)}-\left(\frac{a}{x}\right)^{\frac{1}{n-s}}\right\}.
}
When $a\leq x$, we see
\ali
{
|b(s,a,x)|\leq& \left(\frac{a}{x}\right)^{-s+n-1}2^{s}+\sum_{k=0}^{n-2}\left|\binom{s}{k}\right|\left(\frac{a}{x}\right)^{n-1-k}
\leq\left(\frac{a}{x}\right)^{1-s}c_{\adr{IB1}}
}
with $c_{\adl{IB1}}=c_{\adr{IB1}}(s,a,x)=2^{s}+\sum_{k=0}^{n-2}\left|\binom{s}{k}\right|$.
Consequently, for any $x,a>0$ satisfying that $a/x\leq M(n,s):=\min(1, (1/c_{\adr{IB1}})^{1/(1-s)}, (1/2)^{n-s}\left|\binom{s}{n}\right|/(\left|\binom{s}{n-1}\right|+1))$, we obtain that $|b(s,a,x)|\leq 1$ and
\ali
{
\alpha\geq \left(\frac{a}{x}\right)^{\frac{n-s-1}{n-x}}\frac{1}{2}\left(\frac{\left|\binom{s}{n}\right|}{\left|\binom{s}{n-1}\right|+1}\right)^{1/(n-s)}.
}
Hence the assertion is fulfilled by putting $L(n,s)=\frac{1}{2}\left(\left|\binom{s}{n}\right|/(\left|\binom{s}{n-1}\right|+1)\right)^{1/(n-s)}$.
\proe

\endthebibliography

\begin{thebibliography}{99}
\bibitem{AD} Aaronson, J., Denker, M.: Local limit theorems for partial sums of stationary sequences generated by Gibbs-Markov maps. Stoch. Dyn. \textbf{1}, no. 2, 193237 (2001).
\bibitem{Bowen} Bowen, R.: Equilibrium states and the ergodic theory of Anosov diffeomorphisms, Lecture Notes in Mathematics, 470 (Springer, Berlin, 1975).
\bibitem{BP} Benedetti, R., Petronio, C.: Lectures on Hyperbolic Geometry (Springer, Berlin, 1992).
\bibitem{BS} Buzzi, J., Sarig, O.: Uniqueness of equilibrium measures for countable Markov shifts and multidimensional piecewise expanding maps. Erg. Th. Dynam. Sys., \textbf{23}, 1383--1\adr{g3} (2003).
\bibitem{Constantine_Savits} Constantine, C., Savits, T. A multivariate Fa\`a di Bruno formula with applications, Trans. Amer. Math. Soc. \textbf{348}, 503--520 (1996).
\bibitem{DS} Dunford, N., Schwartz, J.T.: Linear Operators Part I. John Wiley \& Sons, New York (1988).
\bibitem{GL} Gou\"ezel, S., Liverani, C.: Banach spaces adapted to Anosov systems, Erg. Th. Dynam. Sys., \textbf{26} 1, 189--217 (2006).
\bibitem{Haber} Haber, S. An elementary inequality. Internat. J. Math. \& Math, Sci. {\bf 2}, 531--535 (1979).
\bibitem{Haber_Shisha} Haber, S. ; Shisha, O., On the location of the intermediate point in Taylor's Theorem. General Inequalities, 2 (Proc. Second Internat. Conf. Oberwolfach, 1978), 143--144, Birkhauser, Basel, (1980).
\bibitem{Kato} Kato, T.: Perturbation Theory for Linear Operators. Springer, Berlin (1995).
\bibitem{KL} Keller, G., Liverani, C.: Stability of the spectrum for transfer operators. Ann. Scuola Norm. Sup. Pisa Cl. Sci. (4) \textbf{28}, 141--152 (1999).
\bibitem{Kuczma3} Kuczma, M. An Introduction to the Theory of Functional Equations and Inequalities. Warszawa-Krakow-Katowice (1985).
\bibitem{MU1999} Mauldin R. D., Urba\'nski, M.: Conformal iterated function systems with applications to the geometry
of continued fractions. Trans. Amer. Math. Soc. \textbf{351}, 4995--5025 (1999).
\bibitem{MU} Mauldin, R. D., Urba\'nski, M.: Graph Directed Markov Systems : Geometry and dynamics of limit sets, Cambridge (2003).
\bibitem{NPL} R. D. Nussbaum, A. Priyadarshi, and S. Verduyn Lunel, Positive operators and Hausdorff dimension of invariant sets. Trans. Amer. Math. Soc. \textbf{364}, no. 2, 1029--1066 (2012).
\bibitem{PP} Parry, W., Pollicott, M.: Zeta functions and the periodic orbit structure of hyperbolic dynamics. Ast\'erisque, pp. 187--188  (1990).
\bibitem{Priyadarshi} Priyadarshi, A. Infinite Graph-Directed Systems and Hausdorff Dimension. Waves Wavelets Fractals Adv. Anal. 3:84--95 (2017).
\bibitem{RU} Roy, M., Urba\'nski, M.: 
Regularity properties of Hausdorff dimension in infinite conformal iterated function systems,
Ergod. Th. \& Dynam. Sys. \textbf{25}, 1961--1983 (2005).
\bibitem{RU2} Roy, M., Urba\'nski, M.: Real analyticity of Hausdorff dimension for higher dimensional graph directed Markov systems,
Math. Z. \textbf{260}, 153--175 (2008).
\bibitem{Ruelle} Ruelle, D.:
Thermodynamic Formalism, 2nd edn (Cambridge University Press, Cambridge, 2004).
\bibitem{Sar99} Sarig, Omri M.: Thermodynamic formalism for countable Markov shifts, Erg. Th. Dynam. Sys. \textbf{19}, no. 6, 1565--1593 (1999).
\bibitem{Sar03} Sarig, Omri M.: Existence of Gibbs measures for countable Markov shifts, Proc. Amer. Math. Soc. \textbf{131}, no. 6, 1751--1758 (2003).
\bibitem{Sar09} Sarig, Omri M.: Lecture Notes on Thermodynamic Formalism for Topological Markov Shifts, Penn State (2009).
\bibitem{SU} Stratmann, Bernd O.; Urba\'nski, M.: Pseudo-Markov systems and infinitely generated Schottky groups, Amer. J. Math. \textbf{129}, no. 4, 1019-1062 (2007).
\bibitem{T2011} Tanaka, H.: An asymptotic analysis in thermodynamic formalism. Monatsh. Math. \textbf{164}, 467--486 (2011).
\bibitem{T2016} Tanaka, H.: Asymptotic perturbation of graph iterated function systems, Journal of Fractal Geometry, \textbf{3}, 119--161 (2016).
\end{thebibliography}
\end{document}